\documentclass[a4paper,11pt]{article}

\usepackage[title]{appendix}%
\usepackage{textcomp}%
\usepackage{manyfoot}%
\usepackage{listings}%
\usepackage{tcolorbox}
\usepackage{amsfonts,amscd,amssymb,amsmath,amsthm,enumerate}
\usepackage{array,mathtools,mathrsfs,bm}
\usepackage{booktabs,caption,longtable,multicol,multirow}
\usepackage{subcaption,lscape}
\usepackage[linesnumbered,boxruled,lined,commentsnumbered]{algorithm2e}
\usepackage[colorlinks,citecolor=purple,linkcolor=blue]{hyperref}
\usepackage[parfill]{parskip}
\usepackage[margin=1.0in]{geometry}
\usepackage[numbers]{natbib}
\usepackage{authblk}

\newtheorem{theorem}{Theorem}
\newtheorem{proposition}{Proposition}
\newtheorem{corollary}{Corollary}
\newtheorem{lemma}{Lemma}
\newtheorem{assumption}{Assumption}
\newtheorem{definition}{Definition}
\newtheorem{remark}{Remark}
\newtheorem{example}{Example}

\numberwithin{assumption}{section}
\numberwithin{definition}{section}
\numberwithin{equation}{section}
\numberwithin{lemma}{section}
\numberwithin{corollary}{section}
\numberwithin{theorem}{section}
\numberwithin{example}{section}
\numberwithin{condition}{section}

\usepackage{tikz}
\definecolor{ao(english)}{rgb}{0.0, 0.5, 0.0}
\definecolor{cadmiumgreen}{rgb}{0.0, 0.42, 0.24}
\definecolor{darkpastelgreen}{rgb}{0.01, 0.75, 0.24}

\usepackage{adjustbox}
\tikzset{%
    materia/.style={draw, fill=blue!20, text width=6.0em, text centered, minimum height=1.5em,drop shadow},
    etape/.style={materia, text width=8em, minimum width=10em, minimum height=3em, rounded corners, drop shadow},
    texto/.style={above, text width=6em, text centered},
    linepart/.style={draw, thick, color=black!50, -LaTeX, dashed},
    line/.style={draw, thick, color=black!50, -LaTeX},
    ur/.style={draw, text centered, minimum height=0.01em},
    back group/.style={fill=white!20,rounded corners, draw=black!50, dashed, inner xsep=15pt, inner ysep=10pt},
}

\newcommand{\real}{\mathbb R}

\newcommand{\xk}{x_k}
\newcommand{\Hk}{H_k}

\newcommand{\gk}{g_k}

\newcommand{\dk}{d_k}
\newcommand{\gkn}{g_{k+1}}
\newcommand{\xkn}{x_{k+1}}

\newcommand{\diff}[1]{\frac{\mathsf{d}}{\mathsf{d}#1}}
\newcommand{\xmu}{x_\mu}
\newcommand{\tr}{\textrm{tr}}

\newcommand{\betacon}{concordant Lipschitz}

\newcommand{\ahsodm}{\texttt{Adaptive-HSODM}}
\newcommand{\pfhsodm}{\texttt{Homotopy-HSODM}}
\newcommand{\inewton}{\texttt{iNewton-Grad}}
\newcommand{\newtoncg}{\texttt{CG}}
\newcommand{\gmres}{\texttt{GMRES}}
\newcommand{\rgmres}{\texttt{rGMRES}}
\newcommand{\glanczos}{\texttt{Lanczos}}

\newcommand{\newtontrst}{\texttt{Newton-TR-STCG}}

\newcommand{\arc}{\texttt{Cubic-Reg}}

\author[1]{\small Chang He}
\author[1]{\small Yuntian Jiang}
\author[1]{\small Chuwen Zhang\thanks{Corresponding author. This research is partially supported by the National Natural Science Foundation of China (NSFC) [Grant NSFC-72150001, 72225009, 11831002] and the Natural Science Foundation of Shanghai [23ZR1445900].}}
\author[2]{\small Dongdong Ge}
\author[1]{\small Bo Jiang}
\author[2,3]{\small Yinyu Ye}

\affil[1]{\footnotesize School of Information Management and Engineering\\ Shanghai University of Finance and Economics}
\affil[2]{\footnotesize Antai College of Economics and Management, Shanghai Jiao Tong University}
\affil[3]{\footnotesize Department of Management Science and Engineering, Stanford University}

\title{Homogeneous Second-Order Descent Framework: A Fast Alternative to Newton-Type Methods}

\begin{document}

\maketitle
\abstract{
    This paper proposes a homogeneous second-order descent framework (HSODF) for nonconvex and convex optimization based on the generalized homogeneous model (GHM).
    In comparison to the Newton steps, the GHM can be solved by extremal symmetric eigenvalue procedures and thus grant an advantage in ill-conditioned problems. Moreover, GHM extends the ordinary homogeneous model (OHM) \cite{zhang_homogenous_2022} to allow adaptiveness in the construction of the aggregated matrix. Consequently, HSODF is able to recover some well-known second-order methods, such as trust-region methods and gradient regularized methods, while maintaining comparable iteration complexity bounds. We also study two specific realizations of HSODF. One is adaptive HSODM, which has a parameter-free $O(\epsilon^{-3/2})$ global complexity bound for nonconvex second-order Lipschitz continuous objective functions. The other is homotopy HSODM, which is proven to have a global linear rate of convergence without strong convexity. The efficiency of our approach to high-dimensional and ill-conditioned problems is justified by some preliminary numerical results.
}

\clearpage
\tableofcontents
\clearpage

\section{Introduction}

For the unconstrained smooth optimization problem $ \min_x f(x) $, Newton's method is one of the most powerful techniques, particularly due to its local quadratic convergence rate \cite{nocedal_numerical_2006} and its significant role in interior point methods \cite{nesterov_interior-point_1994, ye_interior_1997} and high-order tensor methods \cite{nesterovSuperfastSecondOrderMethods2021,cartisCubicquarticRegularizationModels2025,cartisSecondorderMethodsQuarticallyregularised2025}.

{
For nonconvex functions, \citet{nesterov_cubic_2006} showed the first $O(\epsilon^{-3/2})$ complexity bound of second-order methods (SOM) by cubic regularization (CR). In its original form, the authors focus on theoretical analyses and use an explicit parameter of the Lipschitz constants that may not be available in practice.
By utilizing the basic properties, the benefits of CR also extend not only to standard convex and strongly convex problems but also to star-convex and gradient-dominated functions. The latter is frequently observed in policy optimization; see \cite{masiha_stochastic_2022} for example.
After that, a group of new \emph{adaptive} second-order methods with $O(\epsilon^{-3/2})$ complexity bound emerged by establishing very similar properties, including adaptive cubic regularization (ARC) \cite{cartis_adaptive_2011,cartis_adaptive_2011-1}, modified trust-region methods \cite{curtis_trust_2017,curtis_inexact_2019}, line-search algorithms \cite{royer_complexity_2018}, to name a few\footnote{For a complete monograph on these methods, we refer the readers to \cite{cartis_evaluation_2022}}.
}

A very recent line of research focuses on second-order methods with various regularizations other than the \emph{cubic} one. For example, in the classical Levenberg-Marquardt method, one adopts the following iterate using a sequence of regularization parameters $\{\lambda_k\}_{k=0}^\infty$,
\begin{equation}\label{eq.lm}
    \xkn = \xk -\left(\Hk + \lambda_k I\right)^{-1} \gk,
\end{equation}
where $g_k := \nabla f(x_k)$, and $H_k := \nabla^2 f(x_k)$ are the gradient and Hessian at $x_k$. It is only recently that \citet{mishchenko_regularized_2023} showed the global $O\left(\frac{1}{k^2}\right)$ convergence rate for convex optimization, which was previously known for a cubic regularized method. Later extensions \cite{doikov_gradient_2021, doikov_super-universal_2022} make use of a similar idea based on the gradient norm. These methods are very efficient in convex optimization.

Inspired by the classic trick in quadratic programming \cite{yeApproximatingQuadraticProgramming1999a,sturmConesNonnegativeQuadratic2003,he_quaternion_2022}, \citet{zhang_homogenous_2022} proposed a homogeneous second-order descent method (HSODM). In their method, the iterates are constructed by solving the \emph{Ordinary Homogeneous Model} (OHM) with the aggregated matrix $F_k$:
\begin{equation}\label{eq.homoquadmodel}
    \min_{\|[v; t]\| \le 1} \psi_k(v,t;F_k)  := ~ [v;t]^TF_k[v;t],~ F_k              :=   \begin{bmatrix}
        \Hk   & \gk    \\
        \gk^T & \delta
    \end{bmatrix},
\end{equation}
where variables $v \in \mathbb{R}^n,t \in \mathbb{R}$.
    {
        The method sets a \emph{pre-fixed} $\delta \le 0$ such that the vector $[v;t]$ corresponds to the leftmost eigenvector of the aggregated matrix $F_k$. Then the HSODM undergoes a simple strategy to find a step size $\eta_k$, e.g., by a line-search method, and updates the iterate as $\xkn = \xk + \eta_k (v_k/t_k)$ using the direction $[v_k;t_k]$ generated by OHM. The method is shown to have an $O(\epsilon^{-3/2})$ iteration complexity for nonconvex problems. The numerical performance reported in \cite{zhang_homogenous_2022} showed an advantage over standard second-order methods.
    }
In this paper, we extend the idea of homogenization to a broader context, not limited to nonconvex problems with second-order Lipschitz continuity. We introduce the following \textit{Generalized} Homogeneous Model (GHM)
\begin{equation}\label{eq.ghqm}
    \min_{\|[v; t]\| \le 1} \psi_k(v,t;F_k)  := ~ [v;t]^TF_k[v;t] \quad \mbox{with} \quad   F_k   :=
    \begin{bmatrix}
        \Hk           & \phi_k(\xk) \\
        \phi_k(\xk)^T & \delta_k
    \end{bmatrix},
\end{equation}
albeit now $\delta_k \in \real$ is allowed for some adaptiveness.
Furthermore, we introduce the transformation $\phi_k: \real^n\mapsto\real^n$ in place of the gradient $\gk$ while keeping the advantage of using cheap symmetric eigenvalue procedures.
We show that this flexibility facilitates a machinery to realize other second-order methods and, more importantly, a general homogeneous framework in which new algorithms can be designed.

\paragraph{Connections to prior works}
Historically, the idea of homogenization to solve quadratic constrained quadratic programming problems (QCQP) with iterative schemes dated back to a few earlier works \cite{rojas_new_2001,adachi_solving_2017,adachiEigenvaluebasedAlgorithmAnalysis2019}.  The benefits borne out by eigenvalue procedures were observed
in numerical tests on trust-region subproblems \cite{adachi_solving_2017, rojas_new_2001}. Later literature following such an idea focused on solving the cubic regularization subproblem (CRS). \citet{lieder2020solving} showed that the CRS could be solved globally via a generalized eigenvalue problem. \citet{jia2022solving} independently proved that the CRS could be transformed into a quadratic eigenvalue problem (QEP), which could be further linearized to a generalized eigenvalue problem.
However, the dimension of their matrices typically increases to $2\cdot n + 2$, bringing a heavier computational burden than the GHM of dimension $n+1$.

\paragraph{Contributions}
{
    The contributions of this paper can be summarized as follows.
    We propose a homogeneous second-order descent framework (HSODF) in \autoref{alg.gbasic}, where the linear systems required in subproblems are replaced by the GHMs \eqref{eq.ghqm} that are extremal eigenvalue problems.
    We demonstrate its computational benefits via a discussion on the condition numbers and a few experiments on highly degenerate problems.

    A discussion is given to produce a family of second-order methods using HSODF at the price of conducting a search over $\delta_k$. In this venue, we particularly study an {\it adaptive} HSODM (\autoref{sec.ncvx}) for second-order Lipschitz functions, which is a strengthened version of \cite{zhang_homogenous_2022} that uses line-search methods.  Such a method maintains an $O(\epsilon^{-3/2})$ iteration complexity for nonconvex optimization while transcending the original HSODM to convex problems.  In the appendix, we illustrate an auxiliary $O(\log(1/\epsilon))$ bisection method for locating $\delta_k$, serving as a foundational technique that can be adapted for similar methods.

    Under the same framework, we next discuss a variant that works for \betacon{} instead of standard second-order Lipschitz functions (\autoref{sec.homotopy}).
    In this case, simultaneous adaptations on both  $\delta_k$ and $\phi_k$ are applied to GHMs. We show that the corresponding homotopy HSODM exhibits a global linear rate of convergence through a non-interior path-following technique. For this variant, it is worth mentioning that there is no need for the auxiliary procedure to locate $\delta_k$; at most $2$ GHMs are needed per iteration in sharp comparison to aforementioned ones.
}

\paragraph{Notations}
We introduce the notations used throughout the paper. Denote the standard Euclidean norm in space $\mathbb{R}^n$ by \(\|\cdot\|\). For a matrix \(A \in \mathbb{R}^{n\times n}\), \(\|A\|\) represents the induced \(\ell_2\) norm. Let $A^\star$ denote the pseudo-inverse of the matrix $A$. We let $P_{\mathcal X}$ be the orthogonal projection operator onto a space, where $\mathcal X \subseteq \real^n$. We use \textrm{mod} to denote the binary modulo operation. We say a vector $y$ is orthogonal to a subspace $\mathcal{S}$, i.e. $y \perp \mathcal{S}$ if for any nonzero vector $u \in \mathcal{S}$, $u^T y = 0$.

Next, we introduce the following notations for eigenvalues of Hessian $\Hk$. At an iterate of the algorithm $\xk$, we assume $\Hk$ has $r$ ($1 \le r \le n$) distinct eigenvalues $\{\lambda_1(\Hk), ..., \lambda_r(\Hk)\}$ where $\lambda_1(\Hk) < ... < \lambda_r(\Hk)$ and $\mathcal S_1(\Hk), ..., \mathcal S_r(\Hk)$ are subspaces spanned by corresponding eigenvectors. We sometimes use $\lambda_{\min}, \lambda_{\max}$ as synonyms for $\lambda_1$ and $\lambda_r$, respectively.
We denote the condition number of $\Hk$ as $\kappa(\Hk) = \frac{\lambda_r(\Hk)}{\lambda_1(\Hk)}$.
Since the discussion on eigenvalues is mostly restricted at the iterate $\xk$ only, we sometimes drop the index $k$ for simplicity.

\section{The Homogeneous Second-Order Descent Framework}\label{sec.motivate}

{
    We first give the optimality condition of \eqref{eq.ghqm} in the following. We omit the proof since it follows from the global optimality condition for the standard trust-region subproblem, see \cite{more_computing_1983,sorensen_newtons_1982}.
    \begin{lemma}[Optimality condition]\label{lemma.optimal condition of subproblem}
        $[v_k; t_k]$ is the optimal solution of the subproblem \eqref{eq.ghqm} if and only if there exists a dual variable $\theta_k \ge 0$ such that
        \begin{align}
            \label{eq.homoeig.soc}
                                        & \begin{bmatrix}
                                              \Hk + \theta_k \cdot I & \phi_k            \\
                                              \phi_k ^T              & \delta_k+\theta_k
                                          \end{bmatrix} \succeq 0, \\
            \label{eq.homoeig.foc}
                                        & \begin{bmatrix}
                                              \Hk + \theta_k \cdot I & \phi_k            \\
                                              \phi_k ^T              & \delta_k+\theta_k
                                          \end{bmatrix}
            \begin{bmatrix}
                v_k \\ t_k
            \end{bmatrix} = 0,                                                       \\
            \label{eq.homoeig.norm one} & \theta_k\cdot(\|[v_k; t_k]\| - 1) = 0.
        \end{align}
    \end{lemma}

    The optimal solutions in GHM slightly differ from that of OHM in \citet{zhang_homogenous_2022} in several ways. As adaptiveness of $\delta_k$ is introduced, $[v_k; t_k]$ may not reach the boundary of the unit ball and thus $\theta_k = 0$ indicated by \eqref{eq.homoeig.norm one}. Certainly, this is possible for a convex function and some $\delta_k \gg 0$. To stick to an eigenvalue procedure, we should prevent $\delta_k$ from being too large especially if the iterate is close to a local minimum.

    Another difficult case is when $t_k = 0$ such that a direction cannot be normalized. This dilemma corresponds to the \emph{hard case} of the trust-region subproblem \cite{rojas_new_2001,nocedal_numerical_2006}. However in OHM, this problem can be easily tackled using a fixed-radius strategy or introducing the truncation $\nu$ to update $d_k$ whenever $|t_k| < \nu$. Therefore, no extra operation is needed for the original HSODM, while some nontrivial analysis is required for GHMs as shown in \autoref{sec.hard case}.

    Next, we present the homogeneous second-order descent framework (HSODF) in \autoref{alg.gbasic} by using the GHM \eqref{eq.ghqm} as a subroutine.
}

\begin{figure}[ht]

    \begin{minipage}[t]{0.99\linewidth}%
        \centering
        \small
        \begin{algorithm}[H]
            \caption{Homogeneous Second-Order Descent Framework (HSODF)}\label{alg.gbasic}
            Given initial point $x_1$, controls $\delta_1$, maximum number of iteration $K_{\max}$\;
            \ForAll{$k = 1, \dots, K_{\max}$}{
            \ForAll{$j = 1, \dots, \mathcal T_k$}{
            \label{gbasic.sol}
            Construct a \textbf{GHM}: $F_{k,j} = \begin{bmatrix}
                    { H_{k}}              & \phi_{k,j}(x_{k,j}) \\
                    \phi_{k,j}(x_{k,j})^T & \delta_{k,j}
                \end{bmatrix}$\;
            Obtain $[v_{k,j}; t_{k,j}] = \arg\min\limits_{\|[v; t]\| \le 1}\psi_{k,j}(v,t;F_{k,j})$ (cf. \eqref{eq.ghqm})\;
            Set $d_{k,j} := v_{k,j} / t_{k,j}$\;
            {
            \label{gbasic.eval}  \If{$d_{k,j}$ satisfies (inner) termination criteria}{Set $\dk:=d_{k,j}$;
            \textbf{Break};}
            }
            \label{gbasic.adj}    Adjust $\delta_{k,j}$ and $\phi_{k,j}$.
            }
            Update $\xkn = \xk + \dk$\;
            {
            \label{gbasic.term}    \If{$\xkn$ satisfies (outer) termination criteria}{
                \KwOut{$\xk$}
            }
            }
            }
        \end{algorithm}
    \end{minipage}
    \normalsize
\end{figure}
{
The adaptive strategy in the HSODF excludes the need for Lipschitz constants in the algorithmic parameters. Unlike the original HSODM in \cite{zhang_homogenous_2022}, where $\delta_k \equiv \delta < 0$ is fixed and $ \phi_k = \gk$, each iteration $k$ involves an inner loop labeled by $j$ to search for suitable $\delta_{k,j}$ or $\phi_{k,j}$, terminating when $d_{k,j}$ satisfies certain conditions (cf. \autoref{gbasic.eval}). By properly designing such conditions, we can recover or provide alternatives to known second-order methods with comparable complexity.  Generally, for each $k$ the size of inner iteration $\mathcal T_k$ is around $O(\log(1/\epsilon))$; but in the homotopy HSODM, only 2 GHMs are needed for \betacon{} functions.
}
\subsection{Motivation for the HSODF}\label{sec.warmup}
{
    We motivate the design of HSODF (\autoref{alg.gbasic}) when dealing with high-dimensional data where the problem has nice sparsity and low-rank structures.
    To be specific, our computational findings indicate that solving a GHM can require fewer Krylov iterations than a Newton equation, especially when the Hessian $\Hk$ is degenerate.
    These findings are partly due to the difference between the conditioning of an eigenvalue problem and a linear equation for a Krylov subspace method. Based on this, we conduct further comparisons on the per-iteration cost between HSODF and Newton-type methods.
}


\subsubsection{Comparison of the per-iteration cost between HSODF and Newton-type methods}
To start with, we compare the computational cost of a Newton-type equation and the corresponding homogeneous model, which are solved in each iteration of Newton-type methods and HSODF respectively. Particularly, we consider the case that $\Hk$ is positive definite but ill-conditioned.
{Suppose the following perturbed Newton-type equation is solved at some iterate $\xk \in \real^n$,
\begin{equation}\label{eq.mo.newton}
    (\Hk + \epsilon_{\mathrm{N}} I) \dk = - \gk.
\end{equation}
Without loss of generality, we assume $\epsilon_{\mathrm{N}} \in [0,1)$}. When $n$ is large, the computation relies on iterative methods, such as the conjugate gradient method (\newtoncg{}), generalized minimum residual method (\gmres{}), and the restarted generalized minimum residual method (\rgmres{}) \cite{golub_matrix_2013}.
{
Suppose the HSODM in this case uses the GHM with $\delta_k = -\epsilon_{\mathrm{L}}$,
\begin{equation}\label{eq.mo.ghm}
    F_k := \begin{bmatrix}
        \Hk   & \gk                    \\
        \gk^T & -\epsilon_{\mathrm{L}}
    \end{bmatrix},
\end{equation}
}
the workhorses of which are those for symmetric eigenvalue problems, such as the Lanczos method (\glanczos{}) \cite{saad_numerical_2011}.
It is well known that the complexity of solving the linear system depends on the condition number of the perturbed matrix $H_k + \epsilon_{\mathrm{N}} I$. The number of iterations required by a conjugate gradient method is typically
\begin{equation}
    O(\sqrt{\kappa(\Hk + \epsilon_{\mathrm{N}} I)}\log(1/\varepsilon))\quad\mbox{with}\quad \kappa(\Hk + \epsilon_{\mathrm{N}} I) = \frac{\lambda_{\max}(\Hk) + \epsilon_{\mathrm{N}}}{\lambda_1(\Hk) + \epsilon_{\mathrm{N}}},
\end{equation}
where $\kappa(\cdot)$ is usually referred to as the condition number.

\begin{figure}
    \centering
    \begin{subcaptionblock}[c]{0.445\textwidth}
        \includegraphics[width=0.95\textwidth]{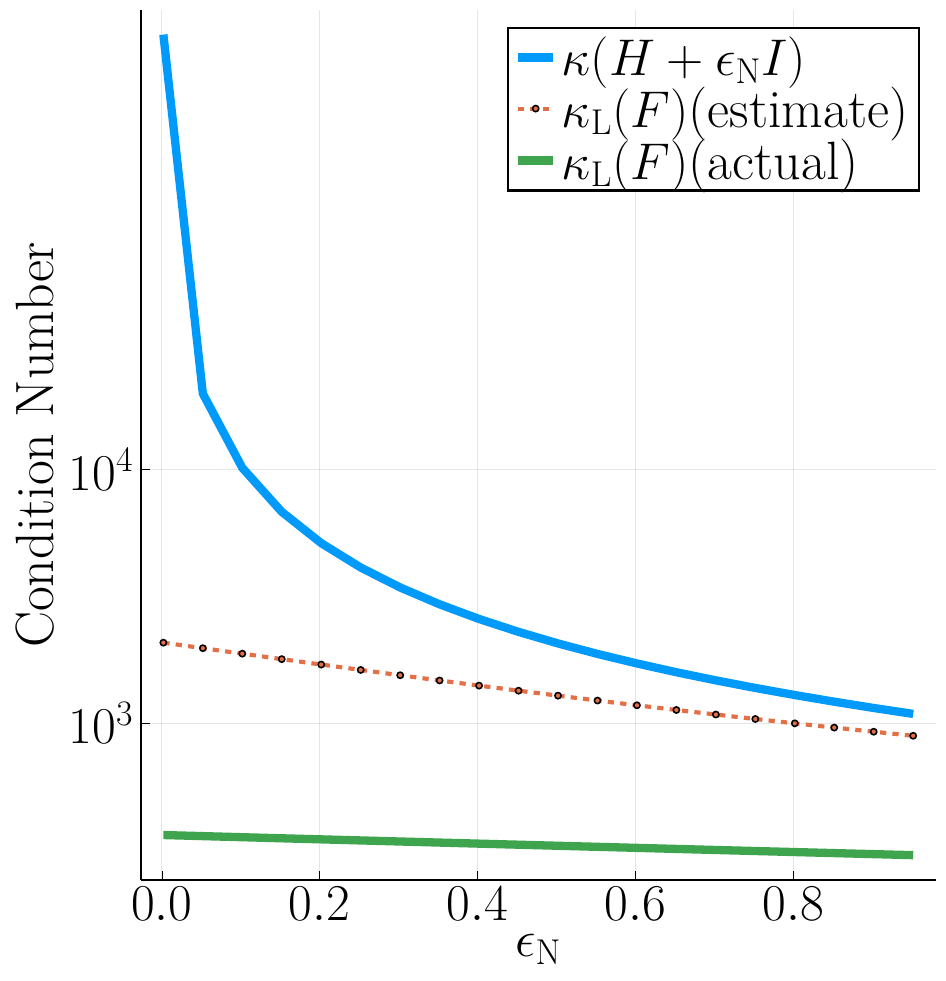}
        \caption{$\gk \not\perp \mathcal S_1$}
    \end{subcaptionblock}
    \begin{subcaptionblock}[c]{0.420\textwidth}
        \includegraphics[width=0.95\textwidth]{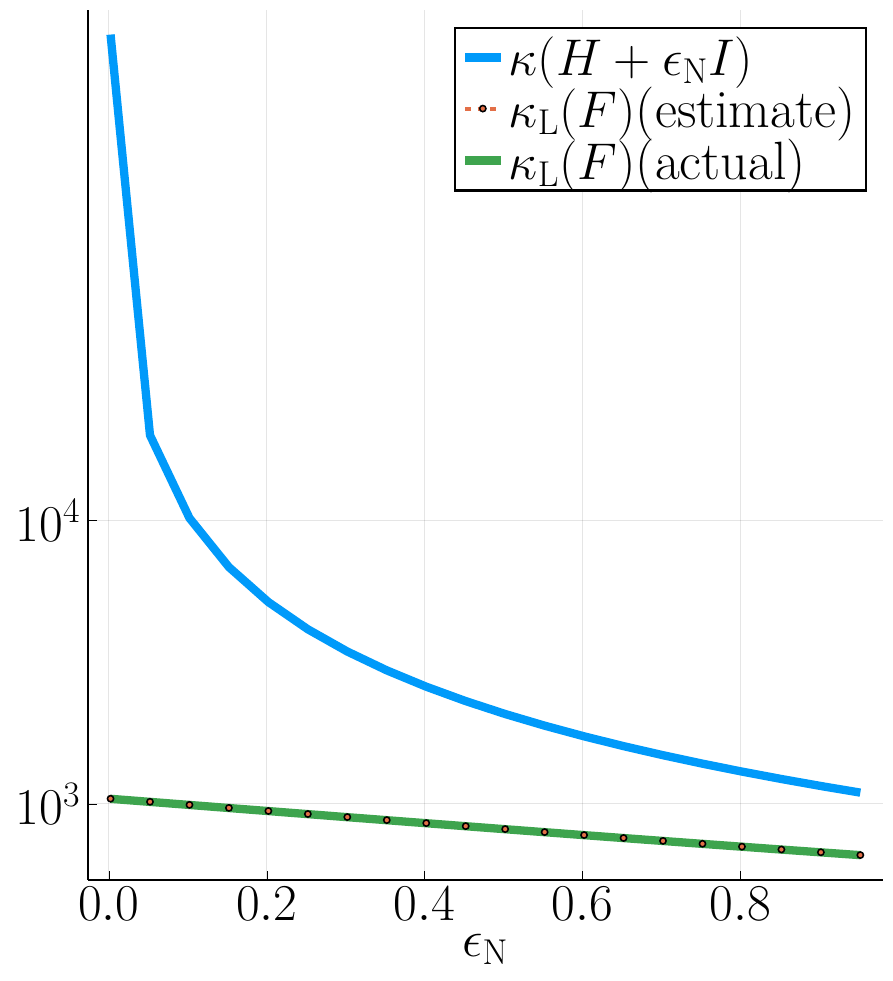}
        \caption{$\gk \perp \mathcal S_1$}
    \end{subcaptionblock}
    \caption{A comparison on $\kappa (\Hk + \epsilon_{\mathrm{N}} I)$ and $\kappa_{\mathrm{L}}(F_k)$ with degeneracy.
    }\label{fig.condno}
\end{figure}

In comparison, for finding the smallest eigenvalue and its associated eigenvector, the iteration complexity bound of the Lanzos method \cite{kuczynski_estimating_1992} is (also see \autoref{lem.eigenvector computation})
\begin{equation}
    O\left(\sqrt{\kappa_{\mathrm{L}}(F_k)}\log\left(1/\varepsilon\right)\right)\quad \mbox{with} \quad \kappa_{\mathrm{L}}(F_k) := \frac{\lambda_{\max}(F_k) - \lambda_1(F_k)}{\lambda_2(F_k) - \lambda_1(F_k)}
\end{equation}
in high probability,
where we accept certain estimate $\xi$ of the smallest eigenvalue such that $\xi - \lambda_1(F_k)\le \varepsilon$.
The following result indicates that $\kappa_{\mathrm{L}}(F_k)$ is always bounded from above, which is in contrast to the unboundness of $\kappa(\Hk + \epsilon_{\mathrm{N}} I)$ when $\epsilon_{\mathrm{N}}$ approaches $0$. Moreover, $\kappa_{\mathrm{L}}(F_k)$ can be much less than $\kappa(\Hk + \epsilon_{\mathrm{N}} I)$ in the degenerate case.
\begin{theorem}\label{thm.eigengap.better}
    For the aggregated matrix $F_k$ in \eqref{eq.mo.ghm}, suppose $U_H := \lambda_{\max}(\Hk) \gg \epsilon_{\mathrm{N}}$, it holds that
    \begin{enumerate}[$(a)$]
        \item For any $\epsilon_{\mathrm{L}} > 0$, $\kappa_{\mathrm{L}}$ is finite, specifically,
              \begin{equation}\label{eq.eigengap.better}
                  \kappa_{\mathrm{L}}(F_k) \le \frac{2(\lambda_{\max}(\Hk) - \epsilon_{\mathrm{L}} - \lambda_1(F_k))}{- U_H + \epsilon_{\mathrm{L}} + \sqrt{(U_H + \epsilon_{\mathrm{L}})^2 + \|\gk\|^2/n}} < \infty.
              \end{equation}
        \item {Furthermore, suppose $\lambda_1(\Hk) = 0$, then
              \begin{equation}
                  \frac{\kappa_{\mathrm{L}}(F_k)}{\kappa(\Hk+\epsilon_{\mathrm{N}} I)}
                  \le O\left(\frac{\epsilon_{\mathrm{N}} }{\frac{\|g_k\|^2}{U_H+\epsilon_{\mathrm{L}}} + \epsilon_{\mathrm{L}}} \right).
              \end{equation}
              }
    \end{enumerate}
\end{theorem}
We postpone the proof of the above theorem to \autoref{proof.thm.eigengap.better} for succinctness. The ratio in the second part of the theorem compares the
two condition numbers in solving the Newton-type equation and GHM respectively. The smaller value of the ratio implies the better condition number in GHM over that in the Newton-type equation.
{
This ratio shows that $\kappa_{\mathrm{L}}$ is generally better and more robust than that of the Newton-type equation, due to the implicit scaling of $\|\gk\|$.
For example, letting $\epsilon_{\mathrm{N}} = \epsilon_{\mathrm{L}} \rightarrow 0$, then the numerator approaches 0 while the denominator remains of constant order, demonstrating that $\kappa_{\mathrm{L}}$ is much smaller in this scenario. The opposite extreme indicates that the two condition numbers become close; the same holds if $\|\gk\| \to 0$.
Furthermore, this analysis echos the gradient-regularized steps where the perturbation is set to $\epsilon_{\mathrm{N}} = \|\gk\|^{1/2}$ \cite{mishchenko_regularized_2023,doikov_gradient_2021}.
}
Numerically, the above facts are visualized in \autoref{fig.condno} where the estimated $\kappa_{\mathrm{L}}$ are computed from \autoref{lem.eigengap.free}. Note that when $\gk \perp \mathcal S_1$, the estimate is tight.

\subsubsection{Further numerics of the per-iteration cost}
Motivated by the analysis above, we conduct experiments to compare the Krylov iterations utilizing the iterative methods. We use \newtoncg{}, \gmres{}, \rgmres{} to denote their corresponding results on \eqref{eq.mo.newton}, respectively. We use the Lanczos method in \eqref{eq.mo.ghm}. For solving linear equations \eqref{eq.mo.newton}, we set the residual at $10^{-5}$, and the Lanczos method in \eqref{eq.mo.ghm} is terminated at a fixed tolerance $10^{-7}$.

\paragraph{Hilbert matrix and a random right-hand-side}

\begin{figure}[h!]
    \centering
    \begin{subcaptionblock}[c]{0.75\textwidth}
        \includegraphics[width=0.95\textwidth]{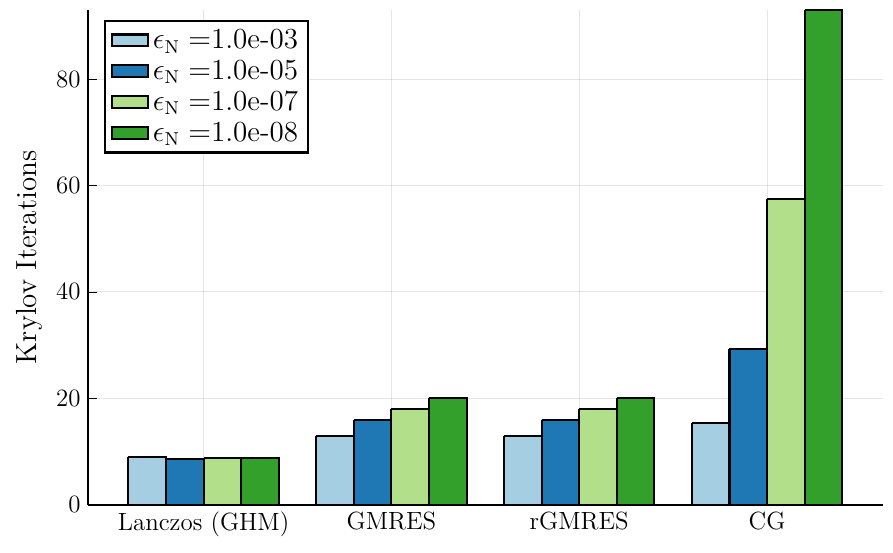}
        \caption{$\|\gk\| = 5.17\times 10^{-1}$}\label{fig.hilbert.large}
    \end{subcaptionblock}
    \begin{subcaptionblock}[c]{0.75\textwidth}
        \includegraphics[width=0.95\textwidth]{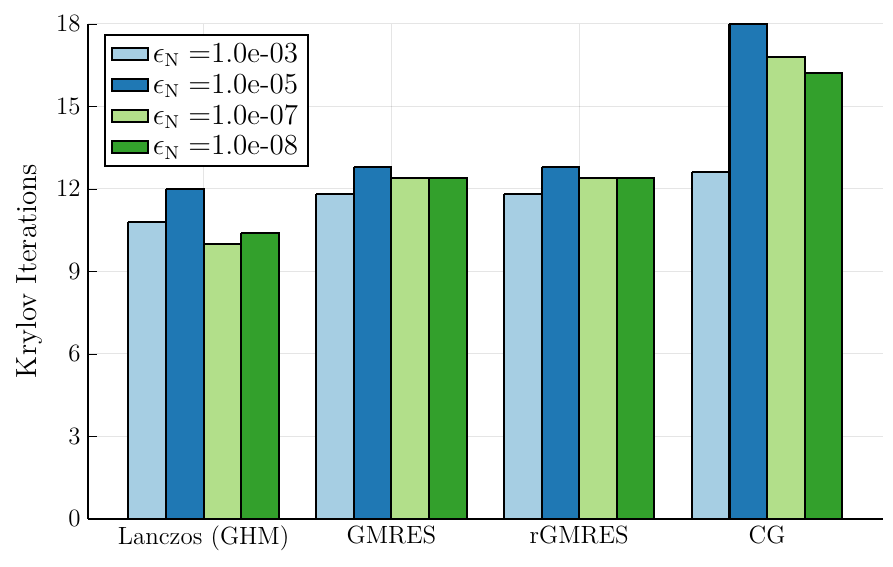}
        \caption{$\|\gk\| = 5.17\times 10^{-7}$}\label{fig.hilbert.small}
    \end{subcaptionblock}
    \caption{{Results of calculating a Newton-type direction for a perturbed Hilbert matrix and a GHM. The blue to green bars represent the Krylov iterations. The numbers shown in the legend correspond to different regularizations $\epsilon_{\mathrm{N}}$.}}
    \label{fig.hilbert}
\end{figure}

Consider the Hilbert matrix \cite{hilbert_beitrag_1970} as $H$ known to be ill-conditioned. A Hilbert matrix in the dimension of $n$ has the following analytic form:
\begin{align*}
    H_{i j}=\frac{1}{i+j-1}, i\le n, j\le n.
\end{align*}
It is known that the condition number $\kappa(H)$ grows in $O((1+{\sqrt {2}})^{4n}/{\sqrt {n}})$. {Under different $\epsilon_{\mathrm{N}} \in \{10^{-3}, 10^{-5},10^{-7},10^{-8}\}$ we compare the four mentioned algorithms ($\epsilon_{\mathrm{L}} = \epsilon_{\mathrm{N}}$). We set $n=300$ and randomly generate $\gk$ with a large norm and a small norm, and then collect the average number of Krylov iterations in \autoref{fig.hilbert}}.

{
The results basically demonstrate that when $\|\gk\|$ is large (see \autoref{fig.hilbert.large}), then for GHM (by \glanczos{}) Krylov iterations remain almost the same under different $\epsilon_{\mathrm{L}}$, while those of \newtoncg{}, \gmres{}, and \rgmres{} can grow with the condition number. This is consistent with our theoretical analysis that the complexity of solving GHMs is less influenced by the condition number of the Hessian matrix.
For the case where $\|\gk\|$ is very small (see \autoref{fig.hilbert.small}), there is not much difference between solving a Newton equation and the GHM.
We also note that for GHM it is quite safe to select a small $\epsilon_{\textrm{L}}$ without causing bad conditioning. This feature is nice in practice, since usually we will adaptively set the regularizer for a Newton method.
}

\paragraph{LIBSVM instances}

Next, consider the Newton-type systems arising from real-world least-square problems:
\begin{equation}\label{eq.linear.ols}
    \min_{\beta} ~ \frac{1}{2N}\sum_{i=1}^N\|x_i^T\beta-y_i\|^2 + \frac{\gamma}{2}\|\beta\|^2.
\end{equation}
where $x_i \in\real^n, \gamma >0, \beta\in\real^n$ and $N$ denotes the number of data points.
It is easy to see that the Newton system at some $\beta$ has the following form,
\begin{equation}\label{eq.linear.ols.newton}
    \left(\frac{1}{N} X^TX + \gamma\cdot I\right)\Delta\beta = - \frac{1}{N} X^T(X\beta-y).
\end{equation}
In this case $\epsilon_{\mathrm{N}} = \gamma$ is the perturbation; we test the above system with $\gamma \in \{10^{-3}, 10^{-4},10^{-5},10^{-6}\}$ on a set of problems in the LIBSVM library\footnote{\texttt{https://www.csie.ntu.edu.tw/\~{}cjlin/libsvmtools/datasets/}}.  Similarly, we keep track of Krylov iterations to find an iterate satisfying the corresponding residuals. We randomly generate five samples of $\beta$ from the uniform distribution and generate the Newton equation and the GHM therein. By computing the average number of Krylov iterations needed, \glanczos{}, \newtoncg{}, \gmres{} and finally \rgmres{} are presented in \autoref{tab.combined}. The results favor the GHMs when degeneracy is present; e.g., see instance \texttt{rcv1} and \texttt{news20}. In these cases, Lanczos has the best performance and does not deteriorate when $\gamma$ decays.
\begin{table}[ht]
    \centering
    \small
    \caption{Details of the datasets and average number of Krylov iterations for calculating one Newton-type direction or GHM for a linear least-square problem.}\label{tab.combined}
    \begin{tabular}{rcc|l|rrrr}
        \toprule
        \multicolumn{3}{c|}{Details}      & \multirow{2}{*}{Method}  & \multicolumn{4}{c}{$\gamma$}                                                             \\
        name                              & $n$                      & $N$                          &           & $10^{-3}$ & $10^{-4}$ & $10^{-5}$ & $10^{-6}$ \\
        \midrule
        \multirow{4}{*}{\texttt{a4a}}     & \multirow{4}{*}{122}     & \multirow{4}{*}{4781}
                                          & \gmres                   & 18.6                         & 29.2      & 31.6      & 32.2                              \\
                                          &                          &                              & \rgmres   & 18.6      & 29.2      & 31.6      & 32.2      \\
                                          &                          &                              & \newtoncg & 20.8      & 46.8      & 77.2      & 80.6      \\
                                          &                          &                              & \glanczos & 6.0       & 6.0       & 6.0       & 6.0       \\
        \midrule
        \multirow{4}{*}{\texttt{a9a}}     & \multirow{4}{*}{123}     & \multirow{4}{*}{32561}
                                          & \gmres                   & 18.6                         & 29.2      & 29.4      & 31.4                              \\
                                          &                          &                              & \rgmres   & 18.6      & 29.2      & 29.4      & 31.4      \\
                                          &                          &                              & \newtoncg & 20.6      & 44.6      & 73.2      & 70.2      \\
                                          &                          &                              & \glanczos & 6.0       & 6.0       & 6.0       & 6.0       \\
        \midrule
        \multirow{4}{*}{\texttt{w4a}}     & \multirow{4}{*}{300}     & \multirow{4}{*}{6760}
                                          & \gmres                   & 19.2                         & 32.0      & 39.0      & 42.8                              \\
                                          &                          &                              & \rgmres   & 19.2      & 32.0      & 39.0      & 42.8      \\
                                          &                          &                              & \newtoncg & 20.6      & 45.0      & 81.0      & 93.6      \\
                                          &                          &                              & \glanczos & 6.0       & 6.0       & 6.0       & 6.0       \\
        \midrule
        \multirow{4}{*}{\texttt{rcv1}}    & \multirow{4}{*}{47236}   & \multirow{4}{*}{20242}
                                          & \gmres                   & 14.0                         & 28.0      & 52.0      & 87.8                              \\
                                          &                          &                              & \rgmres   & 14.0      & 28.0      & 52.0      & 103.8     \\
                                          &                          &                              & \newtoncg & 13.0      & 33.2      & 84.4      & 199.0     \\
                                          &                          &                              & \glanczos & 5.0       & 5.0       & 5.0       & 5.0       \\
        \midrule
        \multirow{4}{*}{\texttt{covtype}} & \multirow{4}{*}{54}      & \multirow{4}{*}{581012}
                                          & \gmres                   & 8.0                          & 9.6       & 9.6       & 9.6                               \\
                                          &                          &                              & \rgmres   & 8.0       & 9.6       & 9.6       & 9.6       \\
                                          &                          &                              & \newtoncg & 7.0       & 11.6      & 14.0      & 14.0      \\
                                          &                          &                              & \glanczos & 6.0       & 6.0       & 6.0       & 6.0       \\
        \midrule
        \multirow{4}{*}{\texttt{news20}}  & \multirow{4}{*}{1355191} & \multirow{4}{*}{19996}
                                          & \gmres                   & 11.0                         & 20.0      & 35.0      & 60.0                              \\
                                          &                          &                              & \rgmres   & 11.0      & 20.0      & 75.0      & 71.2      \\
                                          &                          &                              & \newtoncg & 11.0      & 27.0      & 74.2      & 124.0     \\
                                          &                          &                              & \glanczos & 5.0       & 5.0       & 5.0       & 5.0       \\
        \bottomrule
    \end{tabular}
    \normalsize
\end{table}

{
\subsection{Characterization of primal-dual solutions}\label{sec.basic.pd}
A consequent curiosity arises for using GHMs as subproblems in a generic second-order method under the framework (\autoref{alg.gbasic}).
In preparation for the analysis of the proposed algorithms, we give a thorough characterization of primal-dual solutions to the subproblem \eqref{eq.ghqm} with respect to the control parameters $\delta_k$ and $\phi_k$.
Note that many of these analyses are particularly developed for GHM and are not needed for the ordinary model \cite{zhang_homogenous_2022}. We delay the proofs in this subsection to \autoref{app.appendix a}.
}
\subsubsection{Basic results}

Now we consider the case $\theta_k = 0$ where the equivalence to eigenvalue problems is lost.
\begin{lemma}[Existence of negative curvature]\label{lemma.delta upper bound}
    In \eqref{eq.ghqm}, $\theta_k = 0$ if and only if $F_k \succeq 0$.
\end{lemma}

The success of GHM and HSODF relies on $F_k$ being indefinite. Intuitively, the control variable $\delta_k$ cannot be too big otherwise, the optimal solution to \eqref{eq.ghqm} falls into the interior of the unit ball. We give a comprehensive characterization in \autoref{lem.strict.curv}. Note that $F_k \succeq 0$ occurs only when $H_k \succeq 0$, thus we present a specialized description when Hessian is positive semidefinite.

\begin{lemma}[Sufficiency in convex case]\label{lem.cvx.zero}
    Consider in the homogeneous model \eqref{eq.homoquadmodel}, suppose $\Hk \succeq 0$, then $\lambda_1(F_k) < 0$ as long as $\delta_k < \overline{\delta_k^{\mathsf{cvx}}}:=\phi_k^T H_k^\star \phi_k$.
\end{lemma}

The following Lemma describes the upper bound for $\theta_k$.
\begin{lemma}[Upper bound of $\theta_k$]\label{lemma.uppertheta}
    In GHM, it holds that:
    \begin{equation}\label{eq.lem.bound.theta}
        \theta_k \le \max\{-\delta_k, -\lambda_1(\Hk),0\} + \|\phi_k(x)\|.
    \end{equation}
\end{lemma}
Next, we move to the case where $t_k= 0$.
Let us recall the following lemmas on the spectrum of $F_k$ if $\phi_k \perp \mathcal S_1(\Hk)$.

\begin{lemma}[Lemma 3.1, 3.2, \citet{rojas_new_2001}]\label{lemma.rojas1}
    For any $q \in \mathcal S_j(\Hk), 1\le j \le r$, define
    $$p_j = -\left(\Hk-\lambda_j(\Hk) I\right)^{\star} \phi_k,~\tilde{\alpha}_j = \lambda_j(\Hk) - \phi_k^Tp_j,$$
    then
    \begin{enumerate}[(a)]
        \item $\left(\lambda_j(\Hk), [0; q]\right)$ is an eigenpair of $F_k$ if and only if $\phi_k \perp \mathcal S_j(\Hk)$.
        \item $\left(\lambda_j(\Hk), [1; p_j]\right)$ is an eigenpair of $F_k$ if and only if $\phi_k \perp \mathcal S_j(\Hk)$ and $\delta_k = \tilde \alpha_j$.
    \end{enumerate}
\end{lemma}
The above result basically says the eigenvector of $\Hk$ can be used to construct an eigenvector of $F_k$ if $\phi_k$ is orthogonal to the eigenspace. The dimension of $j$-th eigenspace for $F_k$ depends on whether $\delta_k$ coincides with the critical value $\tilde \alpha_j$ defined above. If $\phi_k \perp \mathcal S_1(\Hk)$, the above result only implies that $\lambda_1(\Hk)$ is an eigenvalue of $F_k$ but not necessarily the smallest one.
$\lambda_1(\Hk) = \lambda_1(F_k)$ holds only for some special value of $\delta_k$.
We are now ready to give the necessary condition to make $t_k$ zero.
\begin{corollary}\label{corr.necc.t0}
    In \eqref{eq.ghqm}, given $\phi_k \perp \mathcal S_1(\Hk)$, if $t_k= 0$ then $\delta_k \ge \tilde \alpha_1$. More specifically, we have
    \begin{enumerate}[(a)]
        \item if $\delta_k < \tilde \alpha_1$, then $\lambda_1(F_k) < \lambda_1(\Hk)$ and $t_k \neq 0$;
        \item if $\delta_k = \tilde \alpha_1$, then $\lambda_1(F_k) = \lambda_1(\Hk)$ and $[0;q]$, $[1;p_1]$ are the eigenvectors associated with $\lambda_1(F_k)$, where $q \in \mathcal{S}_1$ and $p_1$ is defined in \autoref{lemma.rojas1}, which further implies that  $t_k \in [0, 1]$.
        \item if $\delta_k > \tilde \alpha_1$,  then $\lambda_1(F_k) = \lambda_1(\Hk)$ and $t_k = 0$.
    \end{enumerate}
\end{corollary}

The proof is implied by the two lemmas before this corollary. However, it is worth noting that \emph{the reverse} is not true. Literally, if $\phi_k \perp \mathcal S_1(\Hk)$ and $\delta_k = \tilde \alpha_1$, then indeed $\lambda_1(F_k) = \lambda_1(\Hk)$. Nonetheless, $t_k$ may not be zero since the eigenvector can be a linear combination of $[0; q]$ and $[1; p_1]$ for some $q \in \mathcal S_1(\Hk)$. The sufficiency holds only when $\delta_k > \tilde \alpha_j$. Actually, $\phi_k \perp \mathcal S_1(\Hk)$ yet $t_k \neq 0$ is not really an undesired case. In the design of the algorithm, $\delta_k$ should not be increased
whenever $t_k = 0$, which is motivated by the
necessary condition.
As a byproduct, the following lemma characterizes the quantity of $\tilde \alpha_1$.
\begin{lemma}[Ordering of $\tilde \alpha_1$]\label{lemma.ordering of tilde alpha 1}
    Define $p_1, \tilde \alpha_1$ similar to \autoref{lemma.rojas1}, then the following holds for the smallest eigenvalue of our interest, $\lambda_1(\Hk) \le \tilde \alpha_1. $
\end{lemma}

\subsubsection{Continuity of auxiliary functions}\label{sec.aux}

To facilitate discussion, we treat the dual variable $\theta_k$, the norm of the step $\|\dk\|$, and the ratio $\theta_k(\delta_k)/\|\dk\|$ as functions of $\delta_k$, which are referred to as auxiliary functions.
Next, we discuss the continuity and differentiability of these functions.

\begin{definition}[Auxiliary functions of $\delta_k$]\label{def.aux}
    At each iterate $\xk$, consider the GHM with $\delta_k$ and let $v_k, t_k$ be the corresponding solution. $\overline \Delta_k < + \infty$ is an upper bound for the step.
    $$
        \begin{aligned}
            \Delta_k: \real \mapsto \real_+, ~\Delta_k(\delta_k) &
            :=\begin{cases}
                  \|v_k/t_k\|^2      & \textrm{if }~ \delta_k < \tilde \alpha_1 \\
                  \overline \Delta_k & \textrm{o.w. }~                          \\
              \end{cases}                                         \\
            \omega_k: \real \mapsto \real_+, ~\omega_k(\delta_k) & := \theta_k^2                                    \\
            h_k: \real \mapsto \real_+, ~h_k(\delta_k)           & := \frac{\omega_k(\delta_k)}{\Delta_k(\delta_k)}
        \end{aligned}
    $$
\end{definition}
In the sequel, we abbreviate $\Delta_k = \Delta_k(\delta_k), \omega_k = \omega(\delta_k)$ and $h_k = h_k(\delta_k)$. It is understood that these values, with the subscript $k$, correspond to the current iterate $\xk$. We should recall that $t_k \in [0, 1]$ is \emph{set-valued} if $\delta_k = \tilde \alpha_1$ (cf. \autoref{corr.necc.t0}).  In this view, $v_k/t_k$ is not well-defined. However, this occurs only when $\delta_k = \tilde \alpha_1$, so we are not bothered with this issue since any perturbation to $F_k$ will eliminate the set-valued case.

\begin{figure}[h]
    \centering
    \begin{subcaptionblock}[c]{0.47\textwidth}
        \includegraphics[width=0.98\textwidth]{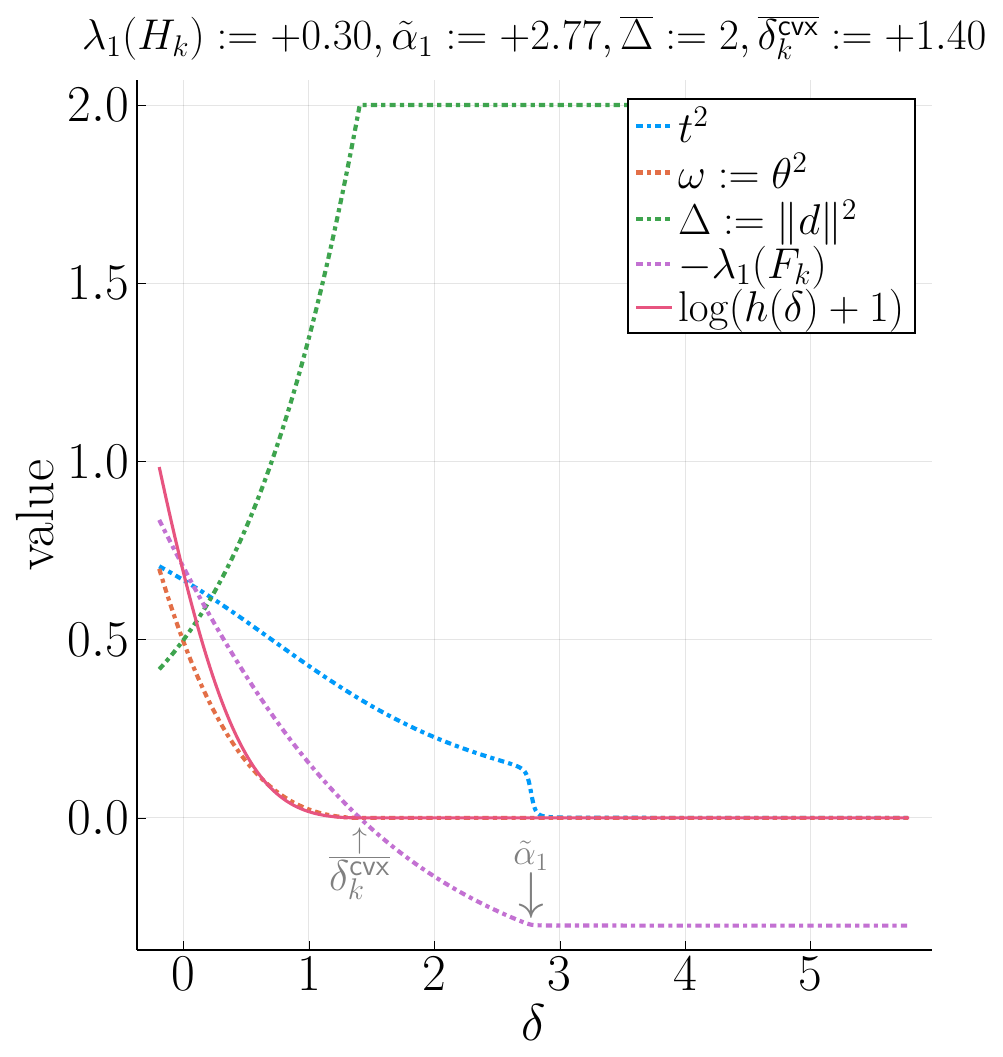}
        \caption{Convex case}
    \end{subcaptionblock}
    \begin{subcaptionblock}[c]{0.47\textwidth}
        \includegraphics[width=0.95\textwidth]{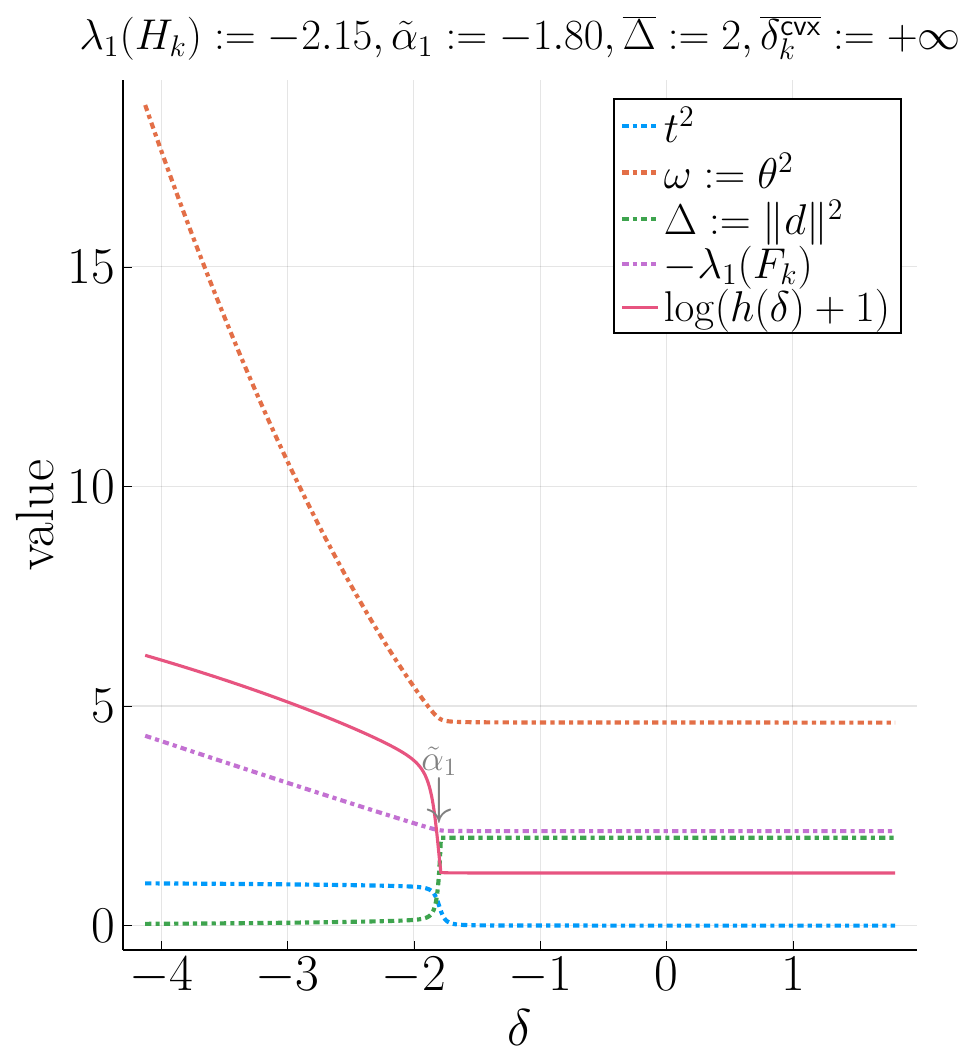}
        \caption{Non-convex case}
    \end{subcaptionblock}
    \caption{An illustration of ``perturbed'' case when $\gk \perp \mathcal S_1(\Hk)$, where $\phi_k = \gk + \varepsilon\cdot u_1, u_1 \in \mathcal S_1(\Hk)$ for better presentation. We also annotate the turning points $\overline{\delta_k^{\mathsf{cvx}}}$ and $\tilde \alpha_1$.}\label{fig.hod}
\end{figure}

In \autoref{fig.hod}, we give convex and nonconvex examples for the case where we \textit{almost} have $\gk \perp \mathcal S_1(\Hk)$. According to this example, we have several observations.
Firstly, in both convex and nonconvex cases, we see that $t_k$ can jump from a value close to $1$ to $0$ as $\delta_k \nearrow \tilde \alpha_1$ as defined in \autoref{lemma.rojas1}. Also, $\Delta_k$ is increasing in both cases and also has a clear jump because of $t_k$. In sharp comparison,
$\theta_k$ (and $\omega_k$) is continuous over $\delta_k \in \real$. Actually, the differentiability of $\theta_k$ actually connects to the case of $t_k = 0$.  We formalize these findings as follows.
\begin{lemma}[Continuity of $\theta_k, \omega_k$]\label{lem.continuity.theteomega}
    For $\phi_k \neq 0$, we have that
    $\theta_k, \omega_k$ is decreasingly convex and continuous for all $\delta_k \in \real$. Moreover, $\theta_k$ is differentiable at $t_k \neq 0$ such that
    \begin{equation}\label{eq.diff.omega}
        \diff{\delta_k} \theta_k = - \frac{1}{\Delta_k + 1}
    \end{equation}
\end{lemma}
Now we move on to the continuity of $\Delta_k$.

\begin{lemma}[Continuity of $\Delta_k$]\label{lem.property.h}
    For every iterate $\xk$, consider $\Delta_k$ as functions of $\delta_k$ in \autoref{def.aux}. Then $\Delta_k$ is continuous in $\delta_k$ if and only if $t_k \neq 0$, namely,
    \begin{enumerate}[(a)]
        \item If $\phi_k \not \perp \mathcal S_1(\Hk)$, then $\Delta_k$ is continuous for all $\delta_k \in \real$.
        \item If $\phi_k \perp \mathcal S_1(\Hk)$, then $\Delta_k$ is discontinuous at $\delta_k = \tilde \alpha_1$.
    \end{enumerate}
\end{lemma}
The proof is similar to \autoref{lem.continuity.theteomega}, which is thus omitted.  From the definition of $h_k$, we can see that its continuity also depends on the continuity of $\Delta_k$. Once $\phi_k \not \perp \mathcal{S}_1(\Hk)$, we show that $h_k$ is differentiable and monotone decreasing.
\begin{lemma}
    \label{lem.differentiablity of h}
    If $\phi_k \perp \mathcal{S}_1(\Hk)$, $h_k(\delta)$ is discontinuous at $\Tilde{\alpha}_1$; otherwise, $h_k(\delta)$ is continuous. Moreover, $h_k(\delta)$ is differentiable in $\delta$ is monotone decreasing.
\end{lemma}

Now we are ready to give the following definition of \emph{hard case}.
\begin{definition}[Hard case for HSODF]
    \label{def.hard case}
    We say the hard case occurs at iteration $k$ if $h_k(\delta)$ is discontinuous at $\tilde{\alpha}_1$.
\end{definition}
{
In HSODF, the \emph{hard case} only matters if $f$ is nonconvex at $\xk$.
As shown in \autoref{fig.hod}, if $f$ is convex near $\xk$, the function $h_k(\delta)$ would be discontinuous at $\tilde \alpha_1$, which is larger than the \emph{turning point} $\overline{\delta_k^{\mathsf{cvx}}}$ where $\theta_k = 0$. This will not affect our algorithm: if the aggregated $F_k$ has a positive leftmost eigenpair, we can decrease $\delta_k$, and the hard case does not matter.
By contrast, if $f$ is locally nonconvex, $t_k = 0$ does obstruct the computations as we can apply a slight perturbation to $\phi_k$ to escape the hard case.
}
\subsection{Using GHMs in a SOM}\label{sec.recover}
{
    We now briefly discuss how to use GHMs to recover a few standard second-order methods. For instance, we show a trust-region subproblem can be solved by using $O(\log(\epsilon^{-1}))$ GHMs via a bisection method. Furthermore, we can also provide an alternative to the gradient-regularized Newton step \cite{mishchenko_regularized_2023,doikov_gradient_2021,doikov_super-universal_2022}, and the search procedure in the inner loop requires $O(\log\log(\epsilon^{-1}))$ GHMs by using Newton's method (instead of the bisection).
}
\paragraph{Recovering trust-region methods}
Consider the trust-region methods for nonconvex optimization solving the following trust-region subproblems (TRS) over $\Delta > 0$:
\begin{equation}\label{eq.trsubp}
    \underset{\|d\| \le \Delta}{\min}~ m_k(d) := f(x_k) + \gk^Td + \frac{1}{2}d^T\Hk d.
\end{equation}
We could make use of GHM by setting $\phi_k$ exactly as $\gk$:
$F_k = [H_k, \gk; \gk^T, \delta_k]$.
Obviously, $\forall v\in\real^n, t\neq 0,$
$m_k(v/t)- f(x_k) + {\frac{1}{2}\delta_k}= {\frac{1}{2t^2}}{\psi_k(v,t;F_k)}.$
This leads to the following updates of HSODF:
$$
    \xkn = \xk - \dk, \; \dk := v_k/t_k.
$$
Using \autoref{alg.gbasic}, we are able to
recover trust-region methods with the same rules to update $\Delta$ in the outer loop $k$ of HSODF and repeatedly solve GHMs
in the inner loop $j$ to obtain the same solution of the trust-region subproblem (TRS). {The basic idea is searching for a proper $\delta_k$ such that $\|d_k\| = \Delta$ at each outer iterate $k$. Subsequently, the dual variable $\theta_k$ associated with GHM plays the same role as the dual variable of the trust-region subproblem.} We give the following result on the upper bound of the inner iteration number $\mathcal T_k$.
\begin{theorem}\label{thm.recover.trs}
    At each iteration $k$, for any radius $\Delta$, by setting $\delta_k$ properly, the solution $[v_k; t_k]$ also solves trust-region subproblem. Furthermore, the inner iteration $\mathcal T_k$ of searching the desired $\delta_k$ is upper bounded by $O(\log(\epsilon^{-1}))$ via a bisection procedure.
\end{theorem}
The proof is deferred to \autoref{proof.thm.recover.trs}.
The idea of using parametric eigenvalue problems to solve TRS has been explored in \cite{rojas_new_2001, adachi_solving_2017}.
However, these approaches limited themselves to solving the subproblem \eqref{eq.trsubp} instead of discussing an improvement of the solution in the iterative algorithm. To achieve the state-of-the-art convergence properties, the method could use the update rules for $\Delta$ in \cite{curtis_trust_2017,curtis_inexact_2019}.

\paragraph{Recovering gradient regularized methods}
The gradient-regularized Newton method was recently studied in \cite{mishchenko_regularized_2023,doikov_gradient_2021} for convex optimization, with the step defined as:
\begin{equation}\label{eq.quadsubp}
    x_{k+1} = x_k - (H_k + \gamma_k \|g_k\|^{1/2})^{-1} g_k,
\end{equation}
where $\gamma_k > 0$. Similarly, as long as $t_k \neq 0$, $d_k = v_k / t_k$ becomes a regularized Newton direction. A straightforward approach is to select a suitable $\delta_k$ such that $\theta_k \approx \Theta(\|g_k\|^{1/2})$, thus establishing an equivalence. We summarize this in the following theorem:
\begin{theorem}\label{thm.recover.reg}
    For any $\gamma_k > 0$, by choosing $\delta_k$ appropriately, the solution $[v_k; t_k]$ generates the same iterate as \eqref{eq.quadsubp}. Furthermore, the inner iteration $\mathcal T_k$ for finding the desired $\delta_k$ is bounded by $O(\log(\epsilon^{-1}))$ using a bisection method or $O(\log\log(\epsilon^{-1}))$ with Newton’s method.
\end{theorem}
The proof of this theorem is deferred to \autoref{proof.thm.recover.reg}. We do not further explore specialized variants here, as they warrant a separate discussion beyond the scope of this paper.

\begin{table}[h]
    \footnotesize
    \centering
    \caption{Using the homogeneous framework to recover or provide an alternative to other second-order methods by GHM with adaptive $\delta_k, \phi_k$.
        The inner complexity $\mathcal T_k$ represents the upper bound for the number of inner iterations associated with the outer iteration $k$.} \label{tab.ghm}
    \begin{tabular}{lllll}
        \toprule
        \multirow{2}{*}{Specific second-order method}                                       & \multirow{2}{*}{Inner complexity {$\mathcal T_k$} } & \multirow{2}{*}{Detailed discussion}                    \\
                                                                                            &                                                                                                               \\
        \midrule
        Trust-Region Method \cite{curtis_trust_2017,curtis_inexact_2019}                    & $O(\log(1/\epsilon))$                               & See \cite{rojas_new_2001} and \autoref{thm.recover.trs} \\
        Regularized Newton's Method \cite{mishchenko_regularized_2023,doikov_gradient_2021} & $O(\log\log(1/\epsilon))$                           & See \autoref{thm.recover.reg}                           \\
        {Adaptive HSODM}                                                                    & $O(\log(1/\epsilon))$                               & See \autoref{sec.ncvx}                                  \\
        Homotopy HSODM \cite{luenberger_linear_2021}                                        & {$\le 2$}                                           & See \autoref{sec.homotopy}                              \\
        \bottomrule
    \end{tabular}
    \normalsize
\end{table}

Broadly speaking, any Newton-type steps, either by regularization or trust region, can be recovered by a series of GHMs with a smart strategy to locate $\delta_k$ (and possibly $\phi_k$).
Besides, the complexity of the inner iteration $\mathcal T_k$ (to local $\delta_k$) depends on the calculus of primal-dual solutions to the specific homogeneous subproblem.

In this venue, we discuss an adaptive HSODM for nonconvex optimization under second-order Lipschitz continuity. The method extends naturally to convex optimization, which complements the original HSODM \cite{zhang_homogenous_2022} that only covers the nonconvex case with the line-search procedure. Within the adaptive HSODM, the dual variable $\theta_k$ can still be interpreted as a regularization term (see the local model \eqref{eq.local model} in \autoref{sec.ncvx}), and the number of inner iteration $\mathcal T_k$ is shown to be prototypical $O(\log(1/\epsilon))$.

Using a pronged approach,
we develop a new homotopy method in \autoref{sec.homotopy}, which has an $O(\log(\epsilon^{-1}))$ iteration complexity in certain structured non-strongly-convex optimization, with the number of inner iterations $\mathcal T_k \le 2$ for each iterate $\xk$.  For a quick view, we summarize the above discussions in \autoref{tab.ghm}.

\section{An Adaptive HSODM for General Nonconvex optimization}\label{sec.ncvx}
In this section, we consider a specific realization of HSODF named {\it adaptive} HSODM to nonconvex optimization problems with second-order Lipschitz continuous objective function. Our motivation comes from well-established results of \cite{nesterov_cubic_2006,cartis_adaptive_2011,cartis_adaptive_2011-1} and the recent monograph \cite{cartis_evaluation_2022}. The goal is to find an $\epsilon$-approximate second-order stationary point defined as follows.
\begin{definition}\label{def.sosp} A point $x$ is called an $\epsilon$-approximate second-order stationary points if it satisfies the following conditions.
    \begin{subequations}
        \begin{align}
            \label{eq.approxfocp} & \|\nabla f(x)\|\leq O(\epsilon)                                       \\
            \label{eq.approxsocp} & \lambda_{1} \left(\nabla^2 f(x)\right) \geq \Omega(-\sqrt{\epsilon}).
        \end{align}
    \end{subequations}
\end{definition}

In addition, we consider a broad class of objective functions satisfying the second-order Lipschitz continuity.
\begin{definition}
    \label{assm.lipschitz}
    We call a function $f$ has \(M\)-Lipschitz continuous Hessian if for all \( x, y \in \mathbb{R}^n\),
    \begin{equation}\label{eq.assm.lipschitz}
        \|\nabla^2 f(x) - \nabla^2 f(y)\| \le M \|x-y\|.
    \end{equation}
\end{definition}

We further make the following assumption on function $\phi_k(\xk)$ at every iterate $\xk$ in GHMs.
\begin{assumption}\label{assumption. phi close}
    Suppose that there exists a uniform constant $\varsigma_\phi>0$. Given an iterate $\xk \in \mathbb{R}^n$, if $t_k\neq 0$, then
    \begin{equation}
        \label{assumption.eq. phi close}  \|\phi_k(x_k) - \gk\| \leq \varsigma_\phi\|\dk\|^2
    \end{equation}
    where $\dk = v_k/t_k$.
\end{assumption}

The purpose of this condition is to ensure that $\phi_k$ remains sufficiently close to the gradient $\gk$, relative to the norm of the direction $\dk$. It is worth noting that the assumption is not restrictive. For example, if we simply choose $\phi_k = \gk$, the left-hand side of \eqref{assumption.eq. phi close} becomes zero, satisfying the condition automatically. In fact, this choice is likely sufficient in most cases. Now we are ready to present the adaptive HSODM for second-order Lipschitz continuous functions in \autoref{alg.gadaptive}.

\begin{figure}[h]
    \scriptsize
    \begin{algorithm}[H]
        \caption{The Adaptive HSODM}\label{alg.gadaptive}
        Initial point $x_0\in \mathbb{R}^n$, $\delta_0 \in \mathbb{R}, I_h=\mathbb{R}, h_{\min} > 0$, parameter $0<\eta_1<\eta_2<1$, $\gamma_1>1,\gamma_3\geq \gamma_2>1, 0<\gamma_4\leq 1$, $\sigma > 0$\;
        \For{$k =0, 1, 2, \dots$}{
        $\phi_k = \gk$, $\delta_{k,0}=\delta_{k-1}$

        \For{$j =0, 1, \dots, \mathcal T_k$}{
        \label{line.solution}  Obtain the solution pair $(\theta_{k,j}, \ [v_{k,j};t_{k,j}])$  of the subproblem
        \begin{equation}\label{eq.ghm.adaptive}
            \min_{\|[v;t]\| \leq 1} \begin{bmatrix}v \\ t\end{bmatrix}^T
            \begin{bmatrix}
                H_{k}          & \phi_{k,j}   \\
                (\phi_{k,j})^T & \delta_{k,j}
            \end{bmatrix}
            \begin{bmatrix}v \\ t\end{bmatrix}
        \end{equation}

        \label{line.perturb}  \If(\tcp*[f]{check hard case, see \autoref{sec.hard case}}){$t_{k,j} = 0$}{Go to \autoref{alg.hardcase};}
        Set $d_{k,j} = v_{k,j}/t_{k,j}$, $h_k(\delta_{k,j}) := \left(\theta_{k,j}/\|d_{k,j}\|\right)^2$\;
        \eIf{$\sqrt{h_k(\delta_k)} \in I_h$ within tolerance $\sigma$}{Set $\dk := d_{k,j}$, $\delta_k = \delta_{k,j}$\;\textbf{Break}}
        (\tcp*[f]{search $\delta_{k,j}$, see \autoref{sec.bisection procedure}}){
            Update $\delta_{k,j}$
        }
        }
        Compute
        \begin{equation*}
            \rho_k : = \frac{f(x_k+\dk)-f(x_k)}{m_k(\dk)-f(\xk)};
        \end{equation*}

        \label{line.ratios} \uIf(\tcp*[f]{very successful iteration}){$\rho_{k}> \eta_2$}{
            $I_h = \left[\max\left\{\sqrt{h_{\min}}, \gamma_4\sqrt{h_{k}(\delta_{k})}\right\},\sqrt{h_{k}(\delta_{k})}\right]$, \ $x_{k+1} = \xk+\dk$
        }
        \eIf(\tcp*[f]{successful iteration}){$\eta_1\leq \rho_{k}\leq \eta_2$}{
            $I_h =\left[\sqrt{h_{k}(\delta_{k})}/\gamma_1,\gamma_2 \sqrt{h_{k}(\delta_{k})}\right]$, \ $x_{k+1}=\xk+\dk$
        }
        (\tcp*[f]{unsuccessful iteration}){
            \label{line.ratioe}    $I_h =\left[\gamma_2 \sqrt{h_{k}(\delta_{k})},\gamma_3 \sqrt{h_{k}(\delta_{k})}\right]$, \ $x_{k+1}=\xk$
        }
        }
    \end{algorithm}
\end{figure}

The adaptive framework is largely motivated from \cite{cartis_adaptive_2011-1,cartis_evaluation_2022}.  We adopt a well-known ratio test on $\rho_k = \frac{f(x_k+\dk)-f(x_k)}{m_k(\dk)-f(\xk)}$ during iteration $k$:
if the ratio $\rho_k \geq \eta_2$, we call it a \textit{very successful iteration}; if $\eta_1 \leq \rho_k \leq \eta_2$, it is a \textit{successful iteration}; otherwise it is an \textit{unsuccessful iteration}. According to the ratio test, we adaptively construct a desired interval $I_h$ for $\sqrt{h_k}$, which serves as a coefficient of the regularizer of the following cubic model,
\begin{equation}\label{eq.local model}
    m_k(d):= f(\xk)+\phi_k^Td+\frac{1}{2}d^T \Hk d +\frac{\sqrt{h_k(\delta_k)}}{3}\|d\|^3 \,\,, \text{where} \,\, h_k(\delta_k) := \frac{\theta_k^2}{\|\dk\|^2}.
\end{equation}
Note again that we use $\phi_k$ to replace the first-order approximation.
    {As before, we focus on the GHM to minimize \eqref{eq.local model}. The following two cases are distinguished by the value of $t_k$ (line \ref{line.perturb}).
    }

    {
        If $t_k\neq 0$, then $h_k$ is continuous and smooth, and thus the solution $\dk$ to \eqref{eq.ghm.adaptive} minimizes \eqref{eq.local model}.
        To position the posterior $\sqrt{h_k(\delta_k)} \in I_h$, we resort to a bisection method based on an inner loop to solve GHMs labeled with $j$.
        A tolerance $\sigma$ is allowed when locating $\sqrt{h_k}\in [\sqrt{\ell},\sqrt{\nu}]$, e.g., $\sqrt{h_k}\in[\sqrt{\ell},\sqrt{\nu+\sigma}]$. We set $\sigma$ either adaptively: $\sigma\leq \Omega(h_k)$, or as a constant satisfying $\sigma < h_{\min}$. The overall convergence of adaptive HSODM is not sensitive to the choice of $\sigma$. As before, the number of inner iterations is bounded above by $O(\log(\epsilon^{-1}))$ from the analysis in \autoref{sec.aux}.
    }

    {
        If $t_k=0$, a perturbation on $\phi_k$ is applied in \autoref{alg.hardcase}. We later show it eventually produces a successful iterate so that the main algorithm (\autoref{alg.gadaptive}) proceeds to the next iteration. For fluency of our presentation, we ignore the ``hard case'' and bisection method for now, and the details will be revisited later.
    }

\subsection{Convergence analysis}

\subsubsection{Global convergence}
We first introduce some technical lemmas to help us identify the sufficient decrease if an iterate is accepted.
To this end, we rewrite the optimality condition of subproblem \eqref{eq.ghm.adaptive} in \autoref{alg.gadaptive} with respect to the function $h_k(\delta_k)$.
\begin{lemma}[Optimality condition for \eqref{eq.ghm.adaptive}]
    \label{lem.second-order cr opt cond}
    Suppose the current iterate $\xk$ does not fall into the hard case defined in \autoref{def.hard case}. Then the optimal primal-dual solution $([v_{k};t_{k}], -\theta_{k})$ of \eqref{eq.ghm.adaptive} satisfies the following conditions:
    \begin{subequations}\label{eq.mod.opt}
        \begin{align}
            \label{eq.firstorderoptc}
            \phi_k+\Hk \dk+\sqrt{h_k(\delta_k)}\|\dk\|\dk = 0 \\
            \label{eq.secondorderoptc}
            \Hk+\sqrt{h_k(\delta_k)}\|\dk\| I \succeq 0,
        \end{align}
    \end{subequations}
    where $d_k = v_k / t_k$. That is, viewing $h_k(\delta_k)$ as a constant here, $d_k$ is the minimizer of $m_k(d)$ defined in \eqref{eq.local model}.
\end{lemma}
\begin{proof}

    These facts can be recognized by substituting $\theta_k$ with $\sqrt{h_k(\delta_k)}\|d_k\|$ into the optimality conditions.
\end{proof}
We give the following estimate of model decrease.
\begin{lemma}\label{lem.model decrease}
    Suppose that $\xk$ does not fall into the hard case, we have the following model decrease
    \begin{equation}
        \label{eq. model succ dec}
        f(\xk)-m_k(\dk) \geq \frac{\sqrt{h_k(\delta_k)}}{6}\|\dk\|^3,
    \end{equation}
    where $m_k(\cdot)$ is defined in \eqref{eq.local model}.
\end{lemma}
\begin{proof}
    Note:
    \begin{align}
        \notag
        f(\xk)-m_k(\dk) & = -\phi_k ^T\dk-\frac{1}{2}\dk^T \Hk\dk-\frac{\sqrt{h_k(\delta_k)}}{3}\|\dk\|^3                \\
        \notag
                        & = \frac{1}{2}\dk^T\Hk\dk + \frac{2\sqrt{h_k(\delta_k)}}{3}\|\dk\|^3                            \\
        \notag
                        & = \frac{1}{2}\dk^T(\Hk+\sqrt{h_k(\delta_k)}\|\dk\|)\dk+\frac{\sqrt{h_k(\delta_k)}}{6}\|\dk\|^3 \\
                        & \geq \frac{\sqrt{h_k(\delta_k)}}{6}\|\dk\|^3.\label{eq.modeldecrease}
    \end{align}
    The second equation is due to \eqref{eq.firstorderoptc}, and the last inequality comes from \eqref{eq.secondorderoptc}.
\end{proof}

Furthermore, if the current iterate $x_k$ is accepted, the gradient norm of the next iteration could be bounded above.

\begin{lemma}\label{lem.step and g}
    {Suppose that $\xk$ does not fall into the hard case and \autoref{assumption. phi close} holds}, if the $k$-th iteration is successful, then we have the following relation between the step $\dk$ and the gradient at the new point $\xkn$:
    \begin{equation}
        \label{eq.relation step and g}
        \|\gkn \| \leq \frac{2\sqrt{h_k(\delta_k)}+M+2\varsigma_\phi}{2}\|\dk\|^2.
    \end{equation}
\end{lemma}
\begin{proof}
    Taking the norm of both sides of the equation \eqref{eq.firstorderoptc} we have
    \begin{equation}
        \label{eq.relation step and apg}
        \| \phi_k+\Hk\dk\| = \sqrt{h_k(\delta_k)}\|\dk\|^2.
    \end{equation}
    Then observe
    \begin{align}
        \notag
        \|\gkn\| & \leq \|\gkn-\gk-\Hk\dk\| + \| \phi_k-\gk\| +\| \phi_k + \Hk\dk\|                               \\
                 & \leq \frac{M}{2}\|\dk\|^2 +\sqrt{h_k(\delta_k)}\|\dk\|^2 +\|\phi_k -\gk\|.\label{eq.final gkn}
    \end{align}
    where the last inequality follows from second-order Lipschitz condition.
    Combining \autoref{assumption. phi close}, we have the desired result.
\end{proof}

Equipped with the \autoref{lem.step and g} and \autoref{lem.model decrease}, the function decreases at iterate $\xk$ can also be estimated in either the gradient or Hessian if the step is accepted.
\begin{lemma}\label{lem.fo decrease}
    Suppose that $\xk$ does not fall into the hard case and \autoref{assumption. phi close} holds, if the $k$-th iteration is successful, then we have the following function decrease:
    \begin{equation}
        \label{eq.function decrease fo}
        f(\xk)-f(\xkn) \geq \eta_1\frac{\sqrt{h_k(\delta_k)}}{12}\bigg( \frac{2}{M+2\sqrt{h_k(\delta_k)}+2\varsigma_\phi}\bigg)^{3/2}\|\gkn\|^{3/2}.
    \end{equation}
\end{lemma}
\begin{proof}
    Plugging \eqref{eq.relation step and g} into \eqref{eq. model succ dec}  we can have
    \begin{equation*}
        f(\xk)-m_k(\dk) \geq \frac{\sqrt{h_k(\delta_k)}}{12}\bigg( \frac{2}{M+2\sqrt{h_k(\delta_k)}+2\varsigma_\phi}\bigg)^{3/2}\|\gkn\|^{3/2} > 0.
    \end{equation*}
    Provided that the $k$-th iteration is successful, we have \eqref{eq.function decrease fo}.
\end{proof}

\begin{lemma}
    \label{lem.so decrease}
    Suppose that $\xk$ does not fall into the hard case, if the $k$-th iteration is successful, then we have the following function decrease
    \begin{equation}
        \label{eq.so decsnet}
        f(\xk)-f(\xk+\dk) \geq -\frac{\eta_1 \left(\lambda_{1} (\Hk)\right)^3}{6h_k(\delta_k)}.
    \end{equation}
\end{lemma}
\begin{proof}
    Note that $h_k(\delta_k) \geq h_{\min} > 0$, then we plug the inequality \eqref{eq.secondorderoptc} into the model decrease \eqref{eq. model succ dec}, which completes the proof.
\end{proof}

Next, we show that as long as $\phi_k$ is close enough to $\gk$ as in \autoref{assumption. phi close}, the number of unsuccessful iterations will be bounded above by a quantity determined by the number of successful iterations.
We begin with a lemma describing the upper bound on the regularization $h_k(\delta_k)$. Note that this also implies RHS of \eqref{eq.relation step and g} is bounded from above.

\begin{lemma}\label{lem.upper bound h}
    In \autoref{alg.gadaptive}, suppose \autoref{assumption. phi close} holds, then $h_k(\delta_k)$ has a uniform upper bound for all $k \ge 1$, i.e.,
    \begin{equation}
        \label{eq.uniform ub}
        h_k(\delta_k) \leq \max \left\{h_{0}(\delta_{0}),9\gamma_3^2\left(\frac{M}{2}+\varsigma_\phi\right)^2 \right\} =: \varsigma_h
    \end{equation}
\end{lemma}
\begin{proof}
    We only need to prove that when $h_k(\delta_k) \geq 9(\frac{M}{2}+\varsigma_\phi)^2$, the iteration $k$ must be very successful. Namely,
    \begin{align}
        \nonumber    m_k(\dk)-f(\xk+\dk) & =(\phi_k-\gk)^T\dk+ \frac{1}{2}\dk^T (\Hk - \nabla^2 f(\xk+\xi\dk))\dk+\frac{\sqrt{h_k(\delta_k)}}{3}\|\dk\|^3 \\
        \nonumber                        & \geq -\varsigma_\phi \|\dk\|^3-\frac{M}{2}\|\dk\|^3+\frac{\sqrt{h_k(\delta_k)}}{3}\|\dk\|^3                    \\
        \label{eq.diff bet modeland f}   & = \left(\frac{\sqrt{h_k(\delta_k)}}{3}-\frac{M}{2}-\varsigma_\phi\right)\|\dk\|^3,
    \end{align}
    where $\xi \in \left[0,1\right]$. Therefore, $m_k(\dk)-f(\xk+\dk)\geq 0$ holds as long as $h_k(\delta_k) \geq 9(\frac{M}{2}+\varsigma_\phi)^2$.

    Then, together with \autoref{lem.model decrease}, the ratio $\rho_k$ follows,
    \begin{align}
        \notag
        \rho_k & = \frac{f(\xk)-f(\xkn)}{f(\xk)-m_k(\dk)}                                   = \frac{f(\xk)-m_k(\dk)+m_k(\dk)-f(\xkn)}{f(\xk)-m_k(\dk)} \\
               & = 1+\frac{m_k(\dk)-f(x_k + d_k)}{f(\xk)-m_k(\dk)}  \geq 1.\label{eq.rhok}
    \end{align}
    Hence iteration $k$ must be very successful, and we get the desired bound.
\end{proof}

\begin{corollary}
    Suppose \autoref{assumption. phi close} holds, then there exists constant ${\varsigma_C}>0$ such that for
    any successful iteration $k$, we have
    \begin{equation}\label{eq.function decrease fo cons}
        f(\xk)-f(\xkn) \geq {\varsigma_C}\|\gkn\|^{3/2},
    \end{equation}
    where ${\varsigma_C} :=\eta_1\frac{\sqrt{h_{\min}}}{12}\left( \frac{2}{M+2\sqrt{h_{\min}}+2\varsigma_\phi}\right)^{3/2}$.
\end{corollary}
\begin{proof}
    Note that $h_k$ is bounded both from below and above. As a function of $h_k$, the expression  $$ \eta_1\frac{\sqrt{h_k(\delta_k)}}{12}\left( \frac{2}{M+2\sqrt{h_k(\delta_k)}+2\varsigma_\phi}\right)^{3/2}$$
    is continuous over the closed interval $[h_{\min},\varsigma_h]$, and thus it has a minimal value ${\varsigma_C}$. Through some basic analysis, we know that the function is increasing and thus attains minimum at $h_k = h_{\min}$. Hence we derive the value of ${\varsigma_C}$.
\end{proof}

\begin{corollary}\label{coro.successive suc iter}
    Suppose \autoref{assumption. phi close} holds, for any two successive successful iterates, without loss of generality, we denote them as $j$ and $j+1$, then exactly one of the following two cases will occur:
    \begin{enumerate}[(a)]
        \item If $x_{j+1}$ satisfies \eqref{eq.approxfocp} and \eqref{eq.approxsocp}, this means $x_{j+1}$ is a point satisfying \autoref{def.sosp}.
        \item {otherwise, $f(x_{j}) - f(x_{j+2}) \ge \Omega(\epsilon^{3/2})$.}
    \end{enumerate}
\end{corollary}
\begin{proof}
    Suppose that the first case does not hold, then we must have:
    $$\|g_{j+1} \| \geq \Omega(\epsilon) ~\text{or}~ \lambda_{1} (\nabla ^2 f(x_{j+1})) \leq O(-\sqrt{\epsilon}).$$ Recall the function value decreases in \eqref{eq.function decrease fo cons} or \eqref{eq.so decsnet}, we have
    \begin{equation*}
        f(x_j)-f(x_{j+2})\geq \max\left\{-\frac{\eta_1 \lambda_1(H_{j+1})^3}{6 \varsigma_h} ,{\varsigma_C} \|g_{j+1}\|^{3/2}\right\}\geq {\varsigma_f} \epsilon^{3/2},
    \end{equation*}
    where ${\varsigma_f} = \min \left\{{\varsigma_C},\frac{\eta_1}{6\varsigma_h}\right\}$, thus the second statement holds.
\end{proof}

Now, we are ready to provide the complexity analysis of the \autoref{alg.gadaptive}.
We define the following sets frequently used in the subsequent analysis:
the index sets of successful iterations $\mathcal{S}_j$ and unsuccessful iterations $\mathcal{U}_j$
up to the given iteration $j$, i.e.,
\begin{equation}
    \label{eq.sucandunsucsetj}
    \mathcal{S}_j: = \{k\leq j: {\rho_k\geq \eta_1}\} \quad \text{and} \quad \mathcal{U}_j:= \{k\leq j: \rho_k < \eta_1\}.
\end{equation}

As a matter of fact, the cardinality of $\mathcal{U}_j$ can be upper bounded by the cardinality of $\mathcal{S}_j$ as described in the following lemma.
\begin{lemma}
    For any $j\geq 0$, let $\mathcal{S}_j$ and $\mathcal{U}_j$ be defined in \eqref{eq.sucandunsucsetj}. Then
    \begin{equation}
        \label{eq.bound unsuc by suc}
        \vert \mathcal{U}_j \vert \leq \frac{1}{\log \gamma_2}\left[\frac{1}{2}\log \frac{\varsigma_h}{h_{0}(\delta_{0})}+\vert \mathcal{S}_j \vert \log \frac{1}{\gamma_4}\right]
    \end{equation}
\end{lemma}
\begin{proof}
    First, by the mechanism of the algorithm, we have
    \begin{equation}\label{eq.relationsuc}
        \gamma_4\sqrt{h_k(\delta_k)} \leq \sqrt{h_{k+1}(\delta_{k+1})}, \quad k\in \mathcal{S}_j,
    \end{equation}
    and
    \begin{equation}\label{eq.relationunc}
        \gamma_2\sqrt{h_k(\delta_k)}\leq \sqrt{h_{k+1}(\delta_{k+1})},\quad k\in \mathcal{U}_j.
    \end{equation}
    Recall the results in \autoref{lem.upper bound h}, there exists an upper bound for $h_k$, by deducing inductively, we conclude that
    \begin{equation}\label{eq.deduce}
        \sqrt{h_{0}(\delta_{0})}\gamma_4^{\vert\mathcal{S}_j\vert}\gamma_2^{\vert\mathcal{U}_j\vert}\leq \sqrt{\varsigma_h},
    \end{equation}
    which is equivalent to
    \begin{equation}
        \label{eq.deducelog}
        \vert \mathcal{S}_j \vert \log \gamma_4+ \vert \mathcal{U}_j \vert \log \gamma_2 \leq \frac{1}{2}\log \frac{\varsigma_h}{h_{0}(\delta_{0})}.
    \end{equation}
    By rearranging terms in \eqref{eq.deducelog} and we have the desired result.
\end{proof}
Next, we assume that $j_f$ is the first index such that the following holds: $$\| g_{j_f+1}\| \leq O(\epsilon), \ \lambda_{1}(H_{j_f+1})\geq \Omega(-\sqrt{\epsilon}),$$
that is, the next iterate is already an $\epsilon$-approximated second-order stationary point.
We present an upper bound for the cardinality of the set $\vert \mathcal{S}_{j_f} \vert$.
\begin{lemma}
    The cardinality of the set $\mathcal{S}_{j_f}$ satisfies
    \begin{equation}
        \label{eq.card of suc}
        \vert \mathcal{S}_{j_f} \vert \leq \frac{f_0-f_{\text{low}}}{{\varsigma_f}}\epsilon^{-3/2} = O(\epsilon^{-3/2}).
    \end{equation}
\end{lemma}
\begin{proof}
    First, we denote $k^s_{j}$ as the $j$-th element in $\mathcal{S}_{j_f}$. That is $\mathcal{S}_{j_f} = \{ k^s_{1},k^s_{2},\ldots,k^s_{\vert \mathcal{S}_{j_f}\vert} \}$.
    Note that for any $k^s_j$, from \autoref{coro.successive suc iter}, we have
    \begin{equation}
        \label{eq.suc dec ult}
        f(x_{k^s_j})-f(x_{k^s_{j+2}})\geq {\varsigma_f} \epsilon^{3/2}.
    \end{equation}

    Therefore, for any two successive successful iterates before $j_f$, the function value must decrease of $O(\epsilon^{3/2})$.
    Without loss of generality, we assume that the first iterate was successful, then we have
    \begin{align}
        \label{eq.telscop}
        f(x_0)-f(x_{j_f+1}) & \geq \sum_{j=\vert\mathcal{S}_{j_f}
            \vert\bmod 2 +1}^{\vert\mathcal{S}_{j_f}\vert-1} [f(x_{k^s_j})-f(x_{k^s_{j+2}})] \geq {\varsigma_f} \vert \mathcal{S}_{j_f}\vert \epsilon^{3/2}.
    \end{align}
    Note in \eqref{eq.telscop} ``\textrm{mod}'' is the modulo operation. Note that ${\varsigma_f}>0$ is a constant depending on the problem parameters such as $M$, hence we conclude that $$\vert \mathcal{S}_{j_f}\vert \leq \frac{f(x_0)-f_{\text{low}}}{{\varsigma_f}}\epsilon^{-3/2} = O(\epsilon^{-3/2}).$$
\end{proof}

\begin{theorem}\label{thm.global}
    The adaptive HSODM takes $O\left(\epsilon^{-3/2}\right)$ iterations to achieve a point $\xk$ satisfying $\|\gk\| \leq O(\epsilon)$ and $\lambda_{1}(\Hk) \geq \Omega(-\sqrt{\epsilon})$.
\end{theorem}
\begin{proof}
    Note that $x_{j_f+1}$ is the point satisfying \autoref{def.sosp}. The total iteration number is the sum of $\vert \mathcal{S}_{j_f} \vert$ and $\vert \mathcal{U}_{j_f} \vert$. The results are a simple combination of \eqref{eq.card of suc} and \eqref{eq.bound unsuc by suc}.
\end{proof}

Apart from the number of iterations consumed by the algorithm, we are also interested in the total number of GHMs solved by the algorithm.
Note that the algorithm has to find an $h_k$ in a desired interval, which means we have to search for an appropriate $\delta_k$. By the analysis in \autoref{sec.aux} we may use a bisection procedure since $h_k$ is continuous in most cases.
In particular, locating $h_k$ brings computations in solving extra GHMs in the order of $O(\log(\epsilon^{-1}))$ (see \autoref{sec.bisection procedure}).

\begin{theorem}
    \label{thm.complexity bisection procedure}
    At an iterate $\xk$ of the adaptive HSODM, the number of iterations of the bisection procedure is
    \begin{equation}
        \label{eq.complexity of bisection}
        O\left(\log \left(\frac{\varsigma_h U_\phi U_H}{h_{\min}\sigma}\right)\right).
    \end{equation}
    If we set the inexactness of the bisection procedure $\sigma=\epsilon$, then the iteration bound becomes
    \begin{equation}
        \label{eq.complexity of bisection 2}
        O\left(\log \left(\frac{\varsigma_h U_\phi U_H}{h_{\min}}\epsilon^{-1}\right)\right).
    \end{equation}
\end{theorem}

We arrive at the conclusion on the total number of  GHMs solved in \autoref{alg.gadaptive}.
\begin{theorem}
    \label{coro.final complexity}
    Let $\mathcal K_{\psi}$ be the total number of evaluations to solve GHMs in \autoref{alg.gadaptive}.
    Suppose the tolerance for the bisection procedure is set as $\sigma = \epsilon$, then we have the following bound before reaching an iterate satisfying \eqref{eq.approxfocp} and \eqref{eq.approxsocp}:
    \begin{equation}
        \label{eq.oracle calls}
        \mathcal K_{\psi} = O\left(\epsilon^{-3/2}\log(\epsilon^{-1})\right).
    \end{equation}
\end{theorem}
\begin{proof}
    It is directly implied by \autoref{thm.global} and \autoref{thm.complexity bisection procedure}.
\end{proof}

\subsubsection{Local convergence}

Now we move on to the analysis of the local performance of the Adaptive HSODM,
where we assume that the adaptive HSODM converges to a nondegenerate local minimum.
\begin{assumption}\label{assm.local}
    The sequence $\{\xk\}$ generated by adaptive HSODM converges to $x^*$ such that
    \begin{equation*}
        H(x^*)\succ \mu I\quad  \mbox{for some}\quad \mu>0.
    \end{equation*}
    where $H(x^*)\succ \mu I$ for some $\mu>0$.
\end{assumption}
Since $\xk \to x^*$ and the  $H(x)$ is Lipschitz continuous, we have that
\begin{equation}
    \label{eq.local strongly cvx}
    H(\xk) \succeq \mu I \quad \text{for some sufficiently large $k$ .}
\end{equation}
Then we are ready to analyze the local performance of adaptive HSODM.
\begin{lemma}
    \label{lem.step and current grad}
    Suppose \autoref{assumption. phi close} and \autoref{assm.local} hold, then for all sufficiently large $k$, the corresponding iteration is very successful, i.e., $\rho_k \geq \eta_2$, and
    \begin{equation}\label{eq.step and current grad}
        \|d_k\| \leq \frac{1}{\mu}\|\phi_k\|.
    \end{equation}
\end{lemma}
\begin{proof}
    First we prove that $\rho_k\geq \eta_2$ for sufficiently large $k$, note that in our algorithm we can always guarantee $m_k(\dk)-f(\xk)<0$ once $d_k \neq 0$, then we can define
    \begin{equation*}
        r_k := f(x_k+\dk)-m_k(\dk)+(1-\eta_2)(m_k(\dk)-f(\xk)),
    \end{equation*}
    $\rho_k\geq \eta_2$ is equivalent to $r_k\leq 0$.  It remains to show that $r_k \leq 0$ as $k \rightarrow +\infty$. For the term $f(x_k+\dk)-m_k(\dk)$, by Taylor expansion, it could be upper bounded as
    \begin{align}
        \notag
        f(\xk+\dk)-m_k(\dk) & = f(\xk)+\gk^T \dk +\frac{1}{2}\dk^T H(\xk+\xi\dk)\dk                                                                \\
        \notag              & \qquad -f(\xk)-\phi_k^T\dk-\frac{1}{2}\dk^T \Hk\dk-\frac{\sqrt{h_k(\delta_k)}}{3}\|\dk\|^3                           \\
        \notag
                            & \leq \frac{1}{2}\dk^T (H(\xk+\xi\dk)-\Hk)\dk+(\gk-\phi_k)^T\dk                                                       \\
                            & \leq \frac{1}{2}\|H(\xk+\xi\dk)-\Hk\| \|\dk\|^2+\varsigma_\phi \|\dk\|^3,\label{eq.difference of function and model}
    \end{align}
    where $\xi \in (0,1)$. For the second term $m_k(\dk)-f(\xk)$, from \eqref{eq.firstorderoptc} we have
    \begin{align}
        \notag
        f(\xk)-m_k(\dk) & = -\phi_k^T \dk-\frac{1}{2}\dk^T \Hk\dk-\frac{\sqrt{h_k(\delta_k)}}{3}\|\dk\|^3  \\
        \notag
                        & = \frac{1}{2}\dk^T \Hk \dk + \frac{2}{3}\sqrt{h_k(\delta_k)}\|\dk\|^3            \\
                        & \geq \frac{1}{2}\mu \|\dk\|^2.\label{eq.difference of function and model step k}
    \end{align}
    Substituting the inequalities \eqref{eq.difference of function and model} and \eqref{eq.difference of function and model step k} into the definition of $r_k$, we conclude that
    \begin{equation}
        \label{eq.rk}
        r_k \leq \frac{1}{2}\|\dk\|^2 \bigg\{ \|H(\xk+\xi\dk)-\Hk\|+2\varsigma_\phi \|\dk\|-(1-\eta_2)\mu\bigg\},\quad \xi\in(0,1).
    \end{equation}
    Since we assume that $\xk \to x^*$, we have $\dk\to 0$ when $k\to +\infty$.
    As a result, we can conclude that $r_k \leq 0$ when $k\to +\infty$, which finishes the first part of the lemma.

    For the second part, from the optimal condition \eqref{eq.firstorderoptc}, we have
    \begin{equation}
        \label{eq.second part lem}
        -\phi_k^T \dk  = \dk^T \Hk \dk +\sqrt{h_k(\delta_k)}\|\dk\|^3.
    \end{equation}
    Therefore, for the sufficiently large $k$, we observe that
    \begin{equation}
        \label{eq.second part lem 2}
        \mu \|\dk\|^2\leq \dk^T \Hk \dk+\sqrt{h_k(\delta_k)}\|\dk\|^3 = -\phi_k^T \dk \leq \|\phi_k\|\|\dk\|.
    \end{equation}
    By rearranging the terms, we complete the proof.
\end{proof}
The following result states that adaptive HSODM has a local quadratic convergence rate.
\begin{theorem}
    \label{thm.local}
    Suppose \autoref{assumption. phi close} and \autoref{assm.local} hold, then the convergence rate of adaptive HSODM is quadratic, i.e., we have
    \begin{equation}
        \label{eq.local quadratic g}
        \lim_{k\to +\infty} \frac{\|\phi_{k+1}\|}{\|\phi_k\|^2} \leq \frac{2\sqrt{\varsigma_h}+M+4\varsigma_\phi}{2\mu^2}.
    \end{equation}
\end{theorem}
\begin{proof}
    First, recall the results of \autoref{lem.step and g} and \autoref{lem.upper bound h}, we obtain that
    $$
        \|g_{k+1}\| \leq \frac{2\sqrt{\varsigma_h}+M+2\varsigma_\phi}{2}\|d_k\|^2.
    $$
    Hence, for the norm of $\phi_{k+1}$, by the \autoref{assumption. phi close}, it implies that
    $$
        \|\phi_{k+1}\| \leq \|\phi_{k+1} - g_{k+1}\| + \|g_{k+1}\| \leq \frac{2\sqrt{\varsigma_h}+M+4\varsigma_\phi}{2}\|d_k\|^2.
    $$
    Therefore, with the inequality \eqref{eq.relation step and g} in place, we have
    \begin{equation}
        \sqrt{\frac{2}{2\sqrt{\varsigma_h}+M+4\varsigma_\phi}}\|\phi_{k+1}\|^{\frac{1}{2}}\leq \|\dk\| \leq \frac{1}{\mu} \|\phi_k\|.
    \end{equation}
    Taking limits of the above, we can have the desired result.
\end{proof}

\subsection{A discussion on the hard case}\label{sec.hard case}
In this part, we complement previous analyses by providing treatment on the hard case (line \ref{line.perturb}) of \autoref{alg.gadaptive}. Recall that, to achieve similar convergence properties in \cite{cartis_adaptive_2011,curtis_trust_2017}, the adaptive HSODM relies on performing a bisection procedure over $\delta_k$ to keep $\sqrt{h_k(\delta_k)}$ lying in a given interval at the iterate $x_k$.
The challenge in the hard case is that the continuity of $h_k$ at $\tilde{\alpha}_1$ (cf. \autoref{lem.differentiablity of h}) is deteriorated, and thus a suitable $\delta_k$ cannot be found.

Fortunately, the flexibility of GHM \eqref{eq.ghqm} provides us with a remedy that utilizes $\phi_k$ with a suitable perturbation over $\gk$. This treatment enables us to escape the hard case and achieve a successful iteration. In particular, such simple perturbation makes two desired outcomes simultaneously: $(a)$ \autoref{assumption. phi close} holds; $(b)$ The function value $h_k(\delta)$ lies in the given interval.
Note that if $\Hk \succeq 0$, then $\tilde \alpha_1 > 0$, in which case the turning point $\overline{\delta^{\mathsf{cvx}}} \leq \tilde \alpha_1$ and $h_k(\Tilde{\alpha}_1) = 0$, as a consequence, $\Tilde{\alpha}_1$ is always out of the target interval. {Thus we only focus on the case that the GHMs are indefinite, since by the design of our method, the hard case will not occur when the Hessian at the current point is positive semi-definite}. Furthermore, we also prevent $\delta_k$ from being too large following the guidance in \autoref{lem.cvx.zero} and \autoref{corr.necc.t0}, which is achieved by using a lower bound such that $h_k > h_{\min} > 0$. Consequently, the hard case occurs at $k$-th iteration \emph{only} if $\Hk \prec 0$ and the preceding $(k-1)$-th iteration is successful.

We now describe the following \autoref{alg.hardcase} that use perturbed gradient $\phi_k$. For better understanding, we denote $k$ as the current iterate, then the subsequent iterates are denoted by $i$: $k,k+1,...,k+i$ and so forth. Similar to the previous discussion, we stick to $\lambda_1 = \lambda_{1}(\Hk)$. Since the hard case appears, $v_k$ is now the leftmost eigenvector. The idea is as follows.  For each of the following iterate $i$, we perturb $\phi_{k+i}$ based on the preceding $h_{k+i-1}$. When gradually increasing $h_{k+i}$, we use the same bisection method indexed by $j$ to find $\delta_{k+i,j}$. We show that it will finally produce a successful iteration. Meanwhile, once an iterate is successful, it must satisfy the \autoref{assumption. phi close}.

\begin{figure}[ht]
    \begin{algorithm}[H]
        \caption{Perturbation for the Hard Case }\label{alg.hardcase}
        Input: Iterate $k$, $x_k\in \mathbb{R}^n$, $\gk$, $\Hk$, $h_{k-1}, \delta_{k-1}$ where $\gk\perp \mathcal{S}_1$, tolerance $\sigma>0$\;
        \For{$i = 0, 1, \dots$}{
        Set
        \begin{equation}
            \label{eq.set phi}
            \phi_{k+i}=\gk + \frac{\varsigma_\phi}{\gamma_3^{2} h_{k+i-1}+\sigma}\lambda_1^2 v_{k}.
        \end{equation}
        Compute $I_h := [\gamma_2 \sqrt{h_{k+i-1}}, \gamma_3 \sqrt{h_{k+i-1}}]$\;
        \Repeat(\tcp*[f]{inner iterates $j$ via bisection (see \autoref{sec.bisection procedure})}){$\sqrt{h_{k+i}} \in I_h$ within tolerance $\sigma$}{
        Obtain the solution $[v_{k+i,j};t_{k+i,j}]$  of the GHM subproblem
        \begin{equation*}
            \min_{\|[v;t]\| \leq 1} \begin{bmatrix}v \\ t\end{bmatrix}^T
            \begin{bmatrix}
                H_{k}          & \phi_{k+i}     \\
                (\phi_{k+i})^T & \delta_{k+i,j}
            \end{bmatrix}
            \begin{bmatrix}v \\ t\end{bmatrix}
        \end{equation*}
        set $d_{k+i,j} = v_{k+i,j}/t_{k+i,j}$, $h_{k+i} := \left(\theta_{k+i,j}/\|d_{k+i,j}\|\right)^2$ \;
        Update $\delta_{k+i,j}$, increase $j = j+1$
        }
        $d_{k+i} = d_{k+i,j}$; $\delta_{k+i}=\delta_{k+i,j}$ \;
        Compute
        \begin{equation*}
            \rho_{k,i}= \frac{f(x_k+d_{k+i})-f(\xk)}{m_k(d_{k+i})-f(x_k)}
        \end{equation*}
        \If{$\rho_{k+i} \geq \eta_1$}{\textbf{break}}
        }
    \end{algorithm}
\end{figure}

Let us first show that if $\phi_{k+i}$ is set in a perturbed fashion \eqref{eq.set phi}, once an iterate is successful, \autoref{assumption. phi close} must hold.
\begin{lemma}\label{lem.perturbedphi.succleadsto}
    {Suppose $\phi_{k+i}$ is perturbed in the following manner:
    $$\phi_{k+i}= \gk + \frac{\varsigma_\phi}{\gamma_3^{2} h_{k+i-1}+\sigma}\lambda_1^2 v_{k},$$
    then $\phi_{k+i}$ and $d_{k+i}$ satisfies
    $
        \|\phi_{k+i}-\gk\| \leq \varsigma_{\phi} \|d_{k+i}\|^2,
    $
    as a result, whenever $d_{k+i}$ is accepted,
    \autoref{assumption. phi close} holds.}
\end{lemma}
\begin{proof}
    Note that the $(k+i)$-th iteration is successful, recalling the optimal condition \eqref{eq.mod.opt}, we have
    \begin{equation}
        H_{k}+\sqrt{h_{k+i}(\delta_{k+i})}\|d_{k+i}\| I\succeq 0. \notag
    \end{equation}
    The above inequality further implies that
    \begin{equation}\label{eq.large stepsize}
        \|d_{k+i}\|^2 \geq \frac{1}{h_{k+i}(\delta_{k+i})}\lambda_1^2 \geq \frac{1}{\gamma_3^2 h_{k+i-1}(\delta_{k+i-1})+\sigma}\lambda_1^2.
    \end{equation}
    Therefore, the gap between $\phi_{k+i}$ and $g_k$ could be bounded as
    $$
        \|\phi_{k+i} - g_k\|  = \frac{\varsigma_\phi}{\gamma_3^{2} h_{k+i-1}+\sigma}\lambda_1^2 \leq \varsigma_\phi\|d_{k+i}\|^2,
    $$
    which means that $\phi_{k+i}$ satisfies the \autoref{assumption. phi close}.
\end{proof}
The rest is to show when the hard case occurs, the \autoref{alg.hardcase} eventually produces a successful iterate.

\begin{theorem}\label{thm.hardcase}
    \autoref{alg.hardcase} takes at most $\lfloor \log_{\gamma_2}\frac{\varsigma_h}{h_{k-1}(\delta_{k-1})}\rfloor +1$ iterations to obtain a successful step. Furthermore, \autoref{assumption. phi close}  must hold.
\end{theorem}

\begin{proof}
    Similar to \autoref{lem.upper bound h}, we have
    \begin{align}
        \nonumber    m_k(d_{k+i})        & -f(\xk+d_{k+i})                                                                                                  \\
                                         &
        =(\phi_{k+i}-\gk)^T d_{k+i}+ \frac{1}{2}d_{k+i}^T (\Hk - \nabla^2 f(\xk+\xi d_{k+i})) d_{k+i}                                                       \\
        \nonumber                        &
        \qquad+\frac{\sqrt{h_{k+i}(\delta_{k+i})}}{3}\|d_{k+i}\|^3                                                                                          \\
                                         & \geq -\|\phi_{k+1}-g_k\|\|d_{k+i}\|-\frac{M}{2}\|d_{k+i}\|^3+\frac{\sqrt{h_{k+i}(\delta_{k+i})}}{3}\|d_{k+i}\|^3 \\
                                         &
        \geq -\kappa_{\phi} \|d_{k+i}\|^3 -\frac{M}{2}\|d_{k+i}\|^3+\frac{\sqrt{h_{k+i}(\delta_{k+i})}}{3}\|d_{k+i}\|^3                                     \\
        \label{eq.diff bet modeland f 2} & = \left(\frac{\sqrt{h_{k+i}(\delta_{k+i})}}{3}-\frac{M}{2}-\varsigma_\phi\right)\|d_{k+i}\|^3,
    \end{align}
    where the second inequality is from \autoref{lem.perturbedphi.succleadsto}. So when $h_{k+i}$ exceeds $\varsigma_h$, the step $d_{k+i}$ must be accepted. Also, by \autoref{lem.perturbedphi.succleadsto} we know that \autoref{assumption. phi close} must hold.
\end{proof}

Now we see the choice of $\phi_k$ and the mechanism of adjusting $h_k$ helps when the hard case occurs.  It can be embedded into \autoref{alg.gadaptive} directly. In fact, \autoref{thm.hardcase} basically implies that the hard case can be escaped by finitely many iterations of \autoref{alg.hardcase}.
As a remark for \autoref{alg.gadaptive}, we emphasize that the analysis for nonconvex functions may also be adapted for convex functions with Lipschitzian Hessians. Notice that the conditions established here \autoref{lem.step and g} and \autoref{lem.fo decrease} are similar to properties for the cubic regularized Newton's method, see for example \cite[Lemma 3, Lemma 4]{nesterov_cubic_2006}. It is thus reasonable for our method to extend to convex functions and other structured functions with similar complexity guarantees. Due to space limitations, we leave this for future study.

    {
        \begin{remark}
            Allowing inexact solutions of the subproblems is a common issue for Newton-type methods \cite{cartis_adaptive_2011,curtis_trust_2017,curtis_worst-case_2022-1}.
            Regarding methods using GHMs as subproblems, this means we can only obtain approximate eigenpair, the Ritz pairs, to update the iterates. This issue has been resolved in the original HSODM \cite{zhang_homogenous_2022}.
            For the sake of conciseness, we do not include a complexity analysis for an inexact {\it adaptive} HSODM in this paper,  but provide the details in an independent technical report \citep{heTechnicalReportHomogeneous}.
            The results show that the adaptive HSODM can still achieve the same complexity as the exact version, provided that the error tolerance in the Lanczos method is sufficiently small.
        \end{remark}
    }
\section{A Homotopy HSODM}\label{sec.homotopy}

In this section, we assume the objective function satisfies the \betacon{} condition introduced in \cite{ye2017second,luenberger_linear_2021}, and propose {\it homotopy} HSODM, a specific realization of HSODF, to solve it.

\subsection{Overview of the homotopy model and HSODM}
We start with the definition as follows.
\begin{definition}[Concordant Second-Order Lipschitz]\label{def.betacon}
    We say a function $f$ is \betacon{} if
    there exists a constant $\beta > 0$ such that for any point $x \in \operatorname{dom}(f)$, we have
    \begin{equation}\label{eq.concordant inequality}
        \|\nabla f(x+d) - \nabla f(x) - \nabla^2 f(x)d\| \leq \beta \cdot d^T \nabla^2 f(x) d,
    \end{equation}
    where $d$ satisfies $\|d\| \leq C$ for some $C > 0$ and $x+d \in \operatorname{dom}(f)$.
\end{definition}
In fact, there are many functions (see \autoref{sec.concordlip}) that satisfy the above the \betacon{} property.
Historically, the \betacon{} condition is a simplification of the Scaled Lipschitz Condition (SLC) introduced for the linearly constrained convex programming  \cite{zhu1992path,kortanek1993polynomial,den1995sufficient}.  Basically, this function class widely appears in machine learning problems, especially when the Hessian matrix degenerates with highly sparse problem data.

\subsubsection{The homotopy model}
We consider the following
\emph{homotopy model} with the objective function satisfying property \eqref{eq.concordant inequality}:
\begin{equation}\label{equation.homotopy model}
    \min_{x \in \real^n} \ f(x) + \frac{\mu}{2}\|x\|^2.
\end{equation}
This model has been used by
a \emph{path-following method} for convex optimization introduced in \cite{luenberger_linear_2021} to tackle degeneracy with a series of well-behaved, strictly convex subproblems.
The homotopy model exhibits several nice properties outlined below.
\begin{lemma}[\citet{luenberger_linear_2021}]\label{lemma. homotopy property}
    Consider a twice continuously differentiable convex function $f$, assume that its value is bounded from below and that $x^*$ is the minimal $L_2$ norm solution of $f(x)$. Denote $\xmu = \arg\min \ \left\{f(x) + \frac{\mu}{2}\|x\|^2\right\}$. Then the following properties hold.
    \begin{enumerate}[(a)]
        \item $\xmu$ is unique for any given $\mu > 0$, and it forms a continuous path as $\mu$ varies;
        \item $f(\xmu)$ is an increasing function of $\mu$, and $\|\xmu\|$ is a decreasing function of $\mu$;
        \item If $\mu \rightarrow 0^+$, then $\xmu$ converges to $x^*$;
        \item If $\mu \rightarrow \infty$, then $\xmu \rightarrow 0$.
    \end{enumerate}
\end{lemma}
Since the authors in \cite{luenberger_linear_2021} did not provide a complete analysis, we give a concise proof in \autoref{sec.lem.homotopy.property}. The above lemma justifies the convergence of the trajectory $\{\xmu\}_{\mu \to 0}$ to an optimal solution $x^*$ as $\mu$ continuously decreases to $0$. Furthermore, based on the above properties, we also conclude the following property, which will be used in the \autoref{subsec.path following complexity}.
\begin{corollary}\label{cor.bounded norm of x_mu_k}
    Given a decreasing sequence $\{\mu_k\}_{k=0}^\infty$ such that $\mu_k \rightarrow 0$ as $k \rightarrow \infty$, then for any $x_{\mu_k} = \arg\min \ \left\{f(x) + \frac{\mu_k}{2}\|x\|^2\right\}$, we have
    $\|x_{\mu_k}\| \leq \|x^*\|$.
\end{corollary}
The result is directly implied by the second and third statements in \autoref{lemma. homotopy property}. Just like the interior-point method \cite{ye_interior_1997}, \citet{ye2017second} designed a method that sequentially minimizes homotopy models defined by a decreasing series $\{\mu_k\}_{k=0}^\infty$ and applied Newton's method to solve the optimal equation at each iteration $k$,
$$
    \nabla f(x) + \mu_k \cdot x = 0.
$$
In an algorithmic scheme, if the penalty parameter $\mu_k$ is decreased at a linear rate, e.g.,
$$
    \mu_{k+1} =  \rho_k \cdot \mu_k, \ 0 < \rho_k <1,
$$
then a global linear rate of convergence can be envisioned. Since it relies on solving linear systems, one motivation naturally goes to design a homotopy HSODM. We observe that a Newton's method with some $\mu$ solves the following second-order quadratic model $m^{H}$  of \eqref{equation.homotopy model} at some point $x \in \real^n$:
\begin{align}\label{eq.homotopy.quadmodel}
    \nonumber m^{H}(x,d)  - f(x) & = \nabla f(x)^Td + \frac{1}{2}d^T \nabla^2 f(x) d + \frac{\mu}{2}\|x + d\|^2                                    \\
                                 & = \left(\nabla f(x) + \mu \cdot x \right)^Td + \frac{1}{2}d^T (\nabla^2 f(x) + \mu I) d + \frac{\mu}{2}\|x\|^2.
\end{align}
In comparison to \eqref{eq.homoquadmodel}, one can see that \eqref{eq.homotopy.quadmodel} has some extra terms regarding the information of the current point $x$, apart from the corresponding function information
$\nabla f(x), \nabla^2 f(x)$.
Therefore, the OHM cannot be applied to \eqref{eq.homotopy.quadmodel}.
Thanks to the flexibility of $\delta, \phi$ in GHM, we can
homogenize \eqref{eq.homotopy.quadmodel} as follows:
\begin{equation}\label{eq.homotopy.homomodel}
    F^H(x) := \begin{bmatrix}
        \nabla^2 f(x)                 & \nabla f(x) + \mu \cdot x \\
        \nabla f(x)^T + \mu \cdot x^T & - \mu
    \end{bmatrix}.
\end{equation}
By constructing $\psi^{H}(v,t;F^{H}) =  [v;t]^TF^{H}[v;t]$, it is equivalent to \eqref{eq.homotopy.quadmodel} up to scaling provided that $\mu$ and $d := v/t$ with $ t\neq 0$:
\begin{align*}
    m^{H}(x,d) - f(x)= \frac{1}{2t^2}\psi^{H}(v, t; F^H) + \frac{\mu}{2}\left( \|x\|^2 + \|[v; t]\|^2\right).
\end{align*}
Similarly, we can solve the GHM subject to a unit ball constraint:
\begin{equation}\label{equation.homotopy inner subproblem}
    \underset{\|[v;t]\|\leq 1}{\min} ~\psi^{H}(v, t; F^H).
\end{equation}
This is essentially a symmetric eigenvalue problem as the optimal solution to \eqref{equation.homotopy inner subproblem} always attains the sphere of the unit ball, cf. \autoref{lem.strict.curv}.

\subsubsection{The homotopy HSODM}

With the specific GHM defined in \eqref{eq.homotopy.homomodel}, we are able to develop a homotopy HSODM (\autoref{alg.homotopy HSODM}) with a series of linearly decreasing $\{\mu_k\}_{k=0}^\infty$.
By an outer iteration $k$ we update $\mu_k$. Then for each one of such intermediate goals, we apply a sequence of GHMs (\autoref{alg.GHMSolver}) to compute an approximate \emph{center}. We use the notation $x_{k,j}$ to denote the iterate in \autoref{alg.homotopy HSODM}, where $k$ is the outer iteration number and $j$ counts the inner GHMs.

\begin{figure}[h]
    \begin{algorithm}[H]
        \caption{Homotopy HSODM}\label{alg.homotopy HSODM}
        \textbf{Initialization:} initial point $x_{0,0} = 0$, iteration number $k=0$, parameter $\mu_0 = 2(\beta+1)(1+\|g_{0,0}\|^2)$; \\

        \For{$k = 0, 1, \cdots, K$}{
            Compute $(x_{k, j}, \rho_k)$ = \textbf{\texttt{iACGHM}}$(x_{k, 0}, \mu_k)$; \\
            Update $\mu_{k+1} = \rho_k \cdot \mu_k$; \\
            Set $x_{k+1, 0} := x_{k, j}$;
        }
        Output $x_{K+1, 0}$.
    \end{algorithm}
\end{figure}
\begin{figure}[h]
    \begin{algorithm}[H]
        \caption{
            \small{
                Inexact Approximate Center by GHMs (\textbf{\texttt{iACGHM}})
            }
        }\label{alg.GHMSolver}
        \textbf{Input:} $x_{k, 0}$, $\mu_k$;

        \For{$j = 0, \dots, \mathcal T_k$}{
        \If{$\|g_{k, j} + \mu_k \cdot x_{k, j}\| \leq \frac{\mu_k}{1+3(\beta+1)}$
        }{
        \label{line.approximatecenter}
        Compute $\rho_k = \frac{3(\beta+1)(1 + \|x_{k,j}\|)}{1 + 3(\beta+1)(1 + \|x_{k,j}\|)}$\;
        Return $(x_{k, j}, \rho_k)$\;
        }

        \Else{
            Obtain the solution $[v_{k,j}; t_{k,j}]$ of the GHM subproblem
            \begin{equation}\label{eq.homotopy GHM T k}
                \min_{\|[v;t]\| \leq 1} \begin{bmatrix}v \\ t\end{bmatrix}^T
                \begin{bmatrix}
                    H_{k, j}                            & g_{k, j} + \mu_k \cdot x_{k, j} \\
                    (g_{k, j} + \mu_k \cdot x_{k, j})^T & -\mu_k
                \end{bmatrix}
                \begin{bmatrix}v \\ t\end{bmatrix};
            \end{equation}
            Set $d_{k, j} = v_{k, j} / t_{k, j}$ and update $x_{T, j+1} = x_{k, j} + d_{k, j}$\;
        }
        }
    \end{algorithm}
\end{figure}

It is worth noting that for each $\mu_k$, \autoref{alg.GHMSolver} stops if the iterate satisfies the \textit{approximate centering condition} gauged by $\mu_k$ and the \betacon{} constant $\beta$ (\autoref{line.approximatecenter} in \autoref{alg.GHMSolver}).  We show that \autoref{alg.GHMSolver} has a quadratic rate of convergence, and moreover it needs at most 2 GHMs for each $\mu_k$. Besides, as an analog to the interior point method, every last inner iterate $x_{k,j}$ associated with $\mu_k$ also lies in the neighborhood of $x_{\mu_k}$. Furthermore, the width of the neighborhood becomes narrow as $\mu_k \rightarrow 0$, and each point in the neighborhood
has a \emph{fixed} deviation bound
with respect to $x_{\mu_k}$.

\subsection{Convergence analysis}\label{subsec.path following complexity}

In this subsection, we turn to \autoref{alg.homotopy HSODM} and its convergence properties. We show that if $\mu_k$ decreases geometrically, it leads to a linear rate of convergence. Besides, for each $\mu_k$, the corresponding homotopy model can be solved via finite many specialized GHMs \eqref{eq.homotopy.homomodel}. In contrast to the variant \autoref{alg.gadaptive} for nonconvex and convex functions with second-order Lipschitz continuity, \autoref{assumption. phi close} is not needed for \autoref{alg.homotopy HSODM}.

We first give the following results on the approximate centering condition.
\begin{lemma}\label{lemma.width of central path}
    For a convex function $f$, suppose the iterate $x_{k,j}$ satisfies the approximate centering condition, i.e.,
    $
        \|g_{k, j} + \mu_k \cdot x_{k, j}\| \leq \frac{\mu_k}{1+3(\beta+1)}.
    $
    Then we have
    $$
        \|x_{k,j} - x_{\mu_k}\| \leq \frac{1}{1+3(\beta+1)},
    $$
    where $x_{\mu_k} = \arg\min \ \left\{f(x) + \frac{\mu_k}{2}\|x\|^2\right\}$.
\end{lemma}
\begin{proof}
    First note that the function $f(x) + \frac{\mu_k}{2}\|x\|^2$ is $\mu_k$-strongly convex, by the Theorem 2.1.10 in \cite{nesterov_lectures_2018}, we obtain that
    $$
        \mu_k\|x_{k, j} - x_{\mu_k}\| \leq \|g_{k, j} + \mu_k \cdot x_{k, j}\|.
    $$
    Together with the approximate centering condition, it implies that
    $$
        \|x_{k, j} - x_{\mu_k}\| \leq \frac{1}{1+3(\beta+1)}.
    $$
\end{proof}

Next, we give the basic characteristics of the specialized GHM \eqref{eq.homotopy GHM T k}.
In particular, statement (a) in the lemma below excludes the possibility of the hard case.

\begin{lemma}\label{lemma.homotopy variant optimal condition}
    Suppose $f$ satisfies the \betacon{} condition. For $\mu_k > 0$, let $([v_{k,j};t_{k,j}], -\theta_{k,j})$ be the optimal primal-dual solution of the GHM \eqref{eq.homotopy GHM T k}, then the following holds.
    \begin{enumerate}[(a)]
        \item $t_{k,j} \neq 0, \theta_{k,j} > 0$;
        \item $\theta_{k,j} - \mu_k \leq \|g_{k,j} + \mu_k \cdot x_{k,j}\|$;
        \item set $d_{k,j} = v_{k,j}/t_{k,j}$, and $d_k$ satisfies
              \begin{align*}
                  \|d_{k,j}\| \leq \frac{\|g_{k,j} + \mu_k \cdot x_{k,j}\|}{\mu_k}, \quad d_{k,j}^T H_{k,j} d_{k,j} \leq \frac{\|g_{k,j} + \mu_k \cdot x_{k,j}\|^2}{\mu_k}.
              \end{align*}
    \end{enumerate}
\end{lemma}
\begin{proof}
    For the statement (a), recall that in the homotopy HSODM with GHM \eqref{equation.homotopy inner subproblem}, the second diagonal term is set as $\delta_{k,j} = -\mu_k < 0$. By the convexity of the function $f(\cdot)$ and applying the result of \autoref{lemma.ordering of tilde alpha 1}, we have $\delta_{k,j} < \tilde \alpha_1$ and hence $t_{k,j} \neq 0$. By \autoref{lem.strict.curv}, the ball constraint is always active, implying $\lambda_1(F_{k,j}) < 0$ and thus $\theta_{k,j} > 0$.

    Statement (b) is a direct consequence of \autoref{lemma.uppertheta}. Finally for (c), by \eqref{eq.homoeig.soc}, $d_k$ satisfies
    $$
        H_{k,j} d_{k,j} + \theta_{k,j} d_{k,j} = -g_{k,j} - \mu_k \cdot x_{k,j},
    $$
    and it ensures that
    $$
        \|H_{k,j}d_{k,j}\|^2 + \theta_{k,j}^2 \|d_{k,j}\|^2 + 2\theta_{k,j} d_{k,j}^T H_{k,j} d_{k,j} = \|g_{k,j} + \mu_k \cdot x_{k,j}\|^2.
    $$
    Since the cross term $d_{k,j}^T H_{k,j} d_{k,j} \geq 0$ and $\theta_{k,j} \geq -\delta_{k,j} = \mu_k > 0$, we obtain that
    $$
        \|H_{k,j} d_{k,j}\| \leq \|g_{k,j} + \mu_k \cdot x_{k,j}\| \ \text{and} \ \theta_{k,j} \|d_{k,j}\| \leq \|g_{k,j} + \mu_k \cdot x_{k,j}\|.
    $$
    Therefore, the norm of $d_{k,j}$ could be upper bounded by
    $$
        \|d_{k,j}\| \leq \frac{\|g_{k,j} + \mu_k \cdot x_{k,j}\|}{\theta_{k,j}} \leq \frac{\|g_{k,j} + \mu_k \cdot x_{k,j}\|}{\mu_k},
    $$
    and it further implies that the cross term $d_{k,j}^T H_{k,j} d_{k,j}$ satisfies
    $$
        d_{k,j}^T H_{k,j} d_{k,j} \leq \|d_{k,j}\| \cdot \|H_{k,j} d_{k,j}\| \leq \frac{\|g_{k,j} + \mu_k \cdot x_{k,j}\|^2}{\mu_k}.
    $$
    This completes the proof.
\end{proof}

In the following two lemmas, we demonstrate that for every fixed $\mu_k> 0$, the corresponding inner problem converges quadratically.  This requires a separate discussion on the initial $\mu_0$ and the following $\mu_k$, $k \ge 1$. Note that at the first iteration, we initially set $x_{0,0} = 0$. By contrast, the rest of iterations will \emph{warm start} at the previous one, i.e., $x_{k+1, 0}:= x_{k, j}$ once $x_{k, j}$ satisfies the approximate centering condition in \autoref{alg.GHMSolver}. In both two cases, the quadratic rate of convergence further leads to finite convergence of the problem in the inner loop.

For the initial iteration, we need to select a proper initial $\mu_0$ to establish the quadratic convergence. From part (c) in \autoref{lemma. homotopy property}, $x^*$ is sufficiently close to 0 as $\mu_k$ tends to infinity.
The following lemma provides guidance on how an initial $\mu_0$ should be chosen.
\begin{lemma}[Quadratic convergence with $\mu_0$]\label{lemma.homotopy initial quadratic convergence}
    Let $x_{0,0} = 0$, and the sequence $\{x_{0,j}\}$ is updated by
    $$
        x_{0, j+1} = x_{0, j} + d_{0, j},~ d_{0, j} = v_{0, j}/t_{0, j}
    $$
    where $[v_{0, j};t_{0, j}]$ solves \eqref{eq.homotopy GHM T k} at iterate $x_{0, j}$.
    If $\mu_0 \geq 2(\beta+1) \cdot \max\left\{1, \|g_{0,0}\|^2\right\}$, then
    the residual error $e_{0,j} = \|g_{0,j} + \mu_0 \cdot x_{0,j}\|$ converges quadratically, that is
    $$
        e_{0,1} \leq \frac{1}{2}, ~ \text{and} \ e_{0,j+1} \leq e_{0,j}^2,~ \forall k \geq 1.
    $$
\end{lemma}
\begin{proof}
    We first prove that $e_{0,1} < 1$. Note that $x_{0,1} = x_{0,0} + d_{0,0} = d_{0,0}$, and due to the \autoref{lemma.homotopy variant optimal condition}, it follows that
    \begin{equation}\label{eq.mu0,fromlem}
        e_{0,0} = \|g_{0,0}\|, 0<\theta_{0,0} - \mu_0 \leq e_{0,0}, \ \|d_{0,0}\| \leq e_{0,0}/\mu_0, \ \text{and} \ d_{0,0}^T H_{0,0} d_{0,0} \leq e_{0,0}^2/\mu_0.
    \end{equation}
    Therefore, in view of the definition of $e_{0,j}$, we see that
    \begin{align}
        \nonumber     e_{0,1} & = \|g_{0,1} + \mu_0 \cdot x_{0,1}\| = \|g_{0,1} + \mu_0 \cdot d_{0,0}\|                                 \\
        \nonumber             & = \|g_{0,1} - g_{0,0} - H_{0,0} d_{0,0} + g_{0,0} + H_{0,0} d_{0,0} + \mu_0 \cdot d_{0,0}\|             \\
        \nonumber             & \leq \|g_{0,1} - g_{0,0} - H_{0,0} d_{0,0}\| + \|g_{0,0} + H_{0,0} d_{0,0} + \mu_0 \cdot d_{0,0}\|      \\
        \nonumber             & \leq \beta \cdot d_{0,0}^T H_{0,0} d_{0,0} + \left|\theta_{0,0} - \mu_0\right| \cdot \|d_{0,0}\|        \\
        \label{eq. bound e1}  & \stackrel{\eqref{eq.mu0,fromlem}}{\le}  (\beta + 1) \cdot \frac{e_{0,0}^2}{\mu_0} \leq \frac{1}{2} < 1,
    \end{align}
    where the last inequality comes from $\mu_0 \geq 2(\beta+1) \cdot \max\left\{1, \|g_{0,0}\|^2\right\}$. As for the second part of this lemma, we could directly obtain the result by using a similar argument,
    \begin{align}
        e_{0,j+1} & = \|g_{0,j+1} + \mu_0 \cdot x_{0,j+1}\|\leq (\beta + 1) \cdot \frac{e_{0,j}^2}{\mu_0} \leq e_{0,j}^2,
    \end{align}
    and therefore, the quadratic speed of convergence is established.
\end{proof}

\begin{corollary}\label{cor.number of iteration in initial epoch}
    In the initial iteration $k=0$, the number of iterates $j$ in \textbf{\texttt{iACGHM}} is upper bound by
    $$
        \mathcal T_0 = \left \lceil
        \log_2 \left(\frac{\max\left\{\log(1+3(1+\beta)) - \log \mu_0, \log2\right\}}{\log 2} \right) \right \rceil + 1 \le 2.
    $$
\end{corollary}
\begin{proof}
    Since $e_{0,1} \leq \frac{1}{2}$ and $e_{0,j+1} \leq e_{0,j}^2$, by setting $e_{0,j} \leq \frac{\mu_0}{1+3(1+\beta)}$, we could conclude that
    $$
        j \ge \log_2 \left(\frac{\max\left\{\log(1+3(1+\beta)) - \log \mu_0, \log2\right\}}{\log 2} \right) + 1,
    $$
    and so that
    $$\mathcal T_0 = \left \lceil
        \log_2 \left(\frac{\max\left\{\log(1+3(1+\beta)) - \log \mu_0, \log2\right\}}{\log 2} \right) \right \rceil + 1.
    $$
    Since
    $$
        \mu_0 \ge 2(\beta+1) \cdot \max\left\{1, \|g_{0,0}\|^2\right\} \ge 2(\beta+1),
    $$
    we have
    \begin{align*}
        \mathcal T_0 & \le \left \lceil
        \max\left\{0, \log_2\left(\frac{\log(1+3(1+\beta)) - \log (2(\beta+1))}{\log 2} \right)\right\} \right \rceil + 1
        \\ &\le \left \lceil
        \log_2\left(\frac{\log(5(1+\beta)) - \log (2(\beta+1))}{\log 2} \right)\right \rceil + 1 = \left \lceil
        \log_2\left(\frac{\log 2.5}{\log 2} \right)\right \rceil + 1 = 2.
    \end{align*}
    This completes the proof.
\end{proof}

Given a large enough $\mu_0$, \autoref{cor.number of iteration in initial epoch} demonstrates that we would naturally obtain an iterate $x_{0,j}$ which satisfies the approximate centering condition. For the subsequent iteration, we prove that the quadratic convergence also holds once the penalty $\mu_{k}$ is linearly decreased.
\begin{lemma}\label{lemma.homotopy inner quadratic convergence}
    For the iteration $k \ge 1$, the sequence $\{x_{k,j}\}$ is updated by
    $$
        x_{k,j+1} = x_{k,j} + d_{k,j}, ~ d_{k,j} = v_{k,j}/t_{k,j}
    $$
    where $[v_{k, j};t_{k, j}]$ solves \eqref{eq.homotopy GHM T k} at iterate $x_{k, j}$. Similarly, letting $e_{k, j} = \|g_{k, j} + \mu_{k} \cdot x_{k, j}\|$, $\frac{\beta + 1}{\mu_{k}}e_{k,j}$ converges quadraticlly, i.e.
    $$
        \frac{\beta + 1}{\mu_{k}}e_{k,0} \leq \frac{2}{3}, ~ \frac{\beta + 1}{\mu_{k}}e_{k,j+1} \leq \left(\frac{\beta + 1}{\mu_{k}}e_{k,j}\right)^2, ~ \forall k \geq 1.
    $$
\end{lemma}
\begin{proof}
    Firstly, by the mechanism of our method, for the initial point $x_{k, 0}$ in the iteration $k \geq 1$, we must have
    \begin{equation}\label{eq.ref.subseq}
        \begin{aligned}
             & \|g_{k, 0} + \mu_{k-1} \cdot x_{k, 0}\|  \leq \frac{\mu_{k-1}}{1+3(\beta+1)},                       \\
             & \rho_{k-1}                     = \frac{3(\beta+1)(1+\|x_{k, 0}\|)}{1 + 3(\beta+1)(1+\|x_{k, 0}\|)}, \\
             & \mu_{k}                  = \rho_{k-1} \cdot \mu_{k-1}.
        \end{aligned}
    \end{equation}
    Note that,
    \begin{equation}\label{eq.rhoratio}
        \frac{1-\rho_{k-1}}{\rho_{k-1}} = \frac{1}{3(\beta+1)(1+\|x_{k, 0}\|)} \\
    \end{equation}
    It follows,
    \begin{subequations}
        \begin{align}
            \nonumber   \frac{\beta + 1}{\mu_{k}}e_0  = & ~\frac{\beta + 1}{\mu_{k}} \cdot \|g_{k,0} + \mu_{k} \cdot x_{k,0}\|                                                                          \\
            \nonumber   =                               & ~\frac{\beta + 1}{\mu_{k}} \cdot \|g_{k,0} + \mu_{k-1} \cdot x_{k,0} - (1-\rho_{k-1})\mu_{k-1} \cdot x_{k,0}\|                                \\
            \nonumber   \leq                            & ~\frac{\beta + 1}{\mu_{k}} \cdot \|g_{k,0} + \mu_{k-1} \cdot x_{k,0}\| + \frac{(\beta + 1)(1-\rho_{k-1})\mu_{k-1}}{\mu_{k}} \cdot \|x_{k,0}\| \\
            \label{eq.sube1quadcase1}   \le             & ~\frac{\beta + 1}{\mu_{k}} \cdot \frac{\mu_{k-1}}{1+3(\beta+1)}   + \frac{(1-\rho_{k-1})\mu_{k-1}}{\mu_k} \cdot (\beta+1)\|x_{k,0}\|          \\
            \le                                         & ~\frac{1}{3} + \frac{\|x_{k,0}\|}{3(1+\|x_{k,0}\|)} \le \frac{2}{3},
        \end{align}
    \end{subequations}
    where \eqref{eq.sube1quadcase1} follows from \eqref{eq.rhoratio}.
    Hence the first part is established. For the second part, by the properties of GHM in \autoref{lemma.homotopy variant optimal condition}, we obtain that
    \begin{equation}\label{eq. bound of thetak dk Hdk}
        \|d_{k,j}\| \leq e_{k,j}/\mu_{k}, \ d_{k,j}^T H_{k,j} d_{k,j} \leq e_{k,j}^2/\mu_{k} , \ \text{and} \ |\theta_{k,j} - \mu_{k}| \leq e_{k,j}.
    \end{equation}
    Plugging the above expression into $\frac{\beta+1}{\mu_{k}}e_{T, j+1}$, we can upper bound the term as,
    \begin{align}
        \nonumber    \frac{\beta+1}{\mu_{k}}e_{k,j+1}
         & = \frac{\beta + 1}{\mu_{k}} \cdot \|g_{k,j+1} + \mu_{k} \cdot x_{k,j+1}\|                                                                                                                 \\
         & \leq \frac{\beta + 1}{\mu_{k}} \cdot \left(\beta \cdot d_{k,j}^T H_{k,j} d_{k,j} + |\theta_{k,j} - \mu_{k}| \cdot \|d_{k,j}\|\right)\leq \left(\frac{\beta + 1}{\mu_{k}}e_{k,j}\right)^2.
    \end{align}
    where the first inequality follows the same technique used in \eqref{eq. bound e1} and the second one is due to \eqref{eq. bound of thetak dk Hdk}.
\end{proof}

\begin{corollary}\label{cor.number of iteration in next epoch}
    In the iteration $k \geq 1$, the number of iterates $j$ in the inner iteration is upper bound by
    $$
        \mathcal T:=\mathcal T_k = \left \lceil \log_2 \left(\frac{\log(1+3(\beta+1)) - \log(\beta+1)}{\log 3 - \log 2}\right) \right \rceil \le 2.
    $$
\end{corollary}
\begin{proof}
    With a little abuse use of notation, for any $k \geq 1$, we denote $E_{k,j} = \frac{\beta+1}{\mu_{k}}e_{k,j}$ and therefore
    $
        E_{k,0} \leq \frac{2}{3}, \ \text{and} \ E_{k,j+1} \leq E_{k,j}^2.
    $
    Recall the approximate centering condition, it is sufficient to set
    $
        E_{k,j} \leq \frac{\beta+1}{1+3(\beta+1)},
    $
    which implies that
    $$
        e_{k,j} = \frac{\mu_{k}}{\beta+1} E_{k,j} \leq \frac{\mu_{k}}{1+3(\beta+1)}.
    $$
    and again we obtain
    $$
        j \geq \log_2 \left(\frac{\log(1+3(\beta+1)) - \log(\beta+1)}{\log 3 - \log 2}\right).
    $$
    And similarly,
    \begin{align*}
        \mathcal T_k & = \left \lceil \log_2 \left(\frac{\log(1+3(\beta+1)) - \log(\beta+1)}{\log 3 - \log 2}\right) \right \rceil                                                                                                                         \\
                     & \le \left \lceil \log_2 \left(\frac{\log(4(\beta+1)) - \log(\beta+1)}{\log 3 - \log 2}\right) \right\rceil \le  \left \lceil \log_2 \left(\frac{\log 4}{\log 3 - \log 2}\right) \right \rceil = \left \lceil 1.77\right \rceil = 2.
    \end{align*}
    which completes the proof.
\end{proof}
In both two cases ($k=0, k\ge 1$), it is shown that the inner loop has finite convergence. Noting that for $k\ge 1$ our estimate is uniform, we may simply let $\mathcal T$ be the number of inquiries of GHMs.

Before we prove the upper bound for the number of iterations, we present a uniform upper bound for $\rho_k$, which plays a significant role in the convergence rate of the homotopy HSODM.
\begin{lemma}\label{lemma.upper bound for rho}
    There exists a constant $\tau \in (0, 1)$ such that for all $k \geq 0$, we have
    $
        \rho_k \leq \tau.
    $
\end{lemma}
\begin{proof}
    Recall the definition of $\rho_k = \frac{3(\beta+1)(1 + \|x_{k,j}\|)}{1 + 3(\beta+1)(1 + \|x_{k,j}\|)}$, it is sufficient to find a uniform bound for $x_{k, j}$. By \autoref{lemma.width of central path} and \autoref{cor.bounded norm of x_mu_k}, we see that
    $$
        \|x_{k,j}\| \leq \|x_{k,j} - x_{\mu_k}\| + \|x_{\mu_k}\| \leq \frac{1}{1+3(\beta+1)} + \|x^*\|.
    $$
    Therefore, set
    $$
        \tau = \frac{3(\beta+1)\left(1 + \frac{1}{1+3(\beta+1)} + \|x^*\|\right)}{1 + 3(\beta+1)\left(1 + \frac{1}{1+3(\beta+1)} + \|x^*\|\right)}
    $$
    and we complete the proof.
\end{proof}

\begin{lemma}[Number of outer iterations]\label{lemma. number of iterations}
    Suppose $f$ satisfies the \betacon{} condition. For any given $\epsilon > 0$, then after at most
    $$
        K = \left\lceil \log_\tau \left(\frac{(1+3(\beta+1))\epsilon}{2(\beta+1)(1+\|\nabla f(0)\|^2)((3\beta+4)\|x^*\|+2)}\right) \right\rceil
    $$
    iterations, where $\tau$ is defined in \autoref{lemma.upper bound for rho}, the output iterate $x_{K+1, 0}$ satisfies $\|\nabla f(x_{K+1, 0})\| \leq \epsilon$.
\end{lemma}
\begin{proof}
    Note that in the \autoref{alg.homotopy HSODM}, $x_{K+1, 0} = x_{K, j}$ and satisfies approximate centering condition, combining with \autoref{cor.bounded norm of x_mu_k}, we have
    \begin{equation}\label{eq.gradient norm for last iterate}
        \begin{aligned}
            \|\nabla f(x_{K+1, 0})\|
             & = \|g_{k, j}\|\leq \|g_{k, j} + \mu_{K} \cdot x_{K, j}\| + \mu_{K}\|x_{K, j}\| \\
             & \leq \frac{\mu_{K}}{1+3(\beta+1)} + \mu_{K}\|x_{K, j}\|                        \\
             & \leq \left(\frac{2}{1+3(\beta+1)} + \|x^*\|\right) \cdot \mu_{K}.
        \end{aligned}
    \end{equation}

    With \autoref{lemma.upper bound for rho} in place, $\mu_{K}$ could be upper bounded as below,
    \begin{equation}
        \mu_{K} = \rho_{K-1} \cdot \mu_{K-1} \leq \tau \cdot \mu_{K-1} \leq \tau^{K} \cdot \mu_{0}.
    \end{equation}
    Plugging the above inequality into \eqref{eq.gradient norm for last iterate}, it follows that
    $$
        \|\nabla f(x_{K+1, 0})\| \leq \left(\frac{2}{1+3(\beta+1)} + \|x^*\|\right) \cdot \tau^{K} \cdot \mu_{0} \leq \epsilon,
    $$
    where the last inequality comes from $K \ge \log_\tau \left(\frac{(1+3(\beta+1))\epsilon}{2(\beta+1)(1+\|\nabla f(0)\|^2)((3\beta+4)\|x^*\|+2)}\right)$ and $\mu_0 = 2(\beta+1)(1+\|\nabla f(0)\|^2)$, which completes the proof.
\end{proof}

As before, we establish the total number of GHMs soloved in the homotopy HSODM.
\begin{theorem}\label{thm.comlexity for homotopy HSODM}
    Suppose $f$ satisfies the \betacon{} condition. For any given $\epsilon > 0$, we let $\mathcal K_{\psi}$ be the total number of needed GHMs in \autoref{alg.homotopy HSODM} to return an approximate global minimum.
    Then,
    $$
        \mathcal K_{\psi} = \left\lceil 2\log_\tau \left(\frac{(1+3(\beta+1))\epsilon}{2(\beta+1)(1+\|\nabla f(0)\|^2)((3\beta+4)\|x^*\|+2)}\right) \right\rceil.
    $$
\end{theorem}
\begin{proof}
    The result is directly implied by \autoref{lemma. number of iterations}, \autoref{cor.number of iteration in initial epoch} and \autoref{cor.number of iteration in next epoch}.
\end{proof}

Generally, the value of concordant second-order Lipschitz parameter $\beta$ is unknown in practice. However, thanks to the mechanism of \autoref{alg.homotopy HSODM}, it is safe to substitute $\beta$ with a larger value while maintaining linear convergence (see \autoref{lemma.upper bound for rho} for more details). Therefore, a straightforward heuristic involves computing an approximate estimator at a given point $x$ with a small perturbation $d$:
\begin{equation*}
    \frac{\|\nabla f(x+d) - \nabla f(x) - \nabla^2 f(x)d\|}{d^T \nabla^2 f(x) d},
\end{equation*}
and constructing an upper bound based on it.

We conclude the discussion of this section with the following remark.
\begin{remark}
    We emphasize that we \emph{do not assume} the sublevel set is bounded, which is widely used in the complexity analysis of convex optimization \cite{hanzely2022damped,nesterov2008accelerating,nesterov_lectures_2018}. Furthermore, the \autoref{alg.homotopy HSODM} finds a stationary point, which further implies convergence of function values. Namely, for the convex function, it holds that
    $$
        f(x) - f^* \leq \|\nabla f(x)\| \cdot \|x - x^*\|.
    $$
    However, conversely given $f(x) - f^* \le \epsilon$, how to make gradients small is nontrivial, see \cite{nesterov2012make,foster2019complexity}. Finally, it is worth noting that the \betacon{} functions are not necessarily strongly convex. Nevertheless,  global linear convergence can be established with the homotopy HSODM.
\end{remark}

\section{Numerical Experiments}\label{sec.numeric}

We implemented our algorithms in Julia on a Mac OS desktop with a 3.2 GHz 6-Core Intel Core i7 processor and 32 GB DDR4 memory. Most of the subroutines can be found in standard Julia packages. For example, we use the line search algorithms in Optim.jl package. Access to the CUTEst benchmark is supported by CUTEst.jl. KrylovKit.jl is used for the conjugate gradient method as well as the Lanczos method. The details of the implementation can be found at \texttt{github.com/bzhangcw/DRSOM.jl}. Briefly, two methods in this paper, \ahsodm{} (\autoref{alg.gadaptive}), \pfhsodm{} (\autoref{alg.homotopy HSODM}), and the inexact Newton method (\inewton{}) are implemented by ourselves. {We set a tolerance to $\min\{10^{-5}, 10^{-2}\|\gk\|\}$ in the Lanczos method when solving GHMs}.

\subsection{CUTEst benchmark}
We test on a subset of CUTEst problems as benchmarking results for nonconvex optimization. We select the problems whose dimension of the decision variable satisfies $500 \le n \le 5000$, which gives 81 instances. We compare to Cubic Regularized Newton Method (\arc{}) and Newton Trust-Region Method (\newtontrst{}) on these problems. The \newtontrst{} uses the Steihaug-Toint conjugate gradient method when solving trust-region subproblems.   We use the state-of-the-art Julia implementation of these two competing algorithms in \cite{dussaultScalableAdaptiveCubic2023} and the default settings therein.

\begin{figure}[ht]
    \subcaptionbox{Profile of iterations}[0.5\textwidth]{
        \includegraphics[width=0.99\linewidth]{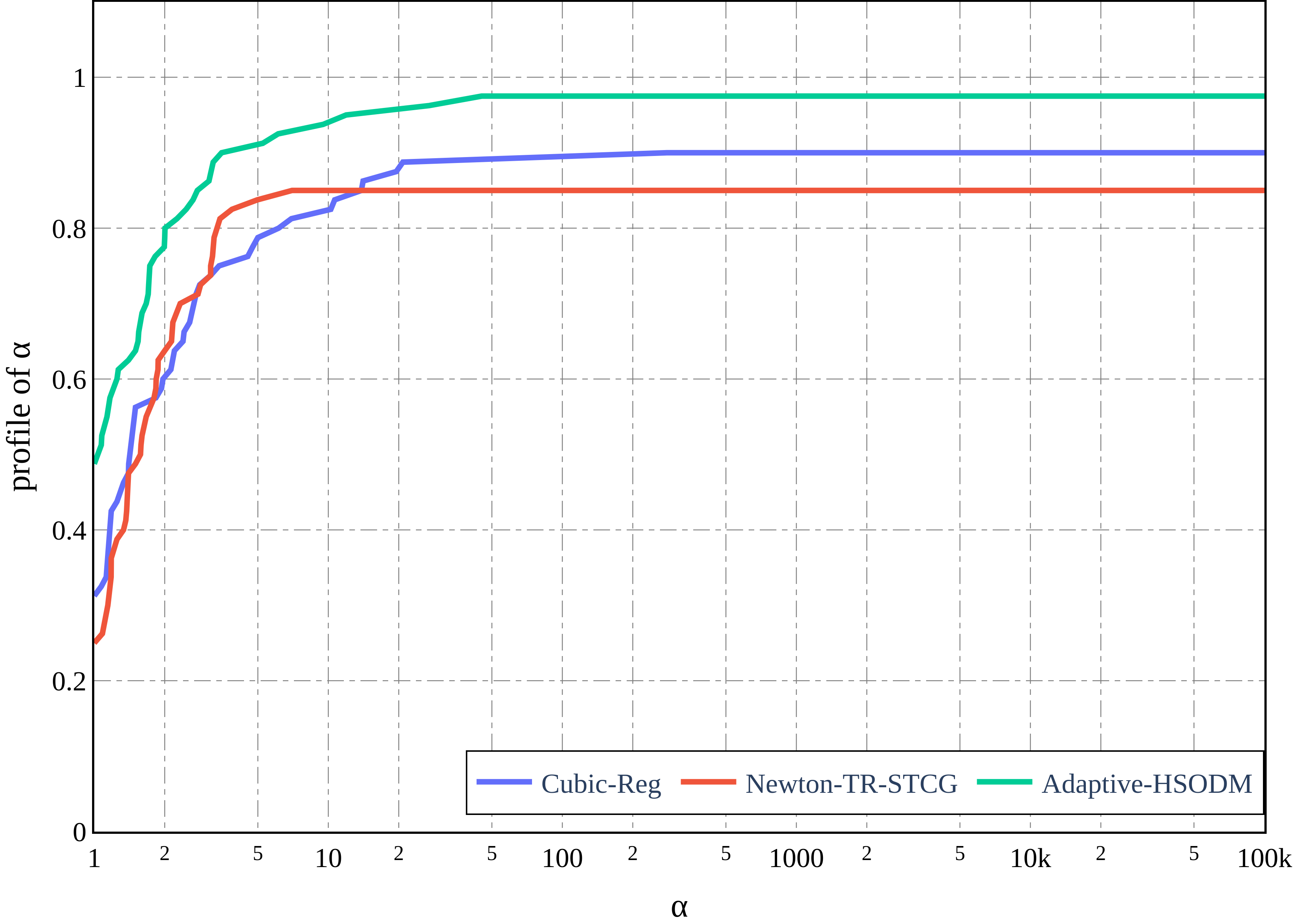}
    }~
    \subcaptionbox{Profile of gradient evaluations}[0.5\textwidth]{
        \includegraphics[width=0.99\linewidth]{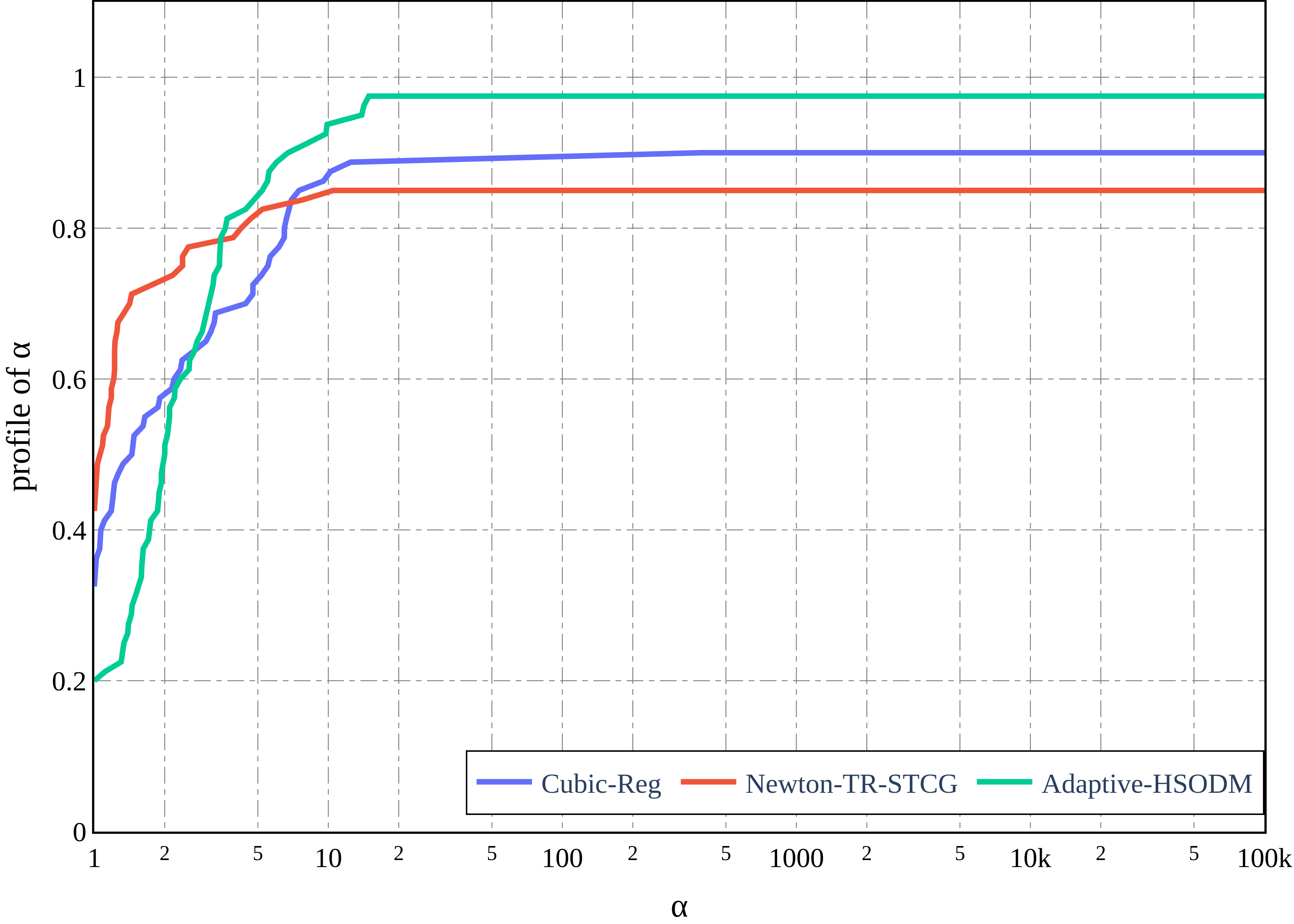}
    }
    \caption{Performance profile on CUTEst problems.}\label{fig.profile}
\end{figure}

We set a termination criterion to find an iterate $\xk$ such that $\|\gk\| \le 10^{-5}$.
The time limit is set to 200 seconds per instance.
If an instance fails, its iteration number and solving time are set to $20,000$. \autoref{tab.perf.geocutest} present a statistical summary of the three algorithms.
    {
        In this table, $\mathcal K$ means the number of successful instances; $\overline t_{G}, \overline k_{G}$ are scaled geometric means (SGM, scaled by 1 second and 50 iterations, respectively) for time and iterations; $\overline k_{G}^f, \overline k_{G}^g$ are SGMs for function and gradient evaluations, respectively, where each Hessian-vector product is counted as two gradient evaluations.
        A closer look at the performance profiles is visualized in \autoref{fig.profile}.
    }
Our implementation of \ahsodm{} seems to be quite robust: it has the best number of successful instances, iteration number, and running time. The three methods are close in function evaluations, and \newtontrst{} has the best performance in gradient evaluations. \ahsodm{} seems to use more gradient evaluations per iteration ($\overline k_G^g/\overline k_G$). This is perhaps due to a higher dependency on the problem dimension.
\begin{table}[ht]
    \centering
    \caption{A summary of CUTEst results (81 instances).}\label{tab.perf.geocutest}
    \begin{tabular}{lrrrrrrl}
        \toprule
        method      & $\mathcal K$ & $\overline t_G$ & $\overline k_G$ & $\overline k_G^f$ & $\overline k_G^g$ \\
        \midrule
        \arc        & 72.00        & 14.20           & 304.94          & 304.94            & 1810.82           \\
        \ahsodm     & 78.00        & 10.49           & 189.70          & 323.47            & 1642.13           \\
        \newtontrst & 68.00        & 21.80           & 353.21          & 353.21            & 1327.13           \\
        \bottomrule
    \end{tabular}
\end{table}
{
Different from the convex (near-generate) case, it seems that the Lanczos method should be modified to adopt better inexactness strategies since the gap-dependent conditioning is missing in the general case (\autoref{sec.warmup}). We leave this improvement for future study.
In all, the results show that our preliminary implementation of \ahsodm{} is competitive with the state-of-the-art packages \cite{dussaultScalableAdaptiveCubic2023} on CUTEst problems.
}

\subsection{$L_2$-regularized logistic regression}
In this subsection, we provide preliminary numerical results of solving $L_2$-regularized logistic regression \eqref{eq.logistic regression}.  We present the results of two degenerate high-dimensional datasets in the LIBSVM library: \texttt{rcv1} and \texttt{news20}.
We randomly start the problem from a normally distributed vector:
$$x_0 \sim \mathbf N(0, 100 \cdot I_n)$$
so that the initial point is far away from a local solution. We terminate the algorithms at iterate $\xk$
when $\|\gk\| \le 10^{-8}$.

Since the problem is convex and highly degenerate, we compare it with the inexact Regularized Newton's method. We use the gradient norm as regularization, i.e., the \inewton{} computes at each iteration
\begin{equation}\label{eq.reg.logistic.eq}
    \left( \Hk + \sigma \|\gk\|^{1/2} I \right) \dk = -\gk.
\end{equation}
For simplicity we try \emph{fixed} regularization: $\sigma \in \{10^{-5},10^{-4},10^{-3}\}$. Due to the fact that $n$ is large and operations on the Hessian matrices may be costly, the conjugate gradient method (CG) is to solve linear systems in \inewton{}, we stop CG if the iterate satisfies
\begin{equation*}
    \left\| \left( \Hk + \sigma \|\gk\|^{1/2} I \right) \dk + \gk \right\| \le \min\{10^{-4}, \zeta \|\gk\|\}, \zeta \approx \Theta(10^3),
\end{equation*}
motivated from \cite{curtis_trust-region_2021,royer_complexity_2018}. The tolerance relative to $\|\gk\|$ is active when it finalizes at a high accuracy. In this case, we slightly tighten the tolerance in case the algorithm gets stuck; otherwise, the looser precision is used.
For the eigenvalue problems in GHMs, we use same tolerance policy in the Lanczos method. All methods are facilitated with a backtrack line-search algorithm.

\begin{figure}[ht]
    \centering
    \footnotesize
    \begin{subfigure}{.49\textwidth}
        \centering
        \includegraphics[width=0.99\linewidth]{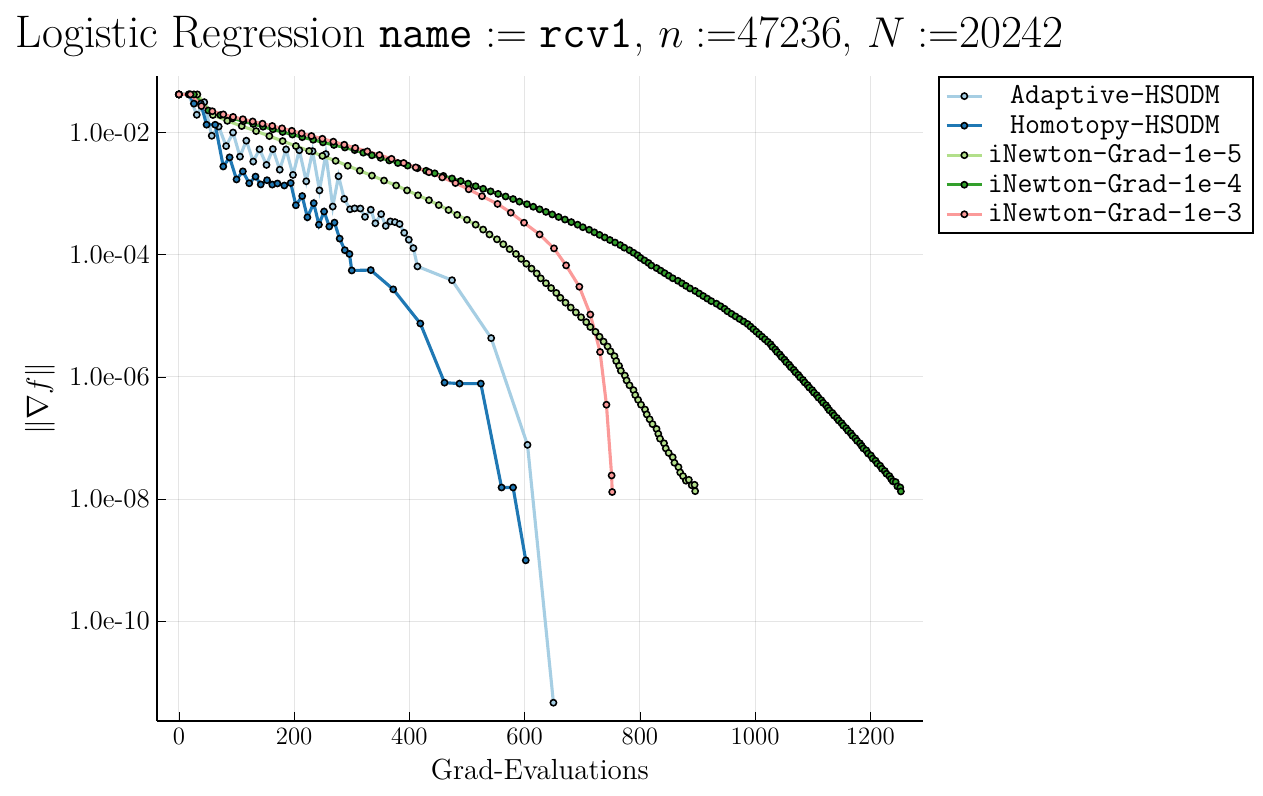}
        \caption{Results of \texttt{rcv1} using $\gamma=10^{-5}$}
    \end{subfigure}
    ~
    \begin{subfigure}{.49\textwidth}
        \includegraphics[width=0.99\linewidth]{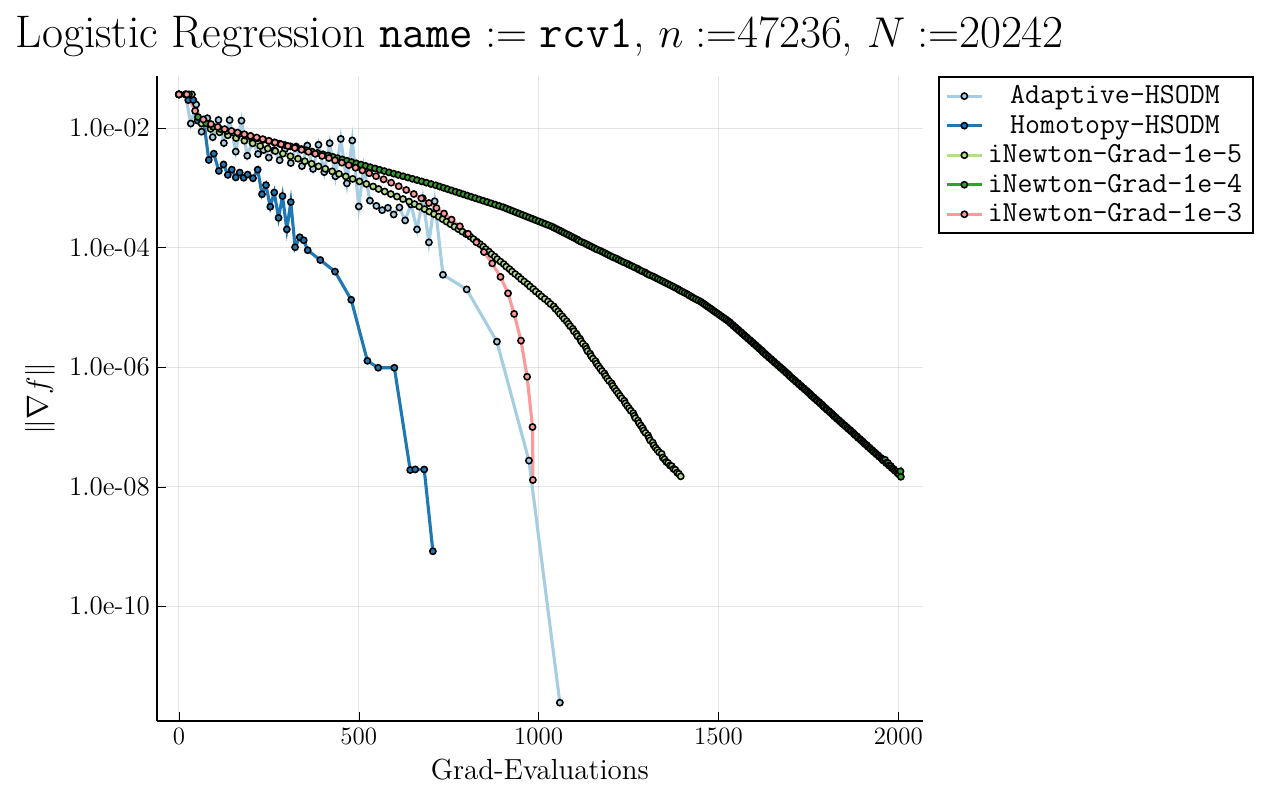}
        \caption{Results of \texttt{rcv1} using $\gamma=10^{-6}$}
    \end{subfigure}

    \begin{subfigure}{.49\textwidth}
        \includegraphics[width=0.99\linewidth]{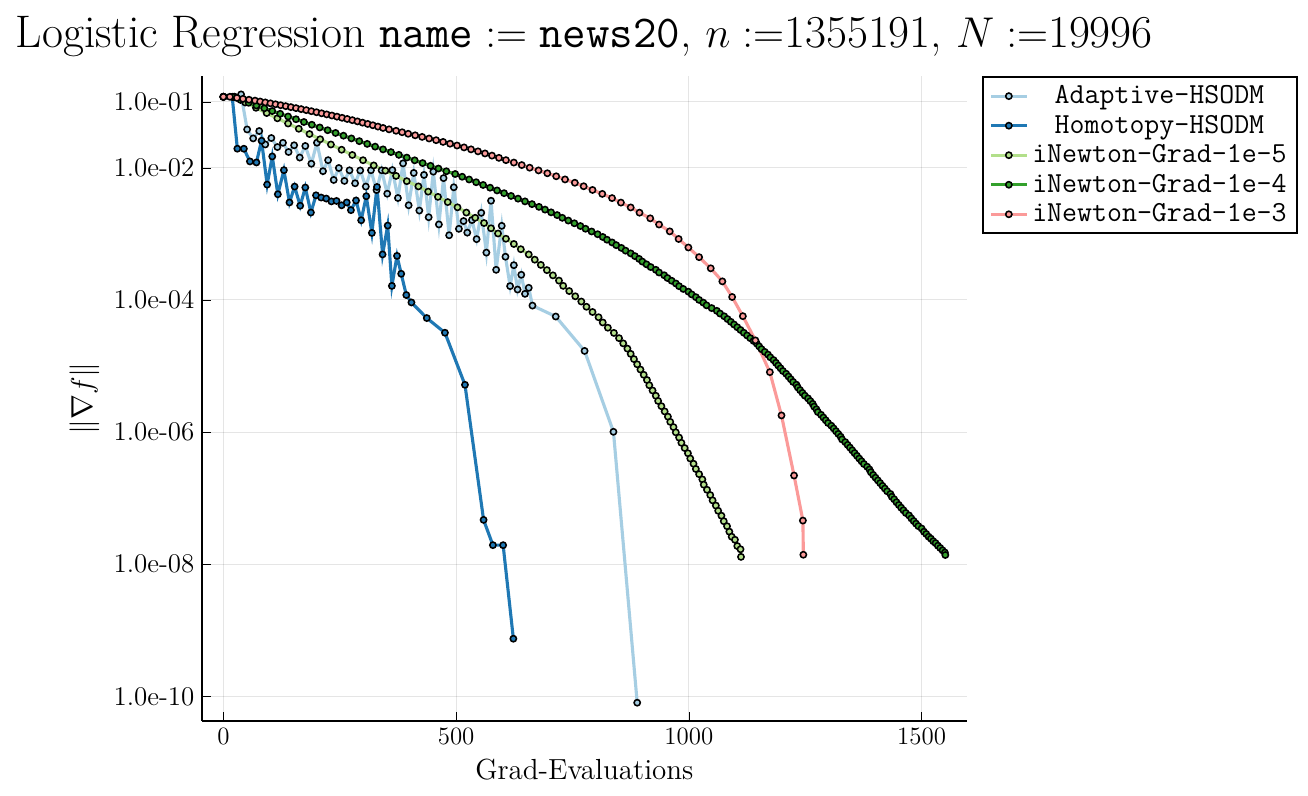}
        \caption{Results of \texttt{news20} using $\gamma=10^{-5}$}
    \end{subfigure}
    ~
    \begin{subfigure}{.49\textwidth}
        \includegraphics[width=0.99\linewidth]{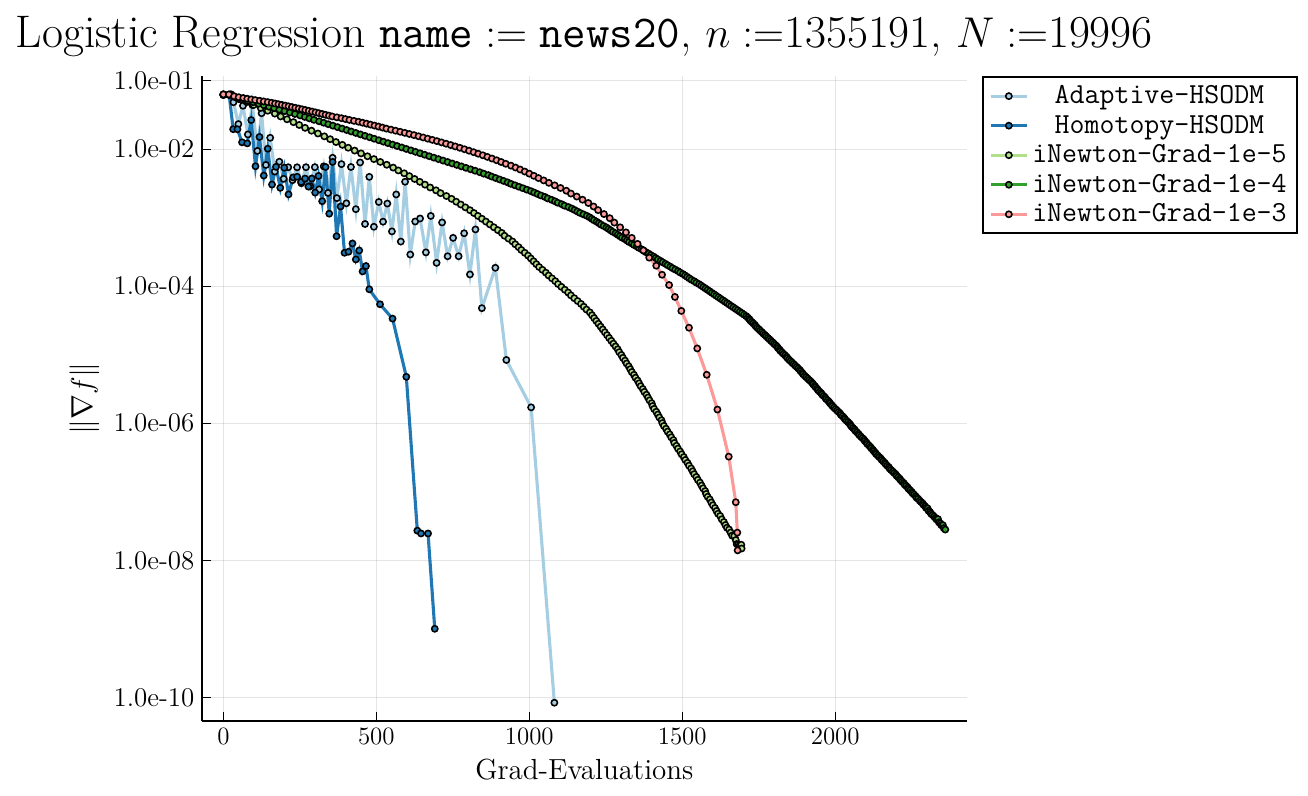}
        \caption{Results of \texttt{news20} using $\gamma=10^{-6}$}
    \end{subfigure}
    \normalsize
    \caption{Performance of a set of SOMs on $L_2$ regularized logistic regression with dataset \texttt{rcv1} and \texttt{news20}. }
    \label{fig.perfprof}
\end{figure}

In \autoref{fig.perfprof}, we illustrate the trajectory of gradient norm versus the number of gradient evaluations. Similarly, we count each Hessian-vector product as two gradient evaluations.  The result show that \pfhsodm{}, \ahsodm{} and \inewton{} are close if the regularization $\gamma$ is sufficiently large and $\sigma$ in \eqref{eq.reg.logistic.eq} is selected properly.  If the problem becomes more degenerate (smaller $\gamma$), all methods become slower, while the \pfhsodm{} seems to be the most resilient. One piece of evidence for this finding is perhaps due to the fact that \pfhsodm{} utilizes the concordance condition \eqref{eq.concordant inequality} instead of the usual Lipschitzian properties.

\paragraph{The benefits of warm-start}
In the \pfhsodm{}, we show that the iterate is continuous in some sense. Since we use the Lanczos method (\glanczos{}) to solve the GHMs, it becomes natural to use the previous solution $[v_{k-1}; t_{k-1}]$ to \emph{warm start} the current eigenvalue problem at the iterate $k$, that is, the Lanczos method use $[v_{k-1}; t_{k-1}]$ as the initial point.

\begin{figure}[h]
    \centering
    \footnotesize
    \begin{subfigure}{.46\textwidth}
        \centering
        \includegraphics[width=0.99\linewidth]{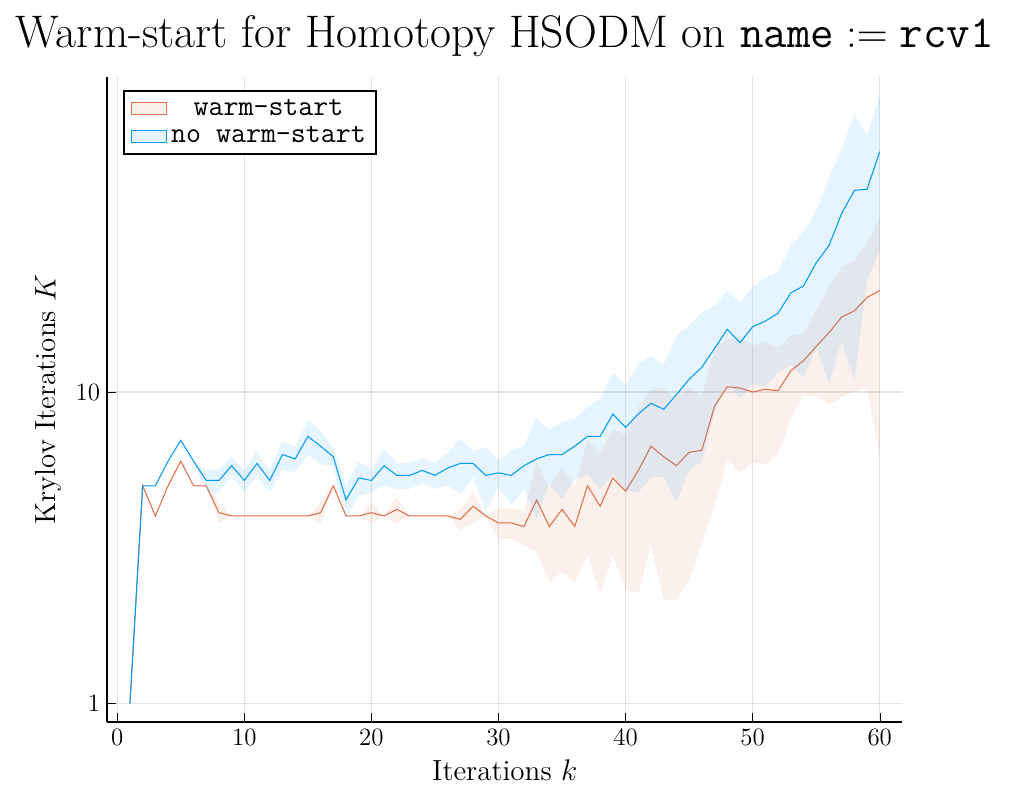}
    \end{subfigure}
    ~
    \begin{subfigure}{.48\textwidth}
        \centering
        \includegraphics[width=0.99\linewidth]{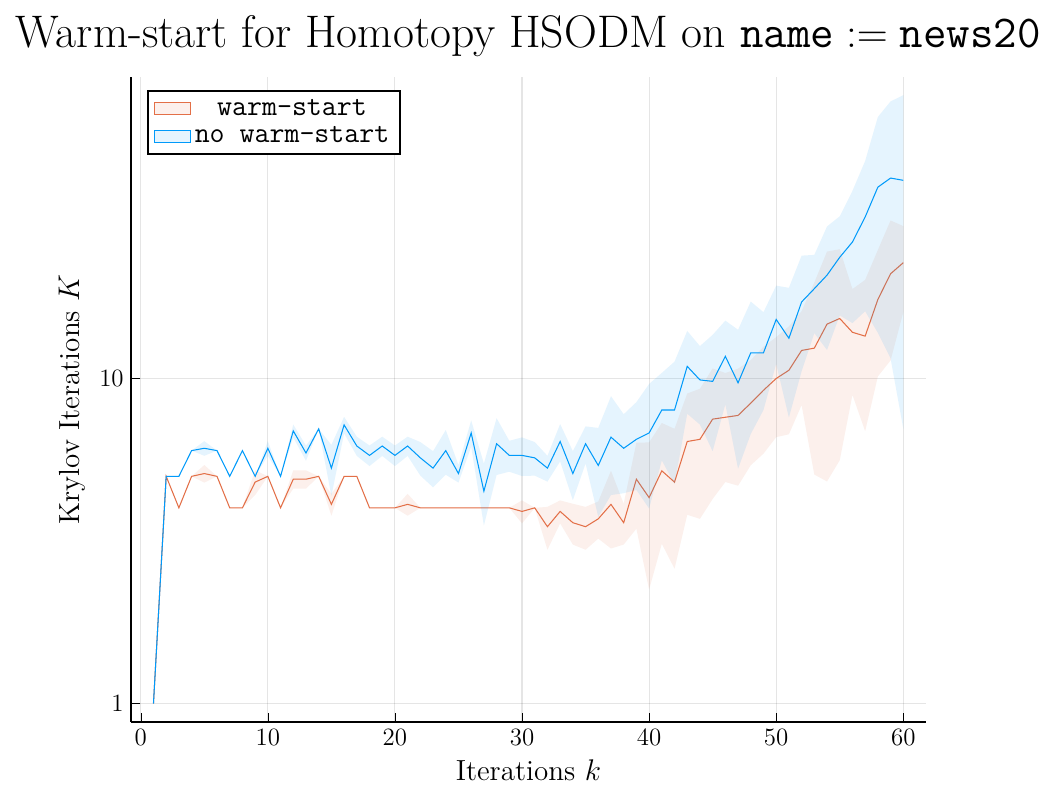}
    \end{subfigure}
    \caption{Performance of warm-starting \pfhsodm{} on $L_2$ regularized logistic regression with dataset \texttt{rcv1} and \texttt{news20}. }
    \label{fig.warmstart}
\end{figure}

We conduct a preliminary result in this direction. For the same set of problems, we compare \pfhsodm{} without and with a warm-starting vector from the last iterate. Similarly, we compare the number of  Krylov iterations $K$ needed to solve the GHM to the desired accuracy $10^{-8}$.
\autoref{fig.warmstart} presents the results on dataset \texttt{rcv1} and \texttt{news20} along the main iterates $\xk$ for every $k$.
Our preliminary results show that \pfhsodm{} can clearly save the number of Krylov iterations from using previous eigenvectors. This proposal is left for future study.

    {
        \section{Conclusions}
        In this paper, we introduce a homogeneous second-order descent framework (HSODF) for unconstrained optimization. We demonstrate that the condition number of GHMs, $\kappa_{\mathrm{L}}(F_k)$, is always bounded from above, providing an advantage in handling degenerate problems. Notably, common second-order methods can be recovered by solving the subproblem via a sequence of GHMs and an auxiliary bisection procedure. Within the proposed framework, we develop serveral variants of HSODM.
        In particular,
        we develop an adaptive HSODM for functions with a second-order continuous Hessian, extending the line-search version from \cite{zhang_homogenous_2022}. For the homotopy HSODM, parameter adjustments are unnecessary, and only 2 GHMs are required per iteration.
        Our next goal is to extend the homogeneous framework to interior-point methods, especially for nonconvex conic optimization.
    }
\section*{Declarations}


\paragraph{Conflict of interest}
The authors have no competing interests to declare that are relevant to the content of this article.

\addcontentsline{toc}{section}{References}

\bibliographystyle{plainnat}
\bibliography{homo}

@article{cartisCubicquarticRegularizationModels2025,
  title   = {Cubic-Quartic Regularization Models for Solving Polynomial Subproblems in Third-Order Tensor Methods},
  author  = {Cartis, Coralia and Zhu, Wenqi},
  year    = {2025},
  month   = jan,
  journal = {Mathematical Programming},
  issn    = {1436-4646},
  doi     = {10.1007/s10107-024-02176-y},
  urldate = {2025-03-09},
  langid  = {english}
}

@article{nesterovSuperfastSecondOrderMethods2021,
  title   = {Superfast {{Second-Order Methods}} for {{Unconstrained Convex Optimization}}},
  author  = {Nesterov, Yurii},
  year    = {2021},
  month   = oct,
  journal = {Journal of Optimization Theory and Applications},
  volume  = {191},
  number  = {1},
  pages   = {1--30},
  issn    = {1573-2878},
  doi     = {10.1007/s10957-021-01930-y},
  urldate = {2023-01-01},
  langid  = {english}
}

@misc{cartisSecondorderMethodsQuarticallyregularised2025,
  title         = {Second-Order Methods for Quartically-Regularised Cubic Polynomials, with Applications to High-Order Tensor Methods},
  author        = {Cartis, Coralia and Zhu, Wenqi},
  year          = {2025},
  month         = jan,
  number        = {arXiv:2308.15336},
  eprint        = {2308.15336},
  primaryclass  = {math},
  publisher     = {arXiv},
  doi           = {10.48550/arXiv.2308.15336},
  urldate       = {2025-05-12},
  archiveprefix = {arXiv}
}

@article{adachi_solving_2017,
  title   = {Solving the {Trust}-{Region} {Subproblem} {By} a
             {Generalized} {Eigenvalue} {Problem}},
  volume  = {27},
  number  = {1},
  journal = {SIAM Journal on Optimization},
  author  = {Adachi, Satoru and Iwata, Satoru and Nakatsukasa, Yuji and
             Takeda, Akiko},
  year    = {2017},
  pages   = {269--291}
}

@article{adachieigenvaluebasedalgorithmanalysis2019,
  title   = {Eigenvalue-based algorithm and analysis for nonconvex
             {QCQP} with one constraint},
  volume  = {173},
  number  = {1-2},
  journal = {Mathematical Programming},
  author  = {Adachi, Satoru and Nakatsukasa, Yuji},
  year    = {2019},
  note    = {Publisher: Springer},
  pages   = {79--116}
}

@inproceedings{agarwal_finding_2017,
  title     = {Finding approximate local minima faster than gradient
               descent},
  booktitle = {Proceedings of the 49th {Annual} {ACM} {SIGACT}
               {Symposium} on {Theory} of {Computing}},
  author    = {Agarwal, Naman and Allen-Zhu, Zeyuan and Bullins, Brian
               and Hazan, Elad and Ma, Tengyu},
  year      = {2017},
  pages     = {1195--1199}
}

@article{cartis_adaptive_2011,
  title      = {Adaptive cubic regularisation methods for unconstrained
                optimization. {Part} {I}: motivation, convergence and
                numerical results},
  volume     = {127},
  shorttitle = {Adaptive cubic regularisation methods for unconstrained
                optimization. {Part} {I}},
  number     = {2},
  journal    = {Mathematical Programming},
  author     = {Cartis, Coralia and Gould, Nicholas I. M. and Toint,
                Philippe L.},
  year       = {2011},
  pages      = {245--295}
}

@article{cartis_adaptive_2011-1,
  title      = {Adaptive cubic regularisation methods for unconstrained
                optimization. {Part} {II}: worst-case function- and
                derivative-evaluation complexity},
  volume     = {130},
  shorttitle = {Adaptive cubic regularisation methods for unconstrained
                optimization. {Part} {II}},
  number     = {2},
  journal    = {Mathematical Programming},
  author     = {Cartis, Coralia and Gould, Nicholas I. M. and Toint,
                Philippe L.},
  year       = {2011},
  pages      = {295--319}
}

@book{cartis_evaluation_2022,
  address    = {Philadelphia, PA},
  title      = {Evaluation {Complexity} of {Algorithms} for {Nonconvex}
                {Optimization}: {Theory}, {Computation} and
                {Perspectives}},
  shorttitle = {Evaluation {Complexity} of {Algorithms} for {Nonconvex}
                {Optimization}},
  author     = {Cartis, Coralia and Gould, Nicholas I. M. and Toint,
                Philippe L.},
  year       = {2022}
}

@article{curtis_inexact_2019,
  title   = {An inexact regularized {Newton} framework with a
             worst-case iteration complexity of for nonconvex
             optimization},
  volume  = {39},
  number  = {3},
  journal = {IMA Journal of Numerical Analysis},
  author  = {Curtis, Frank E. and Robinson, Daniel P. and Samadi,
             Mohammadreza},
  year    = {2019},
  note    = {Publisher: Oxford University Press},
  pages   = {1296--1327}
}

@article{curtis_trust_2017,
  title    = {A trust region algorithm with a worst-case iteration
              complexity of {$\mathcal O(\epsilon^{-3/2})$} for nonconvex
              optimization},
  volume   = {162},
  number   = {1},
  journal  = {Mathematical Programming},
  author   = {Curtis, Frank E. and Robinson, Daniel P. and Samadi,
              Mohammadreza},
  year     = {2017},
  pages    = {1--32}
}

@article{curtis_trust-region_2021,
  title    = {Trust-{Region} {Newton}-{CG} with {Strong} {Second}-{Order} {Complexity} {Guarantees} for {Nonconvex} {Optimization}},
  volume   = {31},
  issn     = {1052-6234},
  number   = {1},
  journal  = {SIAM Journal on Optimization},
  author   = {Curtis, Frank E. and Robinson, Daniel P. and Royer, Clément W. and Wright, Stephen J.},
  month    = jan,
  year     = {2021},
  note     = {Publisher: Society for Industrial and Applied Mathematics},
  pages    = {518--544}
}

@misc{curtis_worst-case_2022-1,
  title    = {Worst-{Case} {Complexity} of {TRACE} with {Inexact}
              {Subproblem} {Solutions} for {Nonconvex} {Smooth}
              {Optimization}},
  author   = {Curtis, Frank E. and Wang, Qi},
  year     = {2022},
  note     = {arXiv:2204.11322 [math]}
}

@article{den1995sufficient,
  title   = {A sufficient condition for self-concordance, with
             application to some classes of structured convex
             programming problems},
  author  = {den Hertog, Dick and Jarre, Florian and Roos, C and
             Terlaky, Tam{\'a}s},
  journal = {Mathematical Programming},
  volume  = {69},
  pages   = {75--88},
  year    = {1995}
}

@article{doikov_super-universal_2022,
  title   = {Super-{Universal} {Regularized} {Newton} {Method}},
  volume  = {34},
  issn    = {1052-6234},
  number  = {1},
  journal = {SIAM Journal on Optimization},
  author  = {Doikov, Nikita and Mishchenko, Konstantin and Nesterov, Yurii},
  month   = mar,
  year    = {2024},
  note    = {Publisher: Society for Industrial and Applied Mathematics},
  pages   = {27--56}
}

@article{doikov_gradient_2021,
  title    = {Gradient regularization of {Newton} method with {Bregman} distances},
  volume   = {204},
  issn     = {1436-4646},
  language = {en},
  number   = {1},
  journal  = {Mathematical Programming},
  author   = {Doikov, Nikita and Nesterov, Yurii},
  month    = mar,
  year     = {2024},
  pages    = {1--25}
}

@article{dussaultScalableAdaptiveCubic2023,
  title    = {Scalable adaptive cubic regularization methods},
  volume   = {207},
  issn     = {1436-4646},
  language = {en},
  number   = {1},
  urldate  = {2024-10-05},
  journal  = {Mathematical Programming},
  author   = {Dussault, Jean-Pierre and Migot, Tangi and Orban, Dominique},
  month    = sep,
  year     = {2024},
  pages    = {191--225}
}

@inproceedings{foster2019complexity,
  title        = {The complexity of making the gradient small in stochastic
                  convex optimization},
  author       = {Foster, Dylan J and Sekhari, Ayush and Shamir, Ohad and
                  Srebro, Nathan and Sridharan, Karthik and Woodworth,
                  Blake},
  booktitle    = {Conference on Learning Theory},
  pages        = {1319--1345},
  year         = {2019},
  organization = {PMLR}
}

@article{gallier2020schur,
  title   = {The Schur complement and symmetric positive semidefinite
             (and definite) matrices (2019)},
  author  = {Gallier, Jean},
  journal = {\url{https://www.cis.upenn.edu/jean/schur-comp.pdf}},
  year    = {2020}
}

@book{golub_matrix_2013,
  address   = {Baltimore},
  edition   = {4},
  title     = {Matrix Computations},
  isbn      = {978-1-4214-0794-4},
  publisher = {The Johns Hopkins University Press},
  author    = {Golub, Gene H. and Van Loan, Charles F.},
  year      = {2013}
}

@article{hanzely2022damped,
  title   = {A Damped Newton Method Achieves Global
             $\mathcal{O}\left(\frac{1}{k^2}\right)$ and Local Quadratic
             Convergence Rate},
  author  = {Hanzely, Slavom{\'\i}r and Kamzolov, Dmitry and
             Pasechnyuk, Dmitry and Gasnikov, Alexander and Richtarik,
             Peter and Takac, Martin},
  journal = {Advances in Neural Information Processing Systems},
  volume  = {35},
  pages   = {25320--25334},
  year    = {2022}
}

@article{he_quaternion_2022,
  title   = {Quaternion matrix decomposition and its theoretical
             implications},
  journal = {Journal of Global Optimization},
  author  = {He, Chang and Jiang, Bo and Zhu, Xihua},
  year    = {2022},
  note    = {Publisher: Springer},
  pages   = {1--18}
}

@article{heTechnicalReportHomogeneous,
  title    = {Technical {Report}: {The} {Homogeneous} {Second}-{Order} {Descent} {Framework} with {Inexact} {Eigenvalue} {Computations}},
  note     = {\url{https://bzhangcw.github.io/assets/pdfs/bisection.pdf}},
  year     = {2024},
  language = {en},
  author   = {He, Chang and Jiang, Yuntian and Zhang, Chuwen and Ge, Dongdong and Jiang, Bo and Ye, Yinyu}
}

@incollection{hilbert_beitrag_1970,
  address   = {Berlin, Heidelberg},
  title     = {Ein {Beitrag} zur {Theorie} des {Legendreschen}
               {Polynoms}},
  booktitle = {Algebra · {Invariantentheorie} {Geometrie}},
  author    = {Hilbert, David},
  editor    = {Hilbert, David},
  year      = {1970},
  pages     = {367--370}
}

@article{jia2022solving,
  title   = {Solving the cubic regularization model by a nested
             restarting Lanczos method},
  author  = {Jia, Xiaojing and Liang, Xin and Shen, Chungen and Zhang,
             Lei-Hong},
  journal = {SIAM Journal on Matrix Analysis and Applications},
  volume  = {43},
  number  = {2},
  pages   = {812--839},
  year    = {2022}
}

@article{kortanek1993polynomial,
  title   = {A polynomial barrier algorithm for linearly constrained
             convex programming problems},
  author  = {Kortanek, Kenneth O and Zhu, Jishan},
  journal = {Mathematics of Operations Research},
  volume  = {18},
  number  = {1},
  pages   = {116--127},
  year    = {1993}
}

@article{kuczynski_estimating_1992,
  title   = {Estimating the {Largest} {Eigenvalue} by the {Power} and
             {Lanczos} {Algorithms} with a {Random} {Start}},
  volume  = {13},
  number  = {4},
  journal = {SIAM Journal on Matrix Analysis and Applications},
  author  = {Kuczyński, J. and Woźniakowski, H.},
  year    = {1992},
  pages   = {1094--1122}
}

@article{lieder2020solving,
  title   = {Solving large-scale cubic regularization by a generalized
             eigenvalue problem},
  author  = {Lieder, Felix},
  journal = {SIAM Journal on Optimization},
  volume  = {30},
  number  = {4},
  pages   = {3345--3358},
  year    = {2020}
}

@book{luenberger_linear_2021,
  address = {Cham},
  series  = {International {Series} in {Operations} {Research} \&
             {Management} {Science}},
  title   = {Linear and {Nonlinear} {Programming}},
  volume  = {228},
  author  = {Luenberger, David G. and Ye, Yinyu},
  year    = {2021}
}

@article{masiha_stochastic_2022,
  title   = {Stochastic {Second}-{Order} {Methods} {Provably} {Beat}
             {SGD} {For} {Gradient}-{Dominated} {Functions}},
  journal = {arXiv preprint arXiv:2205.12856},
  author  = {Masiha, Saeed and Salehkaleybar, Saber and He, Niao and
             Kiyavash, Negar and Thiran, Patrick},
  year    = {2022}
}

@article{mishchenko_regularized_2023,
  title    = {Regularized {Newton} {Method} with {Global} {$\mathcal O(1/k^2)$} {Convergence}},
  volume   = {33},
  issn     = {1052-6234, 1095-7189},
  language = {en},
  number   = {3},
  journal  = {SIAM Journal on Optimization},
  author   = {Mishchenko, Konstantin},
  month    = sep,
  year     = {2023},
  pages    = {1440--1462}
}

@article{more_computing_1983,
  title   = {Computing a {Trust} {Region} {Step}},
  volume  = {4},
  number  = {3},
  journal = {SIAM Journal on Scientific and Statistical Computing},
  author  = {Moré, Jorge J. and Sorensen, D. C.},
  year    = {1983},
  pages   = {553--572}
}

@article{nesterov_cubic_2006,
  title    = {Cubic regularization of {Newton} method and its global
              performance},
  volume   = {108},
  number   = {1},
  journal  = {Mathematical Programming},
  author   = {Nesterov, Yurii and Polyak, B.T.},
  year     = {2006},
  pages    = {177--205}
}

@book{nesterov_interior-point_1994,
  title  = {Interior-{Point} {Polynomial} {Algorithms} in {Convex}
            {Programming}},
  author = {Nesterov, Yurii and Nemirovskii, Arkadii},
  year   = {1994}
}

@book{nesterov_lectures_2018,
  address   = {Cham},
  title     = {Lectures on {Convex} {Optimization}},
  publisher = {Springer International Publishing},
  author    = {Nesterov, Yurii},
  year      = {2018}
}

@article{nesterov2008accelerating,
  title   = {Accelerating the cubic regularization of Newton’s method
             on convex problems},
  author  = {Nesterov, Yu},
  journal = {Mathematical Programming},
  volume  = {112},
  number  = {1},
  pages   = {159--181},
  year    = {2008}
}

@article{nesterov2012make,
  title   = {How to make the gradients small},
  author  = {Nesterov, Yurii},
  journal = {Optima. Mathematical Optimization Society Newsletter},
  number  = {88},
  pages   = {10--11},
  year    = {2012}
}

@book{nocedal_numerical_2006,
  title    = {Numerical optimization},
  author   = {Nocedal, Jorge and Wright, Stephen},
  year     = {2006}
}

@article{rojas_new_2001,
  title    = {A {New} {Matrix}-{Free} {Algorithm} for the
              {Large}-{Scale} {Trust}-{Region} {Subproblem}},
  volume   = {11},
  number   = {3},
  journal  = {SIAM Journal on Optimization},
  author   = {Rojas, Marielba and Santos, Sandra A. and Sorensen, Danny
              C.},
  year     = {2001},
  note     = {Publisher: Society for Industrial and Applied
              Mathematics},
  pages    = {611--646}
}

@article{royer_complexity_2018,
  title   = {Complexity analysis of second-order line-search algorithms
             for smooth nonconvex optimization},
  volume  = {28},
  number  = {2},
  journal = {SIAM Journal on Optimization},
  author  = {Royer, Clément W. and Wright, Stephen J.},
  year    = {2018},
  note    = {Publisher: SIAM},
  pages   = {1448--1477}
}

@article{royer_newton-cg_2020,
  title    = {A {Newton}-{CG} algorithm with complexity guarantees for smooth unconstrained optimization},
  volume   = {180},
  issn     = {1436-4646},
  language = {en},
  number   = {1},
  journal  = {Mathematical Programming},
  author   = {Royer, Clément W. and O’Neill, Michael and Wright, Stephen J.},
  month    = mar,
  year     = {2020},
  pages    = {451--488}
}

@book{saad_numerical_2011,
  address   = {Philadelphia},
  edition   = {Rev. ed},
  series    = {Classics in applied mathematics},
  title     = {Numerical methods for large eigenvalue problems},
  isbn      = {978-1-61197-072-2},
  number    = {66},
  publisher = {Society for Industrial and Applied Mathematics},
  author    = {Saad, Y.},
  year      = {2011}
}

@article{sorensen_newtons_1982,
  title   = {Newton’s method with a model trust region modification},
  volume  = {19},
  number  = {2},
  journal = {SIAM Journal on Numerical Analysis},
  author  = {Sorensen, Danny C.},
  year    = {1982},
  note    = {Publisher: SIAM},
  pages   = {409--426}
}

@article{sturmconesnonnegativequadratic2003,
  title   = {On cones of nonnegative quadratic functions},
  volume  = {28},
  number  = {2},
  journal = {Mathematics of Operations research},
  author  = {Sturm, Jos F. and Zhang, Shuzhong},
  year    = {2003},
  note    = {Publisher: INFORMS},
  pages   = {246--267}
}

@book{ye_interior_1997,
  address    = {New York Weinheim},
  title      = {Interior {Point} {Algorithms}: {Theory} and {Analysis}},
  isbn       = {978-0-471-17420-2},
  shorttitle = {Interior {Point} {Algorithms}},
  publisher  = {Wiley},
  author     = {Ye, Yinyu},
  year       = {1997}
}

@article{ye2017second,
  title  = {A Second-Order Path-Following Algorithm for Unconstrained
            Convex Optimization},
  author = {Ye, Yinyu},
  year   = {2017}
}

@article{yeApproximatingQuadraticProgramming1999a,
  title    = {Approximating quadratic programming with bound and quadratic constraints},
  volume   = {84},
  issn     = {1436-4646},
  language = {en},
  number   = {2},
  journal  = {Mathematical Programming},
  author   = {Ye, Yinyu},
  month    = feb,
  year     = {1999},
  pages    = {219--226}
}

@article{yeCombiningBinarySearch1994,
  title    = {Combining {Binary} {Search} and {Newton}'s {Method} to {Compute} {Real} {Roots} for a {Class} of {Real} {Functions}},
  volume   = {10},
  issn     = {0885-064X},
  number   = {3},
  journal  = {Journal of Complexity},
  author   = {Ye, Yinyu},
  year     = {1994},
  pages    = {271--280}
}

@misc{zhang_homogenous_2022,
  title    = {A {Homogenous} {Second}-{Order} {Descent} {Method} for
              {Nonconvex} {Optimization}},
  author   = {Zhang, Chuwen and Ge, Dongdong and He, Chang and Jiang, Bo
              and Jiang, Yuntian and Xue, Chenyu and Ye, Yinyu},
  year     = {2022},
  note     = {arXiv:2211.08212 [math]}
}

@article{zhu1992path,
  title   = {A path following algorithm for a class of convex
             programming problems},
  author  = {Zhu, Jishan},
  journal = {Zeitschrift f{\"u}r Operations Research},
  volume  = {36},
  pages   = {359--377},
  year    = {1992}
}

\appendix

\section{Proofs in \autoref{sec.warmup}}

\subsection{Technical lemmas}
We introduce several technical lemmas.

\begin{lemma}\label{lem.psdpseudo}
    If $A\in \real^{n\times n}, A \succeq 0$, then its pseudo-inverse satisfies $A^\star \succeq 0$.
\end{lemma}
\begin{proof}
    Take the spectral decomposition of $A = UDU^T$, and assume that there are $d$ strictly positive eigenvalues. We see $A^\star = U_+D^{-1}_+U_+^T $ where $D_+ \in \real^{n\times d}, U_+ \in \real^{n\times d}$ are strictly positive eigenvalues and their corresponding eigenvectors. So we conclude that $A^\star \succeq 0$.
\end{proof}
\begin{lemma}[Strict negative curvature]\label{lem.strict.curv}
    Consider $A\in\real^{n\times n}, b\in\real^n$, and $c\in \real$. Assume that $b\neq 0$ and an aggregated matrix,
    $$
        B = \begin{bmatrix}
            A   & b \\
            b^T & c
        \end{bmatrix}
    $$
    If either one of the following holds, then $\lambda_1(B) < 0$:
    \begin{enumerate}[$(a)$]
        \item $c \le 0$ or $\lambda_1(A) < 0$.
        \item $\lambda_1(A)=0$ and $b \not \perp \mathcal{S}_1$, where $\mathcal{S}_1$ is the subspace spanned by the eigenvectors associated with $\lambda_1(A)$.
    \end{enumerate}
\end{lemma}
\begin{proof}
    For the case $c < 0$ or $\lambda_1(A) < 0$, it is obvious.

    Now we consider the case $c = 0$, for any vector $\xi = [\eta\cdot v; t] \in \real^{n+1}$ which satisfies $v \in \real^n, t\neq 0$. We obtain that
    \begin{align*}
        R_B([\eta\cdot v; t]):=\begin{bmatrix}
                                   \eta  v \\ t
                               \end{bmatrix}^T B \begin{bmatrix}
                                                     \eta  v \\ t
                                                 \end{bmatrix} = \eta^2 \cdot\left( v^TAv\right) + 2 \eta \cdot \left(b^Tv\cdot t\right).
    \end{align*}
    It is easy to choose $v, t$ such that $b^Tv\cdot t \neq 0$, and thus there exists $\eta$ such that $R_B([\eta\cdot v; t]) < 0$, implying $\xi$ corresponds to a negative curvature.

    The only left case is $\lambda_1(A) = 0$ and $b \not \perp \mathcal{S}_1$. For any nonzero vector $u \in \mathcal{S}_1$, $u$ satisfies $b^Tu \neq 0$. Using a similar argument, for the vector $[u; \eta t]$ with $t \neq 0$, we have,
    \begin{align*}
        R_B([u; \eta t]) & :=\begin{bmatrix}
                                 u \\ \eta t
                             \end{bmatrix}^T B \begin{bmatrix}
                                                   u \\\eta t
                                               \end{bmatrix} = \left( u^TAu\right) + 2 \eta \cdot \left(b^Tu\cdot t\right) + \eta^2 \cdot c t^2 \\
                         & =2 \eta \cdot \left(b^Tu\cdot t\right) + \eta^2 \cdot c t^2.
    \end{align*}
    Similarly, since $b^Tu \cdot t \neq 0$, so there exists $\eta$ such that $R_B([u; \eta t]) < 0$. This completes the proof.
\end{proof}

\subsection{Basic results}
\begin{theorem}\label{lem.eigenvector computation}
    For the aggregated matrix $F_k$, the Lanczos method takes
    $$
        \min\left\{O\left(\sqrt{\frac{\|F_k\|}{\varepsilon}}\log\left(\frac{n+1}{\delta}\right)\right), O\left(\sqrt{\frac{\|F_k\|}{\lambda_2(F_k) - \lambda_1(F_k)}}\log\left(\frac{n+1}{\varepsilon\cdot \delta}\right)\right)\right\}
    $$
    iterations to compute an estimate $\xi$ such that $|\xi - \lambda_1(F_k)| \le \varepsilon$ in probability $1-\delta$.
\end{theorem}
\begin{proof}
    We denote $\lambda_{r+1}$ be the maximum eigenvalue (since $F_k$ increases the dimension by 1).
    Computing the minimum eigenvalue of $F_k$ is equivalent to compute the maximum of $B_k \cdot I-F_k$ for some $B_k \ge \lambda_{r+1}(F_k)$, then by \cite[Theorem 4.2 $(a)$]{kuczynski_estimating_1992}, the $k$-th Lanczos iteration outputs the estimate $-\xi$ to $\lambda_{r+1}(B_k \cdot I - F_k)$ such that
    \begin{align*}
        \mathbb{P}\left\{\frac{\left|\lambda_{r+1}(B_k \cdot I -F_k) - (-\xi)\right|}{ (\lambda_{r+1}(B_k \cdot I - F_k) - \lambda_1(B_k \cdot I - F_k))} \ge \tilde \varepsilon\right\} \le 1.648 \sqrt{n+1} \exp\{-\sqrt{\tilde \varepsilon}\cdot(2k-1)\}.
    \end{align*}
    By ceiling the RHS to $\delta$, we have
    then the first argument follows by setting $\tilde \varepsilon = \varepsilon/(\lambda_{r+1}(F_k) - \lambda_1(F_k))$ and observe
    \begin{align*}
         & \lambda_{r+1}(B_k \cdot I - F_k) - \lambda_1(B_k \cdot I - F_k) = \lambda_{r+1}(F_k) - \lambda_1(F_k) \\
         & \lambda_{r+1}(B_k \cdot I -F_k) - (-\xi) =  \xi - \lambda_1(F_k).
    \end{align*}
    as well as $\lambda_{r+1}(F_k) - \lambda_1(F_k) \le 2\|F_k\|$.
    By \cite[Theorem 4.2 $(b)$]{kuczynski_estimating_1992}, it also holds that
    \begin{align*}
         & \mathbb{P}\left\{\frac{\left|\lambda_{r+1}(B_k \cdot I -F_k)
        - (-\xi)\right|}{ (\lambda_{r+1}(B_k \cdot I - F_k) - \lambda_1(B_k \cdot I - F_k))} \ge \tilde \varepsilon\right\}                                                                     \\
         & \qquad\le 1.648 \sqrt{n+1} \left(\frac{1-\sqrt{(\lambda_{r+1} - \lambda_r)/(\lambda_{r+1}-\lambda_1)}}{1+\sqrt{(\lambda_{r+1} - \lambda_r)/(\lambda_{r+1}-\lambda_1)}}\right)^{k-1}.
    \end{align*}
    Similarly, we have,
    \begin{align*}
        k \ge 1 + \log(1.648\sqrt{n+1} / (\delta \sqrt{\tilde \varepsilon})) / \log\left(\frac{1+\sqrt{(\lambda_{r+1} - \lambda_r)/(\lambda_{r+1}-\lambda_1)}}{1-\sqrt{(\lambda_{r+1} - \lambda_r)/(\lambda_{r+1}-\lambda_1)}}\right).
    \end{align*}
    Then observe $\log (1+1 / t) \geq \frac{1}{1 / 2+t}$, and $\lambda_{r+1}(B_k \cdot I -F_k) - \lambda_r(B_k \cdot I -F_k) = \lambda_2(F_k) - \lambda_1(F_k)$ we conclude the second argument.
\end{proof}
The first estimate in the bracket is usually referred to as the ``gap-free" complexity and the latter is called ``gap-dependent".
In nonconvex optimization, where the spectrum of the Hessian distributes along the real line, we can only rely on the ``gap-free" guarantees.  In the original work \cite{zhang_homogenous_2022}, they rigorously established that finding a required ``inexact'' solution of GHM matches the complexity of obtaining inexact Newton-type directions \cite{agarwal_finding_2017,royer_newton-cg_2020,curtis_trust-region_2021}.

\begin{lemma}[Estimation of $F_k$]\label{lem.eigengap.free}
    Consider the the aggregated matrix $F_k$, it holds that
    \begin{align*}
        \lambda_2(F_k) - \lambda_1(F_k) & \ge \lambda_1(H_k) + \max_{j}\left\{\frac{- \lambda_j(H_k) - \delta_k + \sqrt{(\lambda_j(H_k) - \delta_k)^2 + 4\|P_{\mathcal S_j}(\gk)\|^2}}{2}\right\} \\
        \lambda_{r+1}(F_k)              & \le \lambda_r(\Hk) + \delta_k - \lambda_1(F_k)
    \end{align*}
\end{lemma}
\begin{proof}
    By the second-order condition \eqref{eq.homoeig.soc} in \autoref{lemma.optimal condition of subproblem}, we have,
    \begin{align*}
        \Hk + \theta_k I - \frac{\gk\gk^T}{\theta_k + \delta_k} \succeq 0,
    \end{align*}
    by the Schur complement.
    By projecting onto $\mathcal S_j(\Hk), j = 1,...d$, it holds that
    \begin{align*}
        \theta_k + \lambda_j(\Hk) - \frac{\|P_{\mathcal S_j}(\gk)\|^2}{\theta_k + \delta_k} \ge 0, \forall j
    \end{align*}
    Reorganize the terms, we have
    \begin{align*}
        \theta_k^2 + (\lambda_j(\Hk) + \delta_k) \theta_k + \lambda_j(\Hk)\delta_k - \|P_{\mathcal S_j}(\gk)\|^2 \ge 0,
    \end{align*}
    which implies
    \begin{align*}
        \theta_k \ge \max_{j}\left\{\frac{-(\lambda_j(\Hk) + \delta_k) + \sqrt{(\lambda_j(\Hk) - \delta_k)^2 + 4\|P_{\mathcal S_j}(\gk)\|^2}}{2}\right\}.
    \end{align*}
    By Cauchy Interlace Theorem, we have
    \begin{align*}
        \lambda_2(F_k) - \lambda_1(F_k) & = \lambda_2(F_k) + \theta_k \ge \lambda_1(H_k) + \theta_k                                                                                                \\
                                        & \ge  \max_{j}\left\{\frac{2\lambda_1(H_k)-(\lambda_j(\Hk) + \delta_k) + \sqrt{(\lambda_j(\Hk) - \delta_k)^2 + 4\|P_{\mathcal S_j}(\gk)\|^2}}{2}\right\}. \\
    \end{align*}
    Notice
    \begin{align*}
        \tr(F_k) & = \sum_{j=1}^{r+1}\lambda_j(F_k) = \sum_{j=1}^{r}\lambda_j(H_k) + \delta_k.
    \end{align*}
    so we have,
    \begin{align*}
        \lambda_{\max}(F_k) := \lambda_{r+1}(F_k) & = \sum_{j=2}^r \left (\lambda_j(H_k)-\lambda_j(F_k)\right ) + \lambda_1(\Hk) - \lambda_1(F_k) + \delta_k       \\
                                                  & \le \sum_{j=2}^r\left ( \lambda_j(H_k)-\lambda_{j-1}(H_k)\right ) + \delta_k + \lambda_1(\Hk) - \lambda_1(F_k) \\
                                                  & = \lambda_r(\Hk) - \lambda_1(\Hk) + \delta_k + \lambda_1(\Hk) - \lambda_1(F_k)                                 \\
                                                  & = \lambda_r(\Hk) + \delta_k - \lambda_1(F_k)                                                                   \\
                                                  & = \lambda_r(\Hk) + \delta_k + \theta_k.
    \end{align*}
\end{proof}

\subsection{Proof of \autoref{thm.eigengap.better}}\label{proof.thm.eigengap.better}
\begin{proof}
    For part $(a)$, we refer to \autoref{lem.eigengap.free} and set {$\delta_k = -\epsilon_{\mathrm{L}}$}.
    We first note that $\forall \eta \le 1/2\sqrt n$, there exists $j'$ such that $\|P_{\mathcal S_{j'}}(\gk)\|^2 \ge \eta^2 \|\gk\|^2$; otherwise
    \begin{align*}
        \sum_{j=1}^r\|P_{\mathcal S_j}(\gk)\|^2 < \sum_{j=1}^r\eta^2 \|\gk\|^2 \le \frac{\|\gk\|^2}{4},
    \end{align*}
    which contradicts the fact that $\bigcup_j\mathcal S_j = \real^n$. In this view, $\forall \eta \in [0, 1/2\sqrt n]$,
    \begin{subequations}
        \begin{align}
            \lambda_2(F_k) - \lambda_1(F_k) & \ge \lambda_1(H_k) + \frac{- \lambda_{j'}(H_k) + \epsilon_{\mathrm{L}} + \sqrt{(\lambda_{j'}(H_k) + \epsilon_{\mathrm{L}})^2 + 4\eta^2\|\gk\|^2}}{2} \\
            \label{eq.lowerb.agg}           & \ge \lambda_1(H_k) + \frac{- U_H + \epsilon_{\mathrm{L}} + \sqrt{(U_H + \epsilon_{\mathrm{L}})^2 + 4\eta^2\|\gk\|^2}}{2} > 0.
        \end{align}
    \end{subequations}
    Note \eqref{eq.lowerb.agg} holds since the one-dimensional function
    $$
        \ell(x) := \frac{- x + \epsilon_{\mathrm{L}} + \sqrt{(x + \epsilon_{\mathrm{L}})^2 + 4\eta^2\|\gk\|^2}}{2}
    $$
    is decreasing in $x > 0$ and $\lambda_{j'}(H_k) \le U_H$. By definition, we conclude
    \begin{align*}
        \kappa_{\mathrm{L}}(F_k) \le \frac{2(\lambda_r(\Hk) - \epsilon_{\mathrm{L}} - \lambda_1(F_k))}{- U_H + \epsilon_{\mathrm{L}} + \sqrt{(U_H + \epsilon_{\mathrm{L}})^2 + \|\gk\|^2/n}} < \infty
    \end{align*}
    as $\lambda_1(F_k)$ is bounded.

        {
            For part $(b)$, we observe that
            \begin{align*}
                \frac{\kappa_{\mathrm{L}}(F_k)}{\kappa(\Hk+\epsilon_{\mathrm{N}} I)}
                 & \le \frac{2\left\|F_k\right\|}{\lambda_2(F_k)-\lambda_1(F_k)} \cdot \frac{\lambda_1(H_k)+\epsilon_{\mathrm{N}}}{\lambda_r(H_k)+\epsilon_{\mathrm{N}}}                        \\
                 & =\frac{2\left\|F_k\right\|}{\epsilon_{\mathrm{N}} + U_H} \cdot (\lambda_1(H_k)+\epsilon_{\mathrm{N}}) \cdot \underbrace{\frac{1}{\lambda_2(F_k)-\lambda_1(F_k)}}_{\ddagger}.
            \end{align*}
            Note that for $\ddagger$, we have
            \begin{align*}
                \lambda_2(F_k) - \lambda_1(F_k) & \ge \epsilon_{\mathrm{L}} + \frac{- (U_H + \epsilon_{\mathrm{L}}) + \sqrt{(U_H + \epsilon_{\mathrm{L}})^2 + 4\eta^2\|\gk\|^2}}{2}                                  \\
                                                & = \epsilon_{\mathrm{L}} + \frac{U_H + \epsilon_{\mathrm{L}}}{2} \cdot \left(\sqrt{1 + \left(\frac{2\eta\|g_k\|}{U_H + \epsilon_{\mathrm{L}}}\right)^2} - 1\right).
            \end{align*}
            Let $\Gamma = \frac{2\eta\|g_k\|}{U_H + \epsilon_{\mathrm{L}}}$, we have
            \begin{align*}
                \lambda_2(F_k) - \lambda_1(F_k) \ge \epsilon_{\mathrm{L}} + \frac{U_H + \epsilon_{\mathrm{L}}}{2} \cdot \left(\frac{\Gamma^2}{2} - \frac{\Gamma^4}{8} + O(\Gamma^6)\right).
            \end{align*}
            Since $\Gamma = \Theta \left(\|g_k\|\right)$ and $U_H + \epsilon_{\mathrm{L}}$ can be viewed as a constant, it implies that
            \begin{align*}
                \frac{\kappa_{\mathrm{L}}(F_k)}{\kappa(\Hk+\epsilon_{\mathrm{N}} I)}
                \le O\left(\frac{\epsilon_{\mathrm{N}} }{\frac{\|g_k\|^2}{U_H+\epsilon_{\mathrm{L}}} + \epsilon_{\mathrm{L}}} \right).
            \end{align*}
            This completes the proof.
        }
\end{proof}

\section{Proofs in \autoref{sec.basic.pd}}\label{app.appendix a}

The lemma below is a slightly re-organized version of \citet[Lemma 3.3]{rojas_new_2001}.
\begin{lemma} \label{lemma.rojas2}
    Define $p_j, \tilde \alpha_j, 1\le j \le r$ similar to \autoref{lemma.rojas1} and assume $\phi_k \perp \mathcal S_i(\Hk)$,~ for ~$i=1,2, \ldots, \ell$ ~and~ $1 \leq \ell<r$, then,
    \begin{enumerate}[(a)]
        \item if $\delta_k=\tilde{\alpha}_j, 1 \leq j \leq \ell$, then $\lambda_j(F_k)=\lambda_j(\Hk), j=1,2, \ldots, \ell$;
        \item if $\delta_k<\tilde{\alpha}_1$, then $\lambda_1(F_k) < \lambda_1(\Hk)$ and $\lambda_j(F_k)=\lambda_{j-1}(\Hk), j=2, \ldots, \ell+1$;
        \item if $\tilde{\alpha}_{i-1}<\delta_k<\tilde{\alpha}_i, 2 \leq i \leq \ell$,
              then ~$\lambda_j(F_k)=\lambda_j(\Hk)$~ for ~$j=1, \ldots, i-1$,
              and ~$\lambda_j(F_k)=\lambda_{j-1}(\Hk)$~ for ~$j=i+1, \ldots, \ell+1$;
        \item if $\delta_k>\tilde{\alpha}_{\ell}$, then $\lambda_j(F_k)=\lambda_j(\Hk), j=1,2, \ldots, \ell$.
    \end{enumerate}
\end{lemma}
\subsection{Proof of \autoref{lemma.uppertheta}}\label{sec.upper.theta.lemma}
\begin{proof}
    Recall the definition of $\theta_k$ in \eqref{eq.ghqm}, we see that
    \begin{align}
        - \theta_k & = \min_{v^Tv + t^2 \le 1} v^T\Hk v + 2\phi_k^Tv \cdot t + \delta_k t^2                        \notag                            \\
                   & \ge \min_{v^Tv+ t^2 \le 1} v^T\Hk v + \delta_k t^2 +\min_{v^Tv+ t^2 \le 1} 2\phi_k^Tv \cdot t. \label{eq.bound theta two terms}
    \end{align}
    Now we turn to investigate the two terms respectively, for the first term, we have
    \begin{align*}
        \min_{v^Tv+ t^2 \le 1} ~ v^T\Hk v + \delta_k t^2 ~ & = \min_{\|[v; t]\| \le 1} \begin{bmatrix}
                                                                                           v \\ t
                                                                                       \end{bmatrix}^T
        \begin{bmatrix}
            \Hk & 0        \\
            0^T & \delta_k
        \end{bmatrix}
        \begin{bmatrix}
            v \\ t
        \end{bmatrix}                                                                                 \\
                                                           & =
        \begin{cases}
            0                              & \Hk \succeq 0, \delta_k \ge 0 \\
            \lambda_1\left(\begin{bmatrix}
                               \Hk & 0        \\
                               0^T & \delta_k
                           \end{bmatrix}\right) & o.w.
        \end{cases}
    \end{align*}
    Therefore, we conclude that
    \begin{equation}\label{eq.uppertheta.aux1}
        \min_{v^Tv+ t^2 \le 1} ~
        v^T\Hk v + \delta_k t^2 ~  = \min\{\lambda_1,\delta_k, 0\}.
    \end{equation}
    As for the second term, with a similar argument, it is easy to obtain that
    \begin{equation}\label{eq.uppertheta.aux2}
        \min_{v^Tv+ t^2 \le 1} ~
        ~2 \phi_k^Tv \cdot t ~   = \lambda_1 \left(\begin{bmatrix}
                0        & \phi_k \\
                \phi_k^T & 0
            \end{bmatrix}\right) = - \|\phi_k\|.
    \end{equation}
    Substitute the above equations \eqref{eq.uppertheta.aux1} and \eqref{eq.uppertheta.aux2} into \eqref{eq.bound theta two terms} and we complete the proof.
\end{proof}

\subsection{Proof of \autoref{lemma.ordering of tilde alpha 1}}
\begin{proof}
    Since $\Hk-\lambda_1(\Hk) I \succeq 0$, by \autoref{lem.psdpseudo}, $\left(\Hk-\lambda_1(\Hk) I\right)^\star \succeq 0$, then
    $$\tilde\alpha_1 - \lambda_1(\Hk) = - \phi_k^Tp_1 = \phi_k^T\left(\Hk-\lambda_1(\Hk) I\right)^{\star}\phi_k \ge 0.$$
    This completes the proof.
\end{proof}

\subsection{Proof of \autoref{lem.continuity.theteomega}}\label{sec.proof.continuity.theteomega}
\begin{proof}
    We consider two cases as follows.

    \textbf{Case (1) $\phi_k \not \perp \mathcal S_1(\Hk)$.}  In this case, $t_k \neq 0$ for all $\delta_k$. The optimal condition \eqref{eq.homoeig.foc} reduces to a Newton-like equation as $\dk = v_k/t_k$ is well-defined:
    $$
        \left(\Hk + \theta_k I\right)\dk = - \phi_k, ~-\phi_k^T\dk=\theta_k + \delta_k.
    $$
    From second-order condition \eqref{eq.homoeig.soc}, we obtain that
    $$
        \begin{bmatrix}
            \Hk + \theta_k I & \phi_k              \\
            \phi_k^T         & \theta_k + \delta_k
        \end{bmatrix}\succeq 0,
    $$
    and it implies that $\Hk + \theta_k I \succ 0$, otherwise a negative curvature always exists by \autoref{lem.strict.curv}. Therefore, the inverse is well-defined and
    \begin{equation}\label{eq.dk.invert}
        \dk = -\left(\Hk + \theta_k I\right)^{-1}\phi_k,
    \end{equation}

    It follows that $$\theta_k + \delta_k = -\dk^T\phi_k = \phi_k^T\left(\Hk + \theta_k I\right)^{-1}\phi_k.$$ Take the derivative of both sides:
    \begin{align*}
        \diff{\theta_k} (\theta_k + \delta_k)                                    & = \diff{\theta_k}\left(\phi^T_k\left(\Hk + \theta_k I\right)^{-1}\phi_k\right)                                                            & \\
        \Rightarrow  ~1 + \diff{\theta_k}\delta_k                                & = - \phi_k^T\left(\Hk + \theta_k I\right)^{-2}\phi_k                                         \stackrel{\eqref{eq.dk.invert}}{=} -\Delta_k   \\
        \Rightarrow ~                                    \diff{\delta_k}\theta_k & = - \frac{1}{\Delta_k + 1}.
    \end{align*}
    So we have \eqref{eq.diff.omega}, which implies that $\theta_k$ is continuous and decreasing in this case. As for the convexity of $\theta_k$, we first observe that $\|\dk\|$ is decreasing as $\theta_k$ increasing due to that $\|\dk\| = \|\left(\Hk + \theta_k I\right)^{-1}\phi_k\|$. Since $\delta_k \nearrow, \theta_k \searrow$, it follows that $\|\dk\| \nearrow$ and further implies $\diff{\delta_k}\theta_k \nearrow$, and thus we conclude that $\theta_k$ is convex over $t_k \neq 0$. By \autoref{lemma.rojas1}, $t_k \neq 0$ for all $\delta_k \in \real$ as $\phi_k \not \perp \mathcal S_1(\Hk)$. Combining all these facts, we conclude that $\theta_k$ is decreasing, continuous, convex, and differentiable. Since $\theta_k \ge 0$, the same thing holds for $\omega_k$.

    \textbf{Case (2) $\phi_k \perp \mathcal S_1(\Hk)$.} In this case, we split into the following three sub-cases:
    \begin{enumerate}[(i)]
        \item $\lambda_1(H) > 0$. Under this subcase, by Schur complement, we conclude $\lambda_1(F_k) > 0$ if $\delta_k > \phi_k^T\Hk^{-1} \phi_k$, and hence $\theta_k = 0$. Also, note that
              $$\phi_k^T\Hk^{-1}\phi_k < \tilde \alpha_1 := \phi_k^T\left(\Hk - \lambda_1(H_k) I\right)^*\phi_k,$$
              meaning $\theta_k = 0$ long before reaching $\delta_k = \tilde \alpha_1$. This indicates that $t_k \neq 0$ when $\delta_k \leq \phi_k^T\Hk^{-1}\phi_k$, and it further implies that $\theta_k$ is convex and continuous by \textbf{Case (1)}. When $\delta_k > \phi_k^T\Hk^{-1}\phi_k$, $\theta_k = 0$ is a constant.
        \item $\lambda_1(H) = 0$. Based on \autoref{lem.strict.curv}, we know that the strict negative curvature of $F_k$ always exists, which means $\theta_k = -\lambda_1(F_k)$ for all $\delta_k$. First, by \autoref{corr.necc.t0}, for any $\delta_k \le \tilde \alpha_1$, there always exists an associated $t_k \neq 0$. Together the results in \textbf{Case (1)}, we conclude that $\theta_k$ is convex and continuous if $\delta_k \le \tilde \alpha_1$. Moreover, $\theta_k = -\lambda_1(F_k) = -\lambda_1(\Hk)$ for any $\delta_k \ge \tilde \alpha_1$. Therefore, $-\theta_k$ is convex and continuous for all $\delta_k$.

        \item $\lambda_1(H) < 0$. This subcase is similar to the above one, we omit the proof and present the results $\theta_k = -\lambda_1(F_k)$ for all $\delta_k$, $\theta_k$ is convex and continuous when $\delta_k \le \tilde \alpha_1$ and $\theta_k = -\lambda_1(F_k) = -\lambda_1(H_k)$ when $\delta_k \geq \tilde \alpha_1$.
    \end{enumerate}
\end{proof}

\subsection{Proof of \autoref{lem.differentiablity of h}}\label{sec.pf.lem.differentiablity of h}
\begin{proof}
    For the case where $\phi_k \perp \mathcal{S}_1(\Hk)$, if $\delta_k = \tilde\alpha_1$ then by definition \autoref{def.aux}, $\Delta_k$ is not continuous, so we conclude the discontinuity of $h_k$.

    If $\phi_k \not\perp \mathcal{S}_1(\Hk)$, by the optimality condition \eqref{eq.homoeig.soc} and \autoref{lem.strict.curv}, we have
    \begin{equation*}
        \Hk+\theta_k I \succ 0,
    \end{equation*}
    then we apply the spectral decomposition to $\Hk$ and denote
    \begin{equation}\label{eq.spec.decomp}
        \Hk = \sum_{i=1}^r \lambda_i(\Hk) v_i v_i^T,\quad \beta_i =  \phi_k^T v_i, \quad i = 1, \cdots, d.
    \end{equation}
    Therefore, we can rewrite $h_k(\delta)$ in the form:
    \begin{align}
        \nonumber
        h_k(\delta) & = \frac{\theta_k^2}{\|\dk\|^2} = \frac{\theta_k^2}{\phi_k^T(\Hk+\theta_k I)^{-2}\phi_k} \\
        \label{eq.new form h}
                    & = \frac{\theta_k^2}{\sum_{i=1}^n \frac{\beta_i^2}{(\lambda_i+\theta_k)^2}},
    \end{align}
    thus function $h_k(\delta)$ can be viewed as a composite function of $\omega_k$. Recall that $\omega_k(\delta)$ is differentiable by \autoref{lem.continuity.theteomega}, and hence thus $h_k(\delta)$ is differentiable. From the chain rule and the result of \autoref{lem.continuity.theteomega}, we conclude that
    \begin{align}
        \nonumber
        h_k'(\delta) & = \diff{\theta} h_k(\delta) \diff{\delta} \theta                                        \\
        \nonumber
                     & = - \frac{2\theta \|d\|^2+\theta^2 d^T(\Hk+\theta I)^{-1}d}{\|d\|^2}\frac{1}{1+\|d\|^2} \\
                     & = - \frac{2\theta \|d\|^2+\theta^2 d^T(\Hk+\theta I)^{-1}d}{\|d\|^2(1+\|d\|^2)},
        \label{eq.derivative of h}
    \end{align}
    which means that $h_k(\delta)$ is monotone decreasing.

\end{proof}

\section{Proofs in \autoref{sec.recover}}

\subsection{Proof of \autoref{thm.recover.trs}}\label{proof.thm.recover.trs}
\begin{proof}
    For trust-region type methods, the case where $g_k \perp \mathcal{S}_1(H_k)$ usually needs separate treatments. For conciseness, we simply consider the case $g_k \not\perp \mathcal{S}_1(H_k)$. Let $([v_k;t_k],\theta_k)$ be the primal-dual solution of the GHM subproblem, by \autoref{corr.necc.t0}, we obtain that $t_k \neq 0$ for all $\delta_k \in \mathbb{R}$, which leads to that $\Delta_k(\delta_k) = \|v_k/t_k\|^2$ is continuous. Therefore, we could apply a similar bisection procedure in \autoref{sec.bisection procedure} to find a  $\delta_k$ such that
    \[
        \Delta_k(\delta_k) = \Delta^2.
    \]
    Combining the optimal condition in \eqref{eq.homoeig.soc}, set $d_k = v_k/t_k$, it follows that
    $$
        \begin{aligned}
             & H_k + \theta_k \cdot I \succeq 0,   \\
             & (H_k + \theta_k \cdot I)d_k = -g_k, \\
             & \theta_k(\|d_k\| - \Delta) = 0,
        \end{aligned}
    $$
    which implies $(d_k, \theta_k)$ is the solution of the trust-region subproblem \eqref{eq.trsubp}.  By selecting an appropriate $\delta_k$, we can always recover the solution to the trust-region subproblem \eqref{eq.trsubp} based on $[v_k;t_k]$.  Since we already committed to an adaptive HSODM, we only briefly discuss the bracketing stage for $\delta_k$ here.

    Since $\gk \not \perp \mathcal{S}_1(H_k)$, by \autoref{lem.property.h} we know that $\Delta_k$ is continuous in $\delta_k$. Then the key of our analysis is to find an interval $(\delta_{low},\delta_{up})$ such that $$\Delta^2 \in \left [\Delta_k(\delta_{low}),\Delta_k(\delta_{up})\right ].$$
    We focus on the case $\Delta <1$ since the opposite case is similar.

    First, from \autoref{lem.relation t and delta}, we can set $\delta_{up} = \lambda_r(\Hk)$, then $\Delta^2\leq \Delta_k(\delta_{up})$ follows.

    Next, we give an estimate for $\delta_{low}$. In fact, without loss of generality, we can assume $\lambda_1(\Hk)\leq 0$ and set $\delta_{low} = \lambda_1(\Hk)-\|\phi_k\|\frac{1-\Delta^2}{\Delta_2}$. The case $\lambda_1(\Hk)>0$ is very similar and we only need to make a slight change. We will prove $\Delta^2 \geq \Delta_k(\delta_{low})$ in the following.

    Set $\delta_k=\delta_{low}$, similar to \autoref{lem.relation t and delta}, we have
    \begin{align*}
        -\theta_k & = v_k^T \Hk v_k+2 t_k \phi_k^T v_k +\delta_k t_k^2                                                         \\
                  & \geq \lambda_{1}(H_k) \|v_k\|^2 +2 t_k \phi_k^T v_k +\delta_k t_k^2                                        \\
                  & = \lambda_{1}(H_k) \|v_k\|^2 - 2 t_k^2 (\theta_k+\delta_k) +\delta_k t_k^2                                 \\
                  & = (\delta_k+\|\phi_k\|\frac{1-\Delta^2}{\Delta^2}) \|v_k\|^2 - 2 t_k^2 (\theta_k+\delta_k) +\delta_k t_k^2 \\
                  & = \delta_k -2 t_k^2 (\theta_k+\delta_k)+\|\phi_k\|\frac{1-\Delta^2}{\Delta^2}(1-t_k^2),
    \end{align*}
    which in turn leads to
    \begin{equation*}
        (2t_k^2-1)(\theta_k+\delta_k)\geq \|\phi_k\|\frac{1-\Delta^2}{\Delta^2}(1-t_k^2).
    \end{equation*}
    Also note that by \autoref{lemma.uppertheta}, we have
    \begin{equation*}
        (2t_k^2-1)\|\phi_k\| \geq \|\phi_k\|\frac{1-\Delta^2}{\Delta^2}(1-t_k^2),
    \end{equation*}
    rearranging items we can have $\Delta_{low}\leq \Delta^2$, which completes the proof.

    For dealing with $g_k \perp \mathcal{S}_1(H_k)$, we refer readers to Section 3 in \cite{rojas_new_2001} for a different treatment for hard case in trust-region subproblems.
\end{proof}

\subsection{Proof of \autoref{thm.recover.reg}}\label{proof.thm.recover.reg}
\begin{proof}
    To recover the gradient regularized Newton method via GHM, first, we set $\phi_k = g_k$ explicitly.
    Since the gradient regularized Newton method works for the convex objective functions, we need to select $\delta_k$ carefully in order to keep $\theta_k(\delta_k)$ strictly positive, by which we will also collect $t_k \neq 0$ from the ordering given by \autoref{lemma.ordering of tilde alpha 1}.
    By \autoref{lemma.delta upper bound} and Theorem 4.3 in \cite{gallier2020schur}, it follows that $\theta_k = 0$ is equivalent to
    \[
        H_k \succeq 0, \ (I-H_kH_k^{\star})g_k = 0, \ \delta_k \geq \overline{\delta_k^{\mathsf{cvx}}} = g_k^TH_k^{\star}g_k \geq 0.
    \]

    Therefore, at every iteration, we split into the following two cases:
    \begin{enumerate}[(a)]
        \item $(I-H_kH_k^{\star})g_k = 0$. In this case, we need to maintain $\delta_k < \overline{\delta_k^{\mathsf{cvx}}}$ and thus $F_k$ is indefinite.
              Also, we know $\theta_k > -\delta_k$ since $t_k \neq 0$. This justifies the following bracket for searching a proper $\delta_k$:
              $$[-2\gamma_k \|\gk\|^{1/2},  \overline{\delta_k^{\mathsf{cvx}}}), ~\forall \gamma_k > 0.$$
              Combining the continuity of $\theta_k$ (\autoref{lem.continuity.theteomega}), we have $\forall \gamma_k > 0$, $\exists \delta_k$ such that $\theta_k(\delta_k) = \gamma_k\|g_k\|^{1/2}$.
              We may use several ways for computation. Firstly, it can be found by a similar provable bisection procedure in $O(\log(\epsilon^{-1}))$ to the one in \autoref{sec.bisection procedure}.
              Secondly, note that since $t_k \neq 0$, $\omega_k = \theta_k^2$ exhibits a known derivative, cf. \eqref{eq.diff.omega}. This means Newton's method \cite{yeCombiningBinarySearch1994} may also be a possible choice in $O(\log\log(\epsilon^{-1}))$.

              By using either of the above methods, it implies that the HSODM can compute exactly the same steps as the gradient regularized Newton's method:
              $$
                  d_k=-\left(H_k+\theta_k \cdot I\right)^{-1} g_k=-\left(H_k+\gamma_k\left\|g_k\right\|^{1 / 2}\right)^{-1} g_k
              $$
              given $\theta_k(\delta_k) = \gamma_k\|g_k\|^{1/2}$.
        \item $(I-H_kH_k^{\star})g_k \neq 0$. Under this scenario, $F_k$ is always indefinite and hence $\theta_k > 0$ for all $\delta_k$; that is, we may need $\delta_k \to \infty$. However, a simple remedy would suffice. Since $\Hk \succeq 0$, we have:
              $$\tilde \Hk := \Hk + \frac{1}{2}\gamma_k\|\gk\|^{1/2} \succ 0.$$
              To this end, simply construct GHM by $\tilde \Hk$ and then find $\theta_k(\delta_k) = \frac{1}{2}\gamma_k\|g_k\|^{1/2}$. Since now $(I-\tilde \Hk \tilde \Hk^{\star})g_k = 0$ holds and the discussion follows from case (a).
    \end{enumerate}
    Combining two cases, the proof is complete.
\end{proof}
To realize a gradient regularized Newton step, it is not necessary to check whether $(I-H_kH_k^{\star})g_k = 0$ at each iteration. One can always apply the strategy in case (b). However, this finding is more of a theoretical result.

We are aware of other techniques in \cite{mishchenko_regularized_2023,doikov_gradient_2021,doikov_super-universal_2022} beyond simply choosing an appropriate $\gamma_k$.
As a result, designing an efficient way to realize the gradient regularized Newton method remains an interesting future work.

\section{Analysis of the Bisection Procedure}\label{sec.bisection procedure}

\subsection{Bisection using $h_k$}
We discuss the complexity of searching $\delta_k$ to satisfy $\sqrt{h_k} \in I_h$ in \autoref{alg.gadaptive}. Briefly, we make use of the analysis in \autoref{sec.proof.continuity.theteomega} where we justify the continuity and monotonicity of $h_k$. These results enable us to adopt a simple bisection procedure.
We first make the following assumption regarding the gradient and Hessian for the simplicity of the analysis of the bisection procedure

\begin{assumption}\label{assm.bounded g and h}
    Assume that the norm of the approximated gradient $\phi_k$ and the Hessian $\Hk$ are upper bounded along the iterates $\xk$, i.e., there exists two constant $ U_\phi >0$ and $U_H>0$ such that
    $
        \|\phi_k\| \leq U_\phi, \|\Hk\| \leq U_H.
    $
\end{assumption}
Then consider optimal solution $[v_k;t_k]$ in the following cases.

\begin{lemma}
    \label{lem.relation t and delta}
    In the GHM \eqref{eq.ghqm}, suppose $\phi_k \not \perp \mathcal{S}_1$, we have the following statement:
    \begin{enumerate}[(a)]
        \item if $\delta_k \leq \lambda_{1}(H_k)$, then $|t_k|\geq \frac{\sqrt{2}}{2}$.
        \item if $\delta_k \geq \lambda_{d}(H_k)$, then $|t_k|\leq \frac{\sqrt{2}}{2}$, where $\lambda_r(\cdot)$ denotes maximum eigenvalue.
    \end{enumerate}
\end{lemma}
\begin{proof}
    For the first statement, recall the definition of GHM \eqref{eq.ghqm}, we have
    \begin{align}
        \nonumber
        -\theta_k                   & = v_k^T \Hk v_k+2 t_k \phi_k^T v_k +\delta_k t_k^2                         \\
        \nonumber
                                    & \geq \lambda_{1}(H_k) \|v_k\|^2 +2 t_k \phi_k^T v_k +\delta_k t_k^2        \\
        \nonumber
                                    & = \lambda_{1}(H_k) \|v_k\|^2 - 2 t_k^2 (\theta_k+\delta_k) +\delta_k t_k^2 \\
        \label{eq.definition theta} & \geq \delta_k \|v_k\|^2 - 2 t_k^2 (\theta_k+\delta_k) +\delta_k t_k^2
        = \delta_k -2 t_k^2 (\theta_k+\delta_k),
    \end{align}
    where the second equality comes from the optimal condition \eqref{eq.homoeig.soc}, the following inequality is due to $\delta_k \leq \lambda_{1}(H_k)$, and the last equation holds since $[v_k;t_k]$ is a unit vector. Rearranging the terms, it follows that
    $
        (2t_k^2-1)(\theta_k+\delta_k) \geq 0,
    $
    and it further implies $|t_k |\geq \frac{\sqrt{2}}{2}$. By a similar argument, we can prove the second statement.
\end{proof}
Note that during the search procedure for $\delta$, we can always assume that $\phi_k \not \perp \mathcal{S}_1$; otherwise, we can apply a perturbation for $\phi_k$ to escape the case $\phi_k \perp \mathcal{S}_1$ using \autoref{alg.hardcase}. Based on the previous analysis, we know that the function $h_k(\delta)$ is continuous and decreasing with respect to $\delta$ \autoref{lem.differentiablity of h} as long as $\phi_k \not \perp \mathcal{S}_1$. Therefore, at each iteration in \autoref{alg.gadaptive} (Line 4), we can apply the bisection method to search for the proper $\delta_k$.
Without loss of generality, we assume that $I_h= [\sqrt{\ell},\sqrt{\nu}]$.
To apply the bisection method, we need to perform an effective bracketing first, which means finding an interval $I_k:=[\delta_{low}, \delta_{up}]$ that contains the value $\delta_k$ we are searching for. The following results give us an estimate by the fact that $h$ is continuous (see \autoref{lem.differentiablity of h}).
\begin{lemma}
    \label{lem.bracketing}
    Suppose $\phi_k\not \perp \mathcal{S}_1$, and for any given interval $I_h := [\sqrt{\ell},\sqrt{\nu}]$ such that the function value $h_k(\cdot) \in [\ell,\nu]$, if we set
    \small
    \begin{subequations}
        \begin{align}
            \delta_{low}               & = \min\left\{\lambda_{1}(\Hk),-\sqrt{\nu}\right\},                                                                                \\
            \label{eq.dub} \delta_{up} & = \max\left\{\frac{(1+\vert \lambda_{1}(H_k)\vert)^2(1+\lambda_{d}(H_k)+\vert \lambda_{1}(H_k)\vert)}{\ell},\|\phi_k\|^2\right\},
        \end{align}
    \end{subequations}
    \normalsize
    then it follows that
    \begin{equation*}
        [\ell,\nu]\subseteq [h_k(\delta_{up}),h_k(\delta_{low})].
    \end{equation*}
\end{lemma}
\begin{proof}
    To establish $[\ell,\nu]\subseteq [h_k(\delta_{up}),h_k(\delta_{low})]$, it is sufficient to show that $h_k(\delta_{low}) \geq \nu$ and $h_k(\delta_{up}) \leq \ell$.
    By \autoref{lem.relation t and delta}, when $\delta_{low} = \min\{\lambda_{1}(\Hk),-\sqrt{\nu}\}$ we have $|t_k|\geq \frac{\sqrt{2}}{2}$, thus $\frac{t_k^2}{1-t_k^2}\geq 1$. On the other hand, take $\delta_k = \delta_{low}$ and by optimality condition \eqref{eq.homoeig.soc}, we have $\theta_k+\delta_{low} \geq 0$, which leads to $\theta_k \geq \sqrt{\nu}$. Therefore, we conclude
    \begin{equation}
        \label{eq.bracketing left}
        h_k(\delta_{low}) = \frac{t_k^2}{1-t_k^2}\theta_k^2\geq \nu.
    \end{equation}
    Next, we prove $h_k(\delta_{up})\leq \ell$. Similarly, we need to bound two terms $\theta_k^2$ and $\frac{t_k^2}{1-t_k^2}$, respectively. For the term $\theta_k$, we claim that once $\delta_k\geq \|\phi_k\|^2$, $\theta_k$ satisfies $-\lambda_{1}(\Hk) \leq \theta_k \leq 1 + |\lambda_1(\Hk)|$. Recall that when $\phi_k \not \perp \mathcal{S}_1$, it is easy to see that strictly $\theta_k > - \lambda_{1}(H_k)$.

    Also, the value $\theta_k$ can be viewed as the solution of following equation (also see Theorem 3.1, \citet{rojas_new_2001})
    \begin{equation}
        \label{eq.solve theta_k}
        \delta + \theta = \sum_{i=1}^r \frac{\beta_i^2}{\lambda_i(\Hk) + \theta},
    \end{equation}
    where $\beta_i, i = 1,\cdots,r$ is defined in a way similar to the proof \autoref{sec.pf.lem.differentiablity of h}, cf. \eqref{eq.spec.decomp}.
    Now suppose $\theta_k = \vert \lambda_{1}(\Hk)\vert + 1$, in view of \eqref{eq.solve theta_k}, we must have,
    \begin{align}
        \nonumber \delta_k & = - \theta_k + \sum_{i=1}^r \frac{\beta_i^2}{\lambda_i(\Hk)+\theta_k}                                     \\
        \nonumber          & \le - 1 - \vert \lambda_{1}(\Hk)\vert+\frac{\|\phi_k\|^2}{\lambda_{1}(\Hk)+\vert \lambda_{1}(\Hk)\vert+1} \\
        \label{eq.def.l}   & \le \frac{\|\phi_k\|^2}{\lambda_{1}(\Hk)+\vert \lambda_{1}(\Hk)\vert+1} \le \|\phi_k\|^2
    \end{align}
    We conclude that as long as $\delta \ge \|\phi_k\|^2$, $\theta_k \le 1+\vert \lambda_{1}(\Hk)\vert$ since $\theta_k$ is decreasing in $\delta_k$.

    On the other hand, suppose that we set $\delta_k = \lambda_{d}(\Hk)+\alpha$, $\alpha>0$, we only have to give an estimate $\alpha$:
    \begin{align}
        \nonumber
        -\theta_k & = v_k^T \Hk v_k+2 t_k \phi_k^T v_k +\delta_k t_k^2                      \\
        \nonumber
                  & \leq (\delta_k-\alpha)\|v_k\|^2-2 t_k^2(\delta_k+\theta_k)+\delta t_k^2 \\
                  & = \delta_k-\alpha \|v_k\|^2 -2 t_k^2(\delta_k+\theta_k),
        \label{eq.estimate ub}
    \end{align}
    rearranging items we have
    \begin{align}
        \nonumber
        \frac{t_k^2}{1-t_k^2} & \leq 1-\frac{\alpha}{\theta_k+\delta_k}  = 1-\frac{\alpha}{\theta_k+\lambda_{d}(\Hk)+\alpha}                      \\
        \nonumber
                              & \leq 1-\frac{\alpha}{1+\vert \lambda_{1}(\Hk)\vert +\lambda_{d}(\Hk)+\alpha}                                      \\
                              & = \frac{1+\vert \lambda_{1}(\Hk)\vert +\lambda_{d}(\Hk)}{1+\vert \lambda_{1}(\Hk)\vert +\lambda_{d}(\Hk)+\alpha}.
        \label{eq.analysis for t}
    \end{align}
    We note that LHS of \eqref{eq.analysis for t} is decreasing in $\alpha$, and thus if we let
    $$\alpha \geq \alpha_u:= -1 - \vert \lambda_{1}(\Hk)\vert-\lambda_{d}(\Hk) + \frac{(1+\vert \lambda_{1}(\Hk)\vert)^2(1+\lambda_{d}(\Hk)+\vert \lambda_{1}(\Hk)\vert)}{\ell}$$
    we have
    \begin{equation}\label{eq.ubt}
        \frac{t_k^2}{1-t_k^2} \leq \frac{\ell}{(1+\vert \lambda_{1}(\Hk)\vert)^2}.
    \end{equation}
    By construction \eqref{eq.dub} we see that $\delta_{up} \geq \lambda_{d}(\Hk)+\alpha_u$, combining \eqref{eq.def.l}, we have
    \begin{align*}
        h_k(\delta_{up}) & \le \theta_k^2 \cdot \frac{\ell}{(1+\vert \lambda_{1}(\Hk)\vert)^2}                              \\
                         & \le \left(|\lambda_1(\Hk)| + 1\right)^2 \frac{\ell}{(1+\vert \lambda_{1}(\Hk)\vert)^2} \leq \ell
    \end{align*}
    as desired. This completes the proof.
\end{proof}

In the above lemma, we demonstrate the existence of interval for $\delta_k$ such that $h_k(\delta_k) \in [\ell, \nu]$ by providing an explicit $I_k$ for any $[\ell, \nu]$. To examine the complexity of the bisection procedure, we define the following target interval for $\delta_k$ allowing a tolerance $\sigma>0$,
\begin{equation}
    \label{eq.target interval}
    I_k^{\sigma} = [\underline \delta,\overline{\delta}],\quad h_k(\underline \delta) = \nu+\sigma, \quad h_k(\overline{\delta}) = \ell
\end{equation}
Note that if $\sigma>0$ is small enough, the inexactness will not affect the convergence rate of the framework of \autoref{alg.gadaptive}. In the sequel, we focus on estimating the ratio between the length of $I_k^{\sigma}$ and the length of $I_k$, which is essential in the analysis of the complexity of the bisection procedure. Next, we give a lower bound for the length of the interval $I_k^{\sigma}$.

\begin{lemma}
    \label{lem.lb interval}
    In the $k$-th iteration, suppose $\dk$ is the solution of \eqref{eq.ghqm} with $\delta_k\in I_k^\sigma$, then the length of the interval $I_k^\sigma$ is at least
    \begin{equation}
        \label{eq.length interval}
        \frac{\sigma \|\dk\|}{2\sqrt{\nu+\sigma}+(\nu+\sigma)\|\phi_k\|}.
    \end{equation}
\end{lemma}
\begin{proof}
    By the definition of $I_k^\sigma$ in \eqref{eq.target interval}, we have $\ell\leq h_k(\delta)= \frac{\theta^2}{\|d\|^2}\leq \nu+\sigma,$
    thus
    \begin{equation*}
        \ell\|d\|^2\leq \theta^2 \leq (\nu+\sigma)\|d\|^2.
    \end{equation*}
    Plugging the above into \eqref{eq.derivative of h}, we have
    \begin{equation}
        \label{eq.upper bound abs h auxi}
        \vert h_k'(\delta)\vert \leq \frac{2\sqrt{\nu+\sigma}\|d\|^3+(\nu+\sigma)\|d\|^2d^T(\Hk+\theta I)^{-1}d}{\|d\|^4}.
    \end{equation}
    On the other hand, from the optimality condition \eqref{eq.homoeig.soc} and its Schur complement we have
    \begin{equation*}
        \Hk+\theta I - \frac{\phi_k \phi_k^T}{\delta+\theta} \succeq 0,
    \end{equation*}
    as a result,
    \begin{align*}
        \Hk+\theta I \succeq \frac{\phi_k \phi_k^T}{\delta+\theta}
        = -\frac{\phi_k \phi_k^T}{\phi_k^T d}
        \succeq \frac{\phi_k \phi_k^T}{\|\phi_k\|\|d\|},
    \end{align*}
    which further implies
    \begin{equation}
        \label{eq.lower bound h+theta}
        d^T(\Hk+\theta I)^{-1} d\leq \frac{(\phi_k^T d)^2}{\|\phi_k\|\|d\|}.
    \end{equation}
    Plugging \eqref{eq.lower bound h+theta} into \eqref{eq.upper bound abs h auxi}, we have
    \begin{equation}
        \label{eq.upper bound abs h}
        \vert h_k'(\delta)\vert \leq \frac{2\sqrt{\nu+\sigma}+(\nu+\sigma)\|\phi_k\|}{\|d\|}.
    \end{equation}
    By the mean value theorem, it follows that
    \begin{equation*}
        h_k(\underline \delta)- h_k(\overline{\delta}) = h_k'(\xi)(\overline{\delta}-\underline \delta),\quad \xi \in [\underline \delta,\overline{\delta}].
    \end{equation*}
    Together \eqref{eq.target interval} with \eqref{eq.upper bound abs h}, we conclude
    \begin{equation*}
        \overline{\delta}-\underline \delta \geq \frac{\sigma \|\dk\|}{2\sqrt{\nu+\sigma}+(\nu+\sigma)\|\phi_k\|},
    \end{equation*}
    which finishes the proof.
\end{proof}

We remark that during \autoref{alg.gadaptive} we can assume that $\|\dk\| \geq \sqrt{\epsilon}$. Otherwise, by \eqref{eq.relation step and g} and \eqref{eq.secondorderoptc}, we see
\begin{equation*}
    \|\gkn\| \leq O(\epsilon), \quad \lambda_{\min}(\Hk) \geq \Omega(-\sqrt{\epsilon}).
\end{equation*}
By the Lipschitz continuity of $\nabla ^2 f(x)$, we know that $\lambda_{\min}(H_{k+1}) \geq \Omega(-\sqrt{\epsilon})$ and thus $\xkn:= \xk +\dk$ satisfies \eqref{eq.approxfocp} and \eqref{eq.approxsocp} and we can terminate the algorithm at the iterate $x_{k+1}$.

To inspect the ratio between $|I_k|$ and $|I_k^\sigma|$, by \autoref{assm.bounded g and h}, we establish a rough estimate of the length of the interval $I_k^\sigma$:
\begin{align}\label{eq.rough length of I}
    \overline{\delta}-\underline \delta = \Omega\left(\frac{\sigma\sqrt{\epsilon}}{\varsigma_h U_\phi}\right),
\end{align}
which is by \autoref{lem.upper bound h} and the fact that $\nu$ is determined by function value $h_{k-1}(\delta_{k-1})$.

Now we are ready to give the final analysis of the complexity of the bisection procedure.

\subsubsection{Proof of \autoref{thm.complexity bisection procedure}}
\begin{proof}
    Note that we are using the bisection method to search for $I_k^\sigma$ within the initial interval $I_k$. By the mechanism of the bisection method, the complexity of bisection is
    \begin{equation}
        \label{eq.rough complexity}
        O\left(\log \frac{|I_k|}{|I_k^\sigma|}\right)= O\left(\log \frac{\delta_{up}-\delta_{low}}{\overline{\delta}-\underline \delta}\right).
    \end{equation}
    By the inequality \eqref{eq.rough length of I}, we already prove a lower bound for $\overline{\delta}-\underline \delta$. It remains to estimate $\delta_{up}-\delta_{low}$, recall our setting in \autoref{lem.bracketing}, it implies that
    \begin{align}
        \nonumber
        \delta_{up}-\delta_{low} & \leq \max\left\{(U_\phi)^2,\frac{(1+U_H)^2(1+2U_H)}{h_{\min}}\right\}+\max\{U_H,\sqrt{\varsigma_h}\} \\
                                 & \leq O\left(\frac{(U_\phi)^2(U_H)^3\sqrt{\varsigma_h}}{h_{\min}}\right),
        \label{eq.estimate Ib}
    \end{align}
    plugging \eqref{eq.estimate Ib} and \eqref{eq.rough length of I} into \eqref{eq.rough complexity} and we get the desired result.
\end{proof}

\section{Basic Properties of Concordant Second-Order Lipschitz Functions}\label{sec.concordlip}
We first provide some basic properties.
\begin{proposition}
    \label{cor.concordant convex}
    If the function $f$ satisfies the \betacon{}, then $f$ is convex.
\end{proposition}
\begin{proof}
    The inequality \eqref{eq.concordant inequality} indicates that $\nabla^2 f(x) \succeq 0$ for all $x \in \operatorname{dom}(f)$. Hence, $f$ is convex.
\end{proof}
\begin{remark}
    The above corollary implies that the \betacon{} function is \textit{not strictly or strongly convex}, which means we allow the degeneracy of Hessian.
\end{remark}

Now we present more basic properties of \betacon{}, enabling us to derive more ``complicated'' examples. Notably, we show that the space of \betacon{} functions is closed under positive scalar multiplications and summations. This property is also preserved in the affine transformation of variables. For compactness, the proof of the following Lemmas and Examples are left in the appendix.
\begin{lemma}[Sum of \betacon{} functions]\label{prop.sum}
    Let $f_i$ be $\beta_i$-\betacon{} functions satisfying \eqref{eq.concordant inequality}, where $\beta_i \geq 0$ for $i = 1, \cdots, m$. Then $\sum_{i=1}^m f_i$ is a \betacon{} function.
\end{lemma}
\begin{proof}
    By \autoref{def.betacon}, for any point $x \in \bigcap_{i=1}^m \operatorname{dom}(f_i)$, it follows that
    $$
        \begin{aligned}
             & \left\|\sum_{i=1}^m \nabla f_i(x+d) - \sum_{i=1}^m \nabla f_i(x) - ( \sum_{i=1}^m \nabla^2f_i(x)) d\right\| \\
             & \leq \sum_{i=1}^m \left\|\nabla f_i(x+d) - \nabla f_i(x) - \nabla^2 f_i(x)d\right\|                         \\
             & \leq \max_{1 \leq i\leq m} \{\beta_i\} \cdot d^T \sum_{i=1}^m \nabla^2 f_i(x) d,
        \end{aligned}
    $$
    whenever we choose $\|d\| \leq \min_{1 \leq i\leq m} C_i$.
\end{proof}

\begin{lemma}\label{prop.constant coefficient}
    Suppose the function $f$ satisfies the $\beta$-\betacon{} condition, then for any coefficient $c > 0$, the function $c \cdot f$ is $c\beta$-\betacon{}.
\end{lemma}
\begin{proof}
    The definition of \betacon{} functions directly implies the result.
\end{proof}

\begin{lemma}[Composite function]\label{prop.composite}
    The composite function $f(x) = \phi(Ax-b)$ is \betacon{} if $\phi(\cdot)$ is $\beta$-\betacon{}.
\end{lemma}
\begin{proof}
    Note that $\nabla f(x) = A^T \nabla \phi(Ax-b)$ and $\nabla^2 f(x) = A^T \nabla^2 \phi(Ax-b) A$, then for any point $x$ such that $Ax-b \in \operatorname{dom}(\phi)$, we have
    $$
        \begin{aligned}
             & \|\nabla f(x+d) - \nabla f(x) - \nabla^2 f(x)d\|                                 \\
             & = \|A^T (\phi(Ax-b+Ad) - \nabla \phi(Ax-b) - \nabla^2 \phi(Ax-b)Ad) \|           \\
             & \leq \|A^T\| \cdot \|\phi(Ax-b+Ad) - \nabla \phi(Ax-b) - \nabla^2 \phi(Ax-b)Ad\| \\
             & \leq \|A^T\| \cdot \beta \cdot (Ad)^T \nabla^2 \phi(Ax-b) (Ad)                   \\
             & = \|A^T\|\beta \cdot d^T \nabla^2 f(x) d.
        \end{aligned}
    $$
\end{proof}

Moreover, we provide a sufficient condition for ensuring the \betacon{}.
\begin{lemma}\label{prop.sufficient condition}
    Let function $\phi(y)$, $y \in \real^m$ be standard $M$-second-order Lipschitz continuous and $\mu$-strongly convex, then the funtion
    $
        f(x)=\phi(Ax-b)
    $
    is \betacon{} for all $x \in \real^n$ such that $y = Ax - b$ is in the domain of $\phi$, where $A \in \real^{m \times n}$, $m \leq n$ is a constant coefficient matrix with rank $m$.
\end{lemma}
\begin{proof}
    By a similar argument in \autoref{prop.composite}, it follows
    $$
        \begin{aligned}
             & \|\nabla f(x+d) - \nabla f(x) - \nabla^2 f(x)d\|                                                                                                                    \\
             & = \|A^T (\nabla \phi(Ax-b+Ad) - \nabla \phi(Ax-b) - \nabla^2 \phi(Ax-b)Ad)\|                                                                                        \\
             & \leq \|A^T\| \cdot \|\nabla \phi(Ax-b+Ad) - \nabla \phi(Ax-b) - \nabla^2 \phi(Ax-b)Ad\|                                                                             \\
             & \leq \|A^T\| \cdot \frac{M}{2}\|Ad\|^2 \leq \|A^T\| \cdot \frac{M}{2\mu} \cdot d^T (A^T \nabla^2 \phi(Ax-b) A) d = \frac{\|A^T\|M}{2\mu} \cdot d^T \nabla^2 f(x) d.
        \end{aligned}
    $$
    Set $\beta = \frac{\|A^T\|M}{2\mu}$ and we complete the proof.
\end{proof}

\begin{corollary}\label{cor.sufficient condition}
    If the function $f$ is standard second-order Lipschitz continuous and $\mu$-strongly convex, then it is \betacon{}.
\end{corollary}

The above fact easily validates the logistic functions.
\begin{example}[Logistic regression with $L_2$ penalty]\label{example.l2 logistic}
    Consider the following logistic regression function with $l_2$ penalty,
    \begin{equation}\label{eq.logistic regression}
        f(x) = \frac{1}{m} \sum_{i=1}^m \log \left(1 + e^{-b_i \cdot a_i^T x}\right) + \frac{\gamma}{2}\|x\|^2,
    \end{equation}
    where $\gamma > \frac{2}{m}\sum_{i=1}^m \|a_i\|^2$, $a_i \in \real^n$, $b_i \in \{-1, 1\}$, $i = 1, 2, \cdots, m$. Then the function $f(x)$ satisfies the \betacon{} condition.
\end{example}
\begin{proof}
    We first let
    $$\nu = \lambda_{\max}\left(\frac{1}{m}\sum_{i=1}^m b_i^2 \cdot a_i a_i^T\right)$$
    be the maximum eigenvalue.
    Then define the univariate function $\phi(y) = \log\left(1+e^{-y}\right) + \frac{\gamma}{2\nu}\cdot y^2$, $y \in \real$, then $\phi(y)$ is standard second-order Lipschitz continuous and $\frac{\gamma}{\nu}$-strongly convex.
    By \autoref{prop.sufficient condition}, it implies that
    $$
        g_i(x) = \phi_i(b_i\cdot a_i^T x) = \log\left(1+e^{-b_i \cdot a_i^T x}\right) + \gamma\cdot \frac{(b_i \cdot a_i^T x)^2}{2\nu},
    $$
    is \betacon{} for all $i = 1, \dots, m$.
    It remains to see that
    $$
        h(x) = x^T \left(\frac{\gamma}{2}I - \frac{\gamma}{2\nu}\frac{1}{m}\sum_{i=1}^m b_i^2 \cdot a_i a_i^T\right) x
    $$
    is a convex quadratic function and thus also \betacon{}.
    Using the additive rule, we see that the function $f$ is the summation of $g_i(x)$ and $h(x)$, i.e.
    $$
        f(x) = \frac{1}{m} \sum_{i=1}^m g_i(x) + h(x).
    $$
    Combing the \autoref{prop.sum} and \autoref{prop.constant coefficient}, we conclude that $f(x)$ satisfies the \betacon{} condition.
\end{proof}

\section{Extra Proofs in \autoref{sec.homotopy}}
\subsection{Proof of \autoref{lemma. homotopy property}}\label{sec.lem.homotopy.property}
\begin{proof}
    Note that the regulated objective function $f(x) + \frac{\mu}{2}\|x\|^2$ is strongly convex for any given $\mu > 0$ so that its minimizer $\xmu$ is unique, and it is easy to see that $\xmu = -\nabla f(\xmu) / \mu$ is a continuous function due to the continuity of $\nabla f(\xmu)$ and $\frac{1}{\mu}$. Therefore, the first statement holds. As for the second statement, for any $0 < \mu' < \mu$, it follows that
    \begin{equation}\label{eq.mu mu'}
        f(x_{\mu'}) + \frac{\mu'}{2}\|x_{\mu'}\|^2 < f(x_\mu) + \frac{\mu'}{2}\|x_\mu\|^2
    \end{equation}
    and
    \begin{equation}\label{eq.mu' mu}
        f(x_{\mu}) + \frac{\mu}{2}\|x_{\mu}\|^2 < f(x_{\mu'}) + \frac{\mu}{2}\|x_{\mu'}\|^2.
    \end{equation}
    Add the two inequalities \eqref{eq.mu mu'} \eqref{eq.mu' mu} on both sides and rearrange terms, we have
    $$
        \frac{\mu - \mu'}{2}\|x_{\mu'}\|^2 > \frac{\mu - \mu'}{2}\|x_{\mu}\|^2.
    $$
    It implies that $\|x_{\mu'}\| > \|x_\mu\|$ since $\mu - \mu' > 0$, and hence $\|x_\mu\|$ is strictly decreasing function of $\mu$. Substitute $\|x_{\mu'}\| > \|x_\mu\|$ into the inequality \eqref{eq.mu mu'}, we obtain that
    $
        f(x_{\mu'}) < f(x_\mu),
    $
    which completes the proof. Now, we prove the third statement. Recall that $x^*$ is the minimal $L_2$ norm solution of $f(x)$, then $\nabla f(x^*) = 0$, together with $\xmu = \arg\min \ \left\{f(x) + \frac{\mu}{2}\|x\|^2\right\}$, we have
    $$
        \nabla f(\xmu) - \nabla f(x^*) + \mu \xmu  = 0.
    $$
    Multiplying $\xmu - x^*$ on both sides and by the convexity of $f$ (\autoref{cor.concordant convex}), it follows
    $$
        -\mu(\xmu - x^*)^T \xmu = (\xmu - x^*)^T(\nabla f(\xmu) - \nabla f(x^*)) \geq 0,
    $$
    and it further implies that
    $
        \|\xmu\|^2 \leq \xmu^T x^* \leq \|x^*\|\|\xmu\|,
    $
    that is $\|\xmu\| \leq \|x^*\|$ for any $\mu > 0$. By the uniqueness of $x^*$, we conclude that $\lim_{\mu \rightarrow 0^+} \xmu = x^*$. As for the final statement, suppose that $\xmu \rightarrow z \in \mathbb{R}^n \neq 0$ as $\mu \rightarrow \infty$, then we could choose $\mu > \frac{2\left(f(0) - f(z)\right)}{\|x^*\|^2 - \|z\|^2}$, and it implies that
    $$
        f(z) + \frac{\mu}{2}\|z\|^2 > f(0).
    $$
    The above inequality contradicts that $\xmu = \arg\min \ \left\{f(x) + \frac{\mu}{2}\|x\|^2 \right\}$.
\end{proof}

\end{document}